\documentclass[reqno,a4paper]{amsart}

\usepackage{mathrsfs,amsmath,amssymb,graphicx}
\usepackage[unicode]{hyperref}
\usepackage{latexsym}
\usepackage{layout}
\usepackage[english]{babel}
\usepackage{color}
\usepackage{fancyvrb}
\usepackage{amsmath,amscd}
\usepackage{txfonts}
\usepackage{enumitem} 
\usepackage{float}
\usepackage{multicol} 

\usepackage{subcaption}

\usepackage{transparent}

\usepackage[normalem]{ulem} 

\usepackage{tikz}
\usetikzlibrary{calc}
\usetikzlibrary{matrix,arrows,decorations,chains,positioning}

\makeatletter
\@namedef{subjclassname@2020}{\textup{2020} Mathematics Subject Classification}
\makeatother

\theoremstyle{plain}
\newtheorem*{theorem*}{Theorem} 
\newtheorem{theorem}{Theorem}[section]
\newtheorem{lemma}[theorem]{Lemma}
\newtheorem{proposition}[theorem]{Proposition}
\newtheorem{corollary}[theorem]{Corollary}

\theoremstyle{definition}
\newtheorem{definition}{Definition}[section]
\newtheorem{remark}[theorem]{Remark}
\newtheorem{question}{Question}[section]
\newtheorem{problem}[question]{Problem}

\newcommand{\tildeHK}{\widetilde{\rm HK}}
\newcommand{\HK}{{\rm HK}}
\newcommand{\HKA}{{\rm HK_A}}
\newcommand{\LHK}{{\rm HK^\circ_\Gamma}}
\newcommand{\DLHK}{{\rm HK^\circledcirc_\Gamma}}

\newcommand{\SG}{{\Gamma}}
\newcommand{\LSG}{{\Gamma^\circ}}
\newcommand{\DLSG}{{\Gamma^\circledcirc}}
\newcommand{\LA}{{A^\circ_\Gamma}}

\newcommand{\Sbb} {\mathbb S}
\newcommand{\sphere} {\Sbb^3}
\newcommand{\disk}{{\mathcal D}}

\newcommand{\id}{\rm id}
 
\newcommand{\bCompl}[1]{\partial E(#1)}
\newcommand{\Compl}[1]{E(#1)}
\newcommand{\ComplHKA}{\Compl{\HK_A}}
\newcommand{\bComplHKA}{\bCompl{\HK_A}}
\newcommand{\ComplHK}{\Compl \HK}
\newcommand{\bComplHK}{\bCompl\HK}

\newcommand{\snode}{\bullet}
\newcommand{\hnode}{\circ}

\newcommand{\bp}[1]{\smash{\uuline{#1}}}
\newcommand{\cbp}[1]{\smash{\overline{\uuline{#1}}}}

\newcommand{\pbp}[1]{\tilde{\bp #1}}
\newcommand{\bphi}{\cbp \phi}
\newcommand{\elem}[1]{\mathtt{#1}}
\newcommand{\fr}[2]{\partial_{#2}{#1}}

\newcommand{\F}{\mathcal{F}}
\newcommand{\Q}{\mathrm{Q}}

\newcommand{\Li}{\mathrm L}
\newcommand{\Hex}{\mathcal{H}}
\newcommand{\bm}{\partial_M}
\newcommand{\bEhk}{\partial_{\ComplHK}}
\newcommand{\pairtilde}{(\sphere,\tildeHK)}
\newcommand{\pair}{(\sphere,\HK)}
\newcommand{\pairn}{(\sphere,\HK_n)}
\newcommand{\pairA}{(\sphere,\HK_A)}

\newcommand{\pairLHK}{(\sphere,\LHK)}
\newcommand{\pairDLHK}{(\sphere,\DLHK)}

\newcommand{\pairSG}{(\sphere,\SG)}

\newcommand{\pairLSG}{(\sphere,\LSG)}
\newcommand{\pairDLSG}{(\sphere,\DLSG)}
\newcommand{\pairNSG}{(\sphere,\rnbhd{\SG})}

\newcommand{\pairK}{(\sphere,K)}
\newcommand{\pairKtau}{(\sphere,K\cup\tau)}

\newcommand{\pairL}{(\sphere,L)}
\newcommand{\pairLn}{(\sphere,L_n)}

\newcommand{\pairLntau}{(\sphere,L_n\cup\tau)}

\newcommand{\pairfourone}{(\sphere,4_1)}

\newcommand{\lk}[2]{{\ell \mathit k}(#1,#2)}

\newcommand{\hopf}{{\bf h}}
\newcommand{\knot}{{\bf k}}
\newcommand{\link}{{\bf l}}
\newcommand{\EM}{{\bf em}}
\newcommand{\charM}{{\Lambda_{\textsc{M}}}}
\newcommand{\charE}{{\Lambda_{\Compl \HK}}}
\newcommand{\charhk}{{\Lambda_{\textsc{hk}}}}

\newcommand{\Z}{\mathbb{Z}}

\newcommand{\cout}[1]   {}

\newcommand{\op}[1]{\operatorname{#1}}
\newcommand{\til}[1]{{\tilde{#1}}}
\newcommand{\cut}{\op{cut}}

\newcommand{\rnbhd}[1]{\mathfrak N(#1)} 

\newcommand{\openrnbhd}[1]{\mathring{\mathfrak N}(#1)}

\newcommand{\Sym}[2][\sphere]{\mathcal MCG(#1, #2)}
\newcommand{\pSym}[2][\sphere]{\mathcal MCG_+(#1, #2)}
\newcommand{\gSym}[2][\sphere]{\mathcal MCG_{(+)}(#1, #2)}
\newcommand{\MCG}[1]{\mathcal MCG(#1)}
\newcommand{\pMCG}[1]{\mathcal MCG_+(#1)}
\newcommand{\Aut}[1]{\mathcal Homeo(#1)}
\newcommand{\pAut}[1]{\mathcal Homeo_+(#1)}

\newcommand{\rel}{{\rm rel\,}} 
\newcommand{\rot}{\mathtt{r}}
\newcommand{\mir}{\mathtt{m}}
\newcommand{\conj}{\mathtt{c}}

\definecolor{mygray}{rgb}{0.92,0.92,0.92}

\numberwithin{equation}{section}
\numberwithin{figure}{section}

\title{JSJ decomposition for handlebody-knots}
 
\author{Yi-Sheng Wang}
\address{National Sun Yat-sen University, Kaohsiung 804, Taiwan}
\email{yisheng@mail.nsysu.edu.tw}

\date{\today}

\begin{document}
 

\begin{abstract} 
The paper applies the JSJ decomposition
and Koda-Ozawa's annulus classification to 
analyze the annulus configuration in a 
handlebody-knot exterior. 
We introduce the notion of the annulus diagram,
to pack the configuration into a labeled graph, 
and classify genus two handlebody-knots in
terms of their annulus diagrams.
Applications to handlebody-knot symmetries  
are discussed; methods to produce 
handlebody-knots with various types of annulus diagrams are also presented. 
\end{abstract}

\maketitle
 
\section{Introduction}\label{sec:intro}

Let $M$ be an oriented, irreducible, $\partial$-irreducible $3$-manifold. The JSJ decomposition asserts that,
up to isotopy, there is a unique 
surface $S\subset M$ consisting of 
essential annuli and tori such that
\textbf{1.} every component of the exterior $\Compl S:=M-\openrnbhd{S}$ 
is either I-/Seifert fibered or hyperbolic 
and \textbf{2.} the removal of any component of $S$
causes the first condition to fail, where
$\openrnbhd{S}$ is an open regular neighborhood of 
$S\subset M$ \cite{JacSha:79}, \cite{Joh:79} (see also \cite{Bon:01}).
Assign a solid (resp.\ hollow) 
node to each fibered (resp.\ hyperbolic) component of
$\Compl S$, and to each component $N$ 
of $\openrnbhd{S}$ 
assign an edge between nodes corresponding to component(s) of $\Compl S$ that meets/meet the frontier of $N$.
The resulting graph is called 
a \emph{characteristic diagram} $\charM$ of $M$.

The present work concerns the case 
where $M$ is atoroidal, namely, containing no essential tori, and embeddable in 
an oriented $3$-sphere $\sphere$.
Note that by Fox \cite{Fox:48}, $M$ is homeomorphic to 
a handlebody-link exterior---the exterior
of some tangled handlebodies in $\sphere$.
Atoroidality and essentiality of $M$ 
impose strong topological constraints
on its JSJ decomposition. 
If $g(\partial M)=1$, 
there is only one way to embed $M$ in $\sphere$ by Gordon-Luecke \cite{GorLue:89} and its exterior in $\sphere$ is always a solid torus. The characteristic diagram $\charM$ in this case is either Figs.\ \ref{fig:000h} or \ref{fig:100s}. In the former, $M$ is a hyperbolic knot exterior, whereas in the latter $M$ is a torus knot exterior. 
Here we classify the characteristic diagram $\charM$ for $M$ with $g(\partial M)=2$. 
\begin{theorem}[Theorem \ref{teo:possible_shapes}]\label{intro:teo:char_diagram}
Given $M$ with $g(\partial M)=2$,
its characteristic diagram $\charM$ is one of the entries in Table \ref{tab:configurations}.
\end{theorem}
By Thurston's hyperbolization theorem \cite{Thu:82}, 
$M$ is either hyperbolic or cylindrical, namely, $M$ containing an essential annulus; it is the former if and only if $\charM$ is Fig.\ \ref{fig:000h}. 
It is an interesting question as to whether all diagrams in
Table \ref{tab:configurations} can be realized by
such an $M$. 
To the author's knowledge, there is 
currently no known example 
whose characteristic diagram
is Figs.\ \ref{fig:200s}, \ref{fig:301h}, \ref{fig:301s} or \ref{fig:300h}.

Recall that the $W$-system of $M$ introduced by
Neumann-Swarup \cite{NeuSwa:97} is a maximal set of \emph{canonical} annuli in $M$, where an essential annulus is \emph{canonical} if 
any other essential annulus can be isotoped away from it. 
Theorem \ref{intro:teo:char_diagram}, together with 
Theorem \ref{teo:char_three_bigons} and Proposition \ref{prop:properties_char_diagram_Ehk}\ref{itm:I_bundles}, implies the following.
 
\begin{corollary}
A $W$-system of $M$ with $g(\partial M)=2$ coincide with the JSJ decomposition if $\charM$ is not one of Figs.\ \ref{fig:201h}, \ref{fig:301h} and \ref{fig:301s}.
\end{corollary}
 
\begin{corollary}
Given $M$ with $g(\partial M)=2$, 
then up to isotopy, 
$M$ contains four (resp.\ five, and infinitely many) essential annuli if 
$\charM$ is Figs.\ \ref{fig:200s} or \ref{fig:201h} (resp.\ \ref{fig:301h} or \ref{fig:301s}, and \ref{fig:100s}); otherwise, 
$M$ contains 
at most three essential annuli.   
\end{corollary}
%
%
%
In Sections.\ \ref{sec:classification}-\ref{sec:symmetry}, we apply Theorem \ref{intro:teo:char_diagram} to study 
\emph{handlebody-knots of genus $2$}, abbreviated to \emph{handlebody-knots} unless otherwise specified; a genus $g$
handlebody knot $\pair$ is a genus $g$ handlebody $\HK$ in $\sphere$. While, up to isotopy, a genus $1$ handlebody-knot, 
equivalent to a classical knot, 
is determined by its exterior by Gordon-Luecke \cite{GorLue:89}, 
there are infinitely many \emph{inequivalent}, namely non-isotopic, genus $2$ handlebody-knots with 
homeomorphic exteirors by Motto \cite{Mott:90}, Lee-Lee \cite{LeeLee:12}. In particular, 
the characteristic diagram $\charE$ of the 
handlebody-knot exterior
$\Compl\HK$ cannot differentiate them, 
and finer information has to be added.

The present work concerns non-trivial \emph{atoroidal}
handlebody-knots $\pair$---that is, $\Compl\HK$ is atoroidal and not a handlebody; 
they are of particular interest, being precisely those with a finite symmetry group by Funayoshi-Koda \cite{FunKod:20}, where the (positive) symmetry group $\gSym\HK$ of $\pair$, as defined in Koda \cite{Kod:15}, is the (positive) mapping class group of the pair $\pair$.

\begin{figure}[b]
\begin{subfigure}{0.45\linewidth}
\centering
\includegraphics[scale=.15]{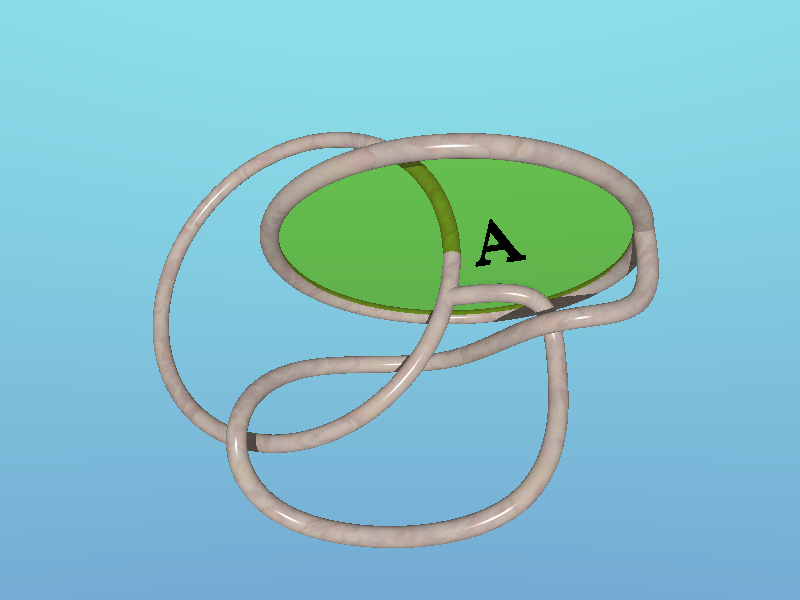}
\caption{$\charhk$ of $(\sphere,5_1)$: 
\raisebox{-.4\height}{\includegraphics[scale=.1]{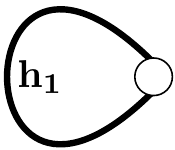}}.
}
\label{intro:fig:ann_fiveone} 
\end{subfigure}
\begin{subfigure}{0.45\linewidth}
\centering
\includegraphics[scale=.15]{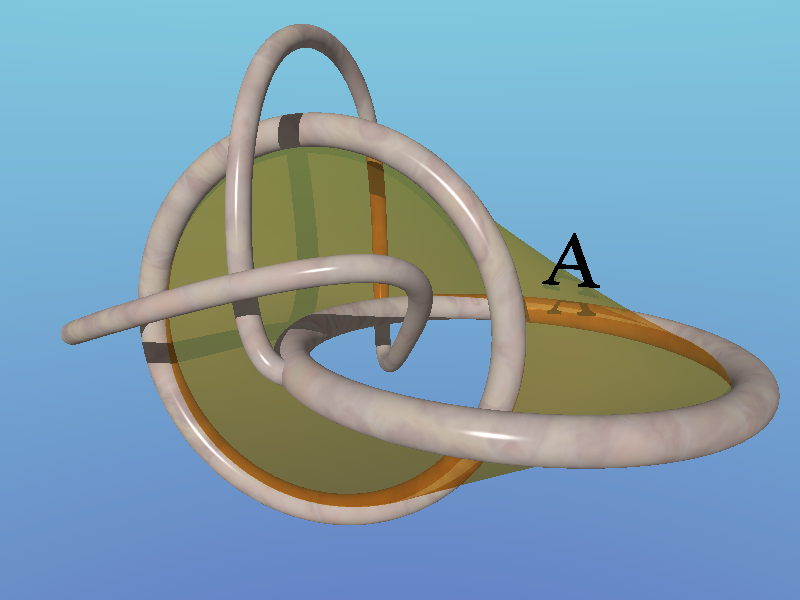} 
\caption{$\charhk$ of $(\sphere,6_4)$: 
\raisebox{-.4\height}{\includegraphics[scale=.1]{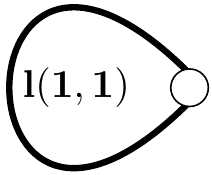}}.
}
\label{intro:fig:ann_sixfour}
\end{subfigure}
\caption{Annulus diagrams.}
\end{figure}

To enhance the characteristic diagram $\charE$, 
we recall that Koda-Ozawa 
\cite{KodOzaGor:15} and Funayoshi-Koda \cite[Lemma $3.2$]{FunKod:20} show that only four types of 
annuli $A$ can occur as essential annuli in an atoroidal 
handlebody-knot exterior $\ComplHK$.
These four types can be described in terms of $\partial A$
in relation to the handlebody $\HK$ \cite[Proof of Theorem $3.3$]{KodOzaGor:15}.
\begin{itemize}
\item[\emph{Type $2$}]: Exactly one component $l_A$ of $\partial A$ bounds a disk $\disk_A$ in $\HK$;  
if the disk $\disk_A$ is non-separating (resp.\ separating) in $\HK$, 
then $A$ is of \emph{type $2$-$1$} (resp.\ \emph{type $2$-$2$}). For an example of a type $2$-$1$ annulus, see Fig.\  \ref{intro:fig:ann_fiveone}.
\item[--] The symbol $\hopf_i$ is reserved for a type $2$-$i$ annulus, $i=1,2$.  
\item[\emph{Type $3$-$2$}]: Components of 
$\partial A$ are \emph{parallel} in $\partial \HK$ 
and bound no disks in $\HK$,
and there exists a unique \emph{non-separating} disk $\disk_A\subset\HK$ disjoint from $\partial A$ \cite{FunKod:20}. Let $V:=\HK-\openrnbhd{D}$. Then
$A$ is of type $3$-$2$i (resp.\ type $3$-$2$ii) if $A$ is essential (resp.\ inessential) in $\Compl V$. 
\item[--] The symbol $\knot_\ast$ is reserved for a type $3$-$2\ast$ annulus.

\item[\emph{Type $3$-$3$}]: Components of 
$\partial A\subset \partial\HK$ are \emph{non-parallel} and bound no disks in $\HK$; there exists a unique \emph{separating} essential disk $\disk_A$ in $\HK$ disjoint from $\partial A$ \cite{Wan:22}.
The disk $\disk_A$ cuts $\HK$ into two solid tori, 
each containing a component of $\partial A$. 
\emph{The slope pair} of $A$ is the slopes of $\partial A$
with respect to the two solid tori. For instance, the handlebody-knot in Fig.\ \ref{intro:fig:ann_sixfour} 
admits a type $3$-$3$ annulus with a slope pair
$(1,1)$.  
\item[--] The symbol $\link(r_1,r_2)$ denotes a type $3$-$3$ annulus with a slope pair $(r_1,r_2)$; if $(r_1,r_2)=(0,0)$, we simply write $\link_0$ and say
$A$ has a \emph{trivial slope pair}. 
The slope pair is of either the form $(\frac{p}{q},\frac{q}{p}),pq\neq 0$ or the form $(\frac{p}{q},pq),q\neq 0$,
where $p,q$ are coprime integers by \cite[Lemma $2.12$]{Wan:22}. 

\item[\emph{Type $4$-$1$}]: Components of $\partial A$
are parallel in $\partial \HK$ and every essential disk in $\HK$ meets $\partial A$. Note that the core of the solid torus cut off by $A$ from $\Compl\HK$ 
is an Eudave-Mu\~noz knot \cite{Eud:97}. 
\item[--] For a type $4$-$1$ annulus the symbol $\EM$ is reserved. 
\end{itemize} 

Label each edge of $\charE$, based on the type 
of the annulus it represents.
Then the resulting edge-labeled diagram, denoted by $\charhk$, is called the \emph{annulus diagram} of $\pair$. The annulus diagram contains finer information; for instance,
$(\sphere,5_1)$ and $(\sphere,6_4)$ in the Ishii-Kishimoto-Moriuchi-Suzuki handlebody-kont table \cite{IshKisMorSuz:12} have homeomorphic exteriors but different annulus diagrams (Figs.\ \ref{intro:fig:ann_fiveone} and \ref{intro:fig:ann_sixfour}). 
By the definition, an essential annulus $A\subset\Compl\HK$ is non-separating if and only if
$A$ is of type $2$ or of type $3$-$3$.


\subsection{Mani results}
We classify the annulus diagrams of atoroidal handlebody-knots admitting an essential annulus of type $2$ or 
of type $3$-$3$ with specific slope pairs.
\begin{theorem}[Theorem \ref{teo:classification_char_diagram_type_two}, Proposition \ref{prop:realization}]\label{intro:teo:classification}
Suppose $\pair$ is atoroidal and $\ComplHK$ 
admits a type $2$ essential annulus $A$. 
\begin{enumerate}[label=\textnormal{(\roman*)}]
\item\label{itm:char_diagram_twoone}  
If $A$ is of type $2$-$1$, then  
$\charhk$ is \raisebox{-.2 cm}{\includegraphics[scale=.13]{typetwoone_ann}}.
\item\label{itm:char_diagram_twotwo} If $A$ is of type $2$-$2$, then $\charhk$ is one of the following:
\begin{figure}[H]
\begin{subfigure}{0.32\linewidth}
\centering
\includegraphics[scale=.12]{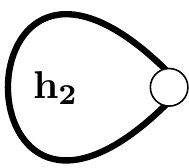}
\end{subfigure}
\begin{subfigure}{0.32\linewidth}
\centering
\includegraphics[scale=.13]{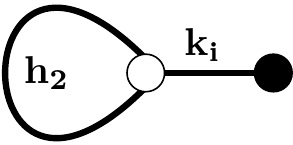} 
\quad 
\raisebox{.2 cm}{\footnotesize{$i=1$ or $2$}}
\end{subfigure}
\begin{subfigure}{0.32\linewidth}
\centering
\includegraphics[scale=.12]{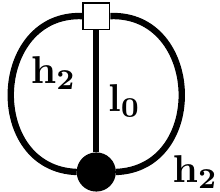} 
\quad \raisebox{.4 cm}{\footnotesize{$\square=\snode$ or $\hnode$}}.
\end{subfigure}
\end{figure} 
\item\label{itm:typetwo_realization} Every diagram in \ref{itm:char_diagram_twoone} and
\ref{itm:char_diagram_twotwo} can be realized by 
some atoroidal handlebody-knot.
\end{enumerate}
\end{theorem}
In the case the characteristic diagram $\charE$ 
is of $\theta$-shape, we show that the annulus diagram
$\charhk$ is determined by $\charE$, and obtain
a characterization of the simplest non-trivial atoroidal handlebody-knot in terms of the characteristic diagram.
%
\begin{theorem}[Theorems \ref{teo:char_three_bigons}, \ref{teo:characterization_fourone}]\label{intro:teo:characterization_theta_char}
Suppose $\pair$ is atoroidal. 
\begin{enumerate}[label=\textnormal{(\roman*)}]
\item 
If $\charE$ is \raisebox{-.2cm}{\includegraphics[scale=.12]{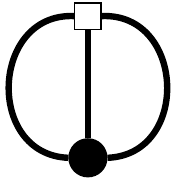}}, then 
the annulus diagram $\charhk$ is 
\raisebox{-.2cm}{\includegraphics[scale=.12]{typetwotwo_ann_3}}, 
where {\footnotesize $\square=\hnode$ or $\snode$}.
\item\label{itm:hkfourone_char} 
If $\charE$ is \raisebox{-.2cm}{\includegraphics[scale=.11]{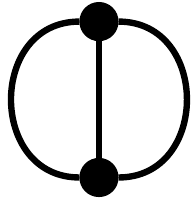}}, then $\pair$ is equivalent to  $\pairfourone$
in the handlbody-knot table \cite{IshKisMorSuz:12}. 
\end{enumerate}
\end{theorem}

For a type $3$-$3$ annulus $A$, we have the following
partial classification.   
\begin{theorem}[Corollaries \ref{cor:unique_annulus}, \ref{cor:char_annuli}, Lemma \ref{lm:typethreethree_trivial_slope}]\label{intro:teo:typethreethree} 
Suppose $\pair$ is atoroidal, and 
$A\subset \ComplHK$ a type $3$-$3$ essential annulus. 
\begin{enumerate}[label=\textnormal{(\roman*)}]
\item\label{itm:typethree_rational} If $A$ has a boundary slope pair of $(\frac{p}{q},\frac{q}{p})$, $pq\neq 0$, then $\charhk$ is
\raisebox{-.25cm}{\includegraphics[scale=.11]{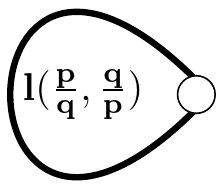}}. 
\item\label{itm:typethree_trivial} If $A$ has a trivial slope pair, then $\charhk$ is 
\raisebox{-.25cm}{\includegraphics[scale=.11]{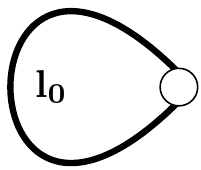}} or 
\raisebox{-.2cm}{\includegraphics[scale=.12]{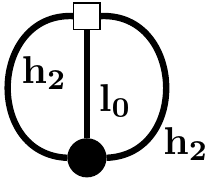}}.
\end{enumerate}
\end{theorem}
We remark that \ref{itm:typethree_rational} is Corollary \ref{cor:unique_annulus}, and \ref{itm:typethree_trivial}
follows from Theorem \ref{intro:teo:classification}, Lemma \ref{lm:typethreethree_trivial_slope}
and Corollary \ref{cor:char_annuli}. 
Also, Theorem \ref{intro:teo:char_diagram}, Theorem \ref{intro:teo:characterization_theta_char}, and 
Lemma \ref{lm:nonchar_annuli} imply
that $\Compl\HK$ can admit at most two type $3$-$3$
essential annuli, up to isotopy, and should this happen, both would have the same boundary slope pair $(\frac{p}{q},pq)$ with $\vert p\vert$ greater than $1$.

Applying Theorem \ref{intro:teo:classification}, 
we compute the symmetry group for atoroidal handlebody-knots
whose exteriors contain a type $2$ annulus.
\begin{theorem}[Theorems \ref{teo:symmetry_type_2_1}
--\ref{teo:symmetry_two_type_2_2}]\label{intro:teo:symmetry}
Suppose $\pair$ is atoroidal and $A\subset\ComplHK$ a type $2$ essential annulus. 
\begin{enumerate}[label=\textnormal{(\roman*)}]
\item If $A$ is of type $2$-$1$,
then $\pSym\HK<\Z_2$ and $\Sym\HK<\Z_2\times \Z_2$.
\item If $A$ is the unique type $2$-$2$ 
annulus in $\ComplHK$, up to isotopy,
then $\pSym\HK\simeq \{1\}$ and $\Sym\HK<\Z_2$.
\item If $A$ is the unique type $2$-$2$ annulus,
but not the unique annulus in $\ComplHK$, up to isotopy, 
then $\pSym\HK\simeq \{1\}\simeq \Sym\HK$.
\item\label{intro:itm:two_typetwotwos} If $A$ is not
the unique type $2$-$2$ annulus, up to isotopy, then
$\pSym\HK<\Z_2$ and $\Sym\HK<\Z_2\times \Z_2$.
\end{enumerate}
\end{theorem}

\begin{figure}
\centering
\includegraphics[scale=.15]{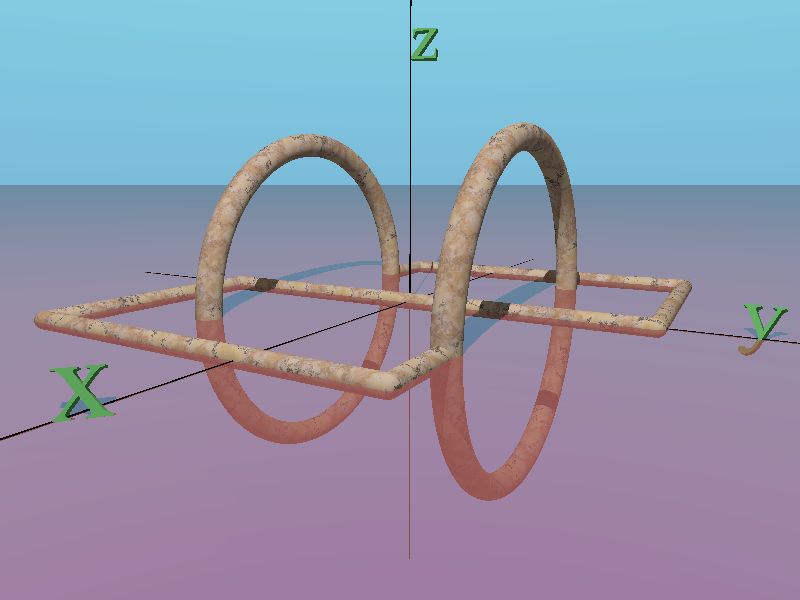}
\caption{Rigid symmetries of $\pairfourone$.}
\label{intro:fig:rigid_symmetry_hk4_1}
\end{figure}

Note the difference between
\emph{``unique annulus''} and \emph{``unique type XXX annulus''}: in the latter, annuli of other types might exist. 
Theorem \ref{intro:teo:symmetry} implies $\Sym{4_1}\simeq \Z_2\times \Z_2$ and $\pSym{4_1}\simeq \Z_2$ as  
the reflection against the xy-plane and 
rotation around the z-axis by $\pi$ in
Fig.\ \ref{intro:fig:rigid_symmetry_hk4_1} 
represent two non-trivial mapping classes. 
To our knowledge, $\pairfourone$ is the only known example that
attains the upper bound in Theorem \ref{intro:teo:symmetry}
\ref{intro:itm:two_typetwotwos}; 
on the other hand, no handlebody-knot admitting 
a unique type $2$ annulus has been found 
to have a non-trivial symmetry group so far.
We speculate the following sharper statements both are true.

\begin{problem}
Under the same assumption as in Theorem \ref{intro:teo:symmetry},
$\Sym\HK\simeq \Z_2\times \Z_2$ 
if and only if  $\pair$ is equivalent to $\pairfourone$.
\end{problem}
\begin{problem}
Under the same assumption as in Theorem \ref{intro:teo:symmetry},
suppose $A$ is the unique type $2$ annulus in $\ComplHK$, up to isotopy. 
Then $\Sym\HK\simeq \{1\}$.
\end{problem}
 
The rigid motions shown in Fig.\ \ref{intro:fig:rigid_symmetry_hk4_1} 
suggest a variant of the Nielsen realization problem.
\begin{problem}\label{intro:pro:nielsen_realization}
Let $\pair$ be a non-trivial atoroidal handlebody-knots.
Then there exists a subgroup $G<\Aut{\sphere,\HK}$
such that $\pi_0:\Aut{\sphere,\HK}\rightarrow \Sym\HK$
restricts to an isomorphism on $G$.
\end{problem}
 
Handlebody-knot symmetry is itself a topic of independent interest. 
To our knowledge, apart from $\pairfourone$, 
the symmetry group is computed for only five other handlebody-knots 
in the table \cite{IshKisMorSuz:12}: 
\begin{multline*}
\Sym{5_1}\simeq \Sym{6_1}\simeq \Sym{6_{11}}\simeq \{1\},\\
\Sym{5_2}\simeq \pSym{5_2}\simeq \mathbb{Z}_2\times \mathbb{Z}_2,\quad\Sym{6_4}\simeq \pSym{6_4}\simeq \mathbb{Z}_2.
\end{multline*} 
The first two are computed by Koda
\cite{Kod:15} using results from Motto \cite{Mott:90} and Lee-Lee \cite{LeeLee:12}, 
while the third follows from \cite{Wan:21}
and Theorem \ref{intro:teo:classification};
the last two are computed in \cite{Wan:22}.
They all can be realized as subgroups of the homeomorphism groups.

To prove Theorem \ref{intro:teo:classification}\ref{itm:typetwo_realization}, 
we need to produce 
\emph{atoroidal} handlebody-knots admitting a type $2$
\emph{essential} annulus---a type $2$ 
annulus is not necessarily essential 
by the definition. Sections \ref{sec:irre_atoro} and \ref{sec:examples} develop essentiality and
atoroidality tests and present a systematical approach,
via spatial graphs, to generate atoroidal handlebody-knots admitting a type $2$ essential annulus.
 
Our tests make use of an \emph{unknotting operation}: given a type $2$ annulus $A\subset\Compl\HK$, 
then the union $\HKA:=\HK\cup\rnbhd{A}$
induces a handlebody-knot $\pairA$,
where $\rnbhd{A}\subset\Compl\HK$ is a regular neighborhood of $A$. 
The frontier of $\rnbhd{A}\subset\Compl\HK$
consists of two annuli in $\partial\HKA$, whose 
cores we denote by $l_+,l_-\subset\partial\HKA$.
Recall also that a set of disjoint simple loops
$\{l_1,\dots,l_n\}$ in the boundary of a $3$-manifold $M$ is \emph{primitive}
if there exists a set of disjoint disks $\{D_1,\dots,D_n\}$
in $M$ 
such that $l_i\cap \partial D_j$ is a point when $i=j$ and empty otherwise. Our essentiality and atoroidality criteria are stated as follows.
  
\begin{theorem}[Propositions \ref{prop:irre_atoro_type2_1} and \ref{prop:irre_atoro_type2_2})]\label{intro:teo:irre_atoro_criteria}
Given a handlebody-knot $\pair$,
and a type $2$ annulus $A\subset\ComplHK$.
\begin{enumerate}[label=\textnormal{(\roman*)}]
\item
Suppose $A$ is of type $2$-$1$.
Then $\pair$ is atoroidal and $A$ is essential 
if and only if $\pairA$ is trivial 
with $\{l_+,l_-\}$ not primitive in $\ComplHKA$
or is non-trivial and atoroidal.
\item 
Suppose $A$ is of type $2$-$2$.
Then $\pair$ is atoroidal and $A$ is essential 
if and only if $\pairA$ 
is trivial with $l_+,l_-$ not homotopically trivial in $\ComplHKA$ or is non-trivial and atoroidal. 
\end{enumerate}
\end{theorem}
 
\subsection{Convention} 
We work in the piecewise linear category.
Given a subpolyhedron $X$ of $M$, we denote by $\overline{X}$, 
$\mathring{X}$, $\mathfrak{N}(X)$, and 
$\fr{X}{M}$ the closure, the interior, a regular neighborhood, and 
the frontier of $X$ in $M$, respectively.
The \emph{exterior} $\Compl X$ of $X$ in $M$ is 
defined to be 
the complement of $\openrnbhd{X}$ 
if $X\subset M$ is of positive codimension, 
and defined to be the closure of $M-X$ otherwise. 
Submanifolds of a manifold $M$ are assumed to be
proper and in general position except in some obvious cases
where submanifolds are in $\partial M$.  
A surface $S$ other than a disk
in a three-manifold $M$ is \emph{essential} 
if it is incompressible and $\partial$-incompressible.
A disk $D\subset M$ is \emph{essential} 
if $D$ does not cut a $3$-ball off from $M$.
When $M$ is a handlebody-knot, 
an essential disk is also called a \emph{meridian} disk. $3$-manifolds here are assumed to be orientable.
Given a loop $l$ in a space $X$, $[l]$ denotes the homology class represented by $l$ in $H_1(X)$.  
We denote by $(\sphere,X)$ an embedding 
of $X$ in the oriented $3$-sphere $\sphere$.


\begin{table}[t]
\caption{Complete list of characteristic diagrams of $\ComplHK$.}
\label{tab:configurations}
\begin{subfigure}[l]{0.3\textwidth}
\centering
\vspace*{.5cm}{
\begin{tikzpicture}
      \tikzset{enclosed/.style={draw, circle, inner sep=0pt, minimum size=.15cm, fill=black}}
      \tikzset{hollow/.style={draw, circle, inner sep=0pt, minimum size=.15cm, fill=white}}

      \node[hollow
      ] (E) at (0,1) {};
\end{tikzpicture}

\caption{$(0,0,0,\hnode)$}
\label{fig:000h}
} 
\end{subfigure}
\begin{subfigure}[l]{0.3\textwidth}
\centering
\vspace*{-.5cm}{
\begin{tikzpicture}
      \tikzset{enclosed/.style={draw, circle, inner sep=0pt, minimum size=.15cm, fill=black}}
      \tikzset{hollow/.style={draw, circle, inner sep=0pt, minimum size=.15cm, fill=white}}

      \node[hollow
      ] (E) at (0,1) {};
  \draw (E) to [out=140,in=-140, looseness=40] (E);
\end{tikzpicture}
}
\vspace*{-.7cm}{\caption{$(1,1,0,\hnode)$}
\label{fig:110h}
} 
\end{subfigure}
\begin{subfigure}[l]{0.3\textwidth}
\centering
\begin{tikzpicture}
      \tikzset{solid/.style={draw, circle, inner sep=0pt, minimum size=.15cm, fill=black}}
      \tikzset{hollow/.style={draw, circle, inner sep=0pt, minimum size=.15cm, fill=white}}
      \node[hollow
      ] (main) at (0.4,0.3) {};
      \node[solid] (node1) at (-.4,-.3) {};
      \draw (main) to (node1);
\end{tikzpicture}
\caption{$(1,0,0,\hnode)$}
\label{fig:100h}
\end{subfigure}
\begin{subfigure}[l]{0.3\textwidth}
\centering
\begin{tikzpicture}
      \tikzset{solid/.style={draw, circle, inner sep=0pt, minimum size=.15cm, fill=black}}
      \tikzset{hollow/.style={draw, circle, inner sep=0pt, minimum size=.15cm, fill=white}}
      \node[solid
      ] (main) at (.4,.3) {};
      \node[solid] (node1) at (-.4,-.3) {};
      \draw (main) to (node1);
\end{tikzpicture}
\caption{$(1,0,0,\snode)$}
\label{fig:100s}
\end{subfigure}
\begin{subfigure}[l]{0.3\textwidth}
\centering
\begin{tikzpicture}
      \tikzset{solid/.style={draw, circle, inner sep=0pt, minimum size=.15cm, fill=black}}
      \tikzset{hollow/.style={draw, circle, inner sep=0pt, minimum size=.15cm, fill=white}}

      \node[hollow
      ] (main) at (0,0) {};
      \node[solid] (node) at (1,0) {};
  \draw (main) to [out=140,in=-140, looseness=30] (main);
  \draw (main) to (node);
\end{tikzpicture}
\vspace*{-.3cm}{
\caption{$(2,1,0,\hnode)$}
\label{fig:210h}
}
\end{subfigure}
\begin{subfigure}[l]{0.3\textwidth}
\centering
\begin{tikzpicture}
      \tikzset{solid/.style={draw, circle, inner sep=0pt, minimum size=.15cm, fill=black}}
      \tikzset{hollow/.style={draw, circle, inner sep=0pt, minimum size=.15cm, fill=white}}

      \node[hollow
      ] (main) at (0.5,0.5) {};
      \node[solid] (node) at (-.5,-.5) {};
  \draw (main) to [out=180,in=90] (node);
  \draw (main) to [out=-90,in=0] (node);
\end{tikzpicture}
\caption{$(2,0,1,\hnode)$}
\label{fig:201h}
\end{subfigure}
 
\begin{subfigure}[l]{0.3\textwidth}
\centering
\begin{tikzpicture}
      \tikzset{solid/.style={draw, circle, inner sep=0pt, minimum size=.15cm, fill=black}}
      \tikzset{hollow/.style={draw, circle, inner sep=0pt, minimum size=.15cm, fill=white}}
      \node[hollow
      ] (main) at (0,-.7) {};
      \node[solid] (node1) at (-.5,0) {};
      \node[solid] (node2) at (.5,0) {};
      \draw (main) to (node1);
      \draw (main) to (node2);
\end{tikzpicture}
\caption{$(2,0,0,\hnode)$}
\label{fig:200h}
\end{subfigure}
\begin{subfigure}[l]{0.3\textwidth}
\centering
\begin{tikzpicture}
      \tikzset{solid/.style={draw, circle, inner sep=0pt, minimum size=.15cm, fill=black}}
      \tikzset{hollow/.style={draw, circle, inner sep=0pt, minimum size=.15cm, fill=white}}
      \node[solid
      ] (main) at (0,-.7) {};
      \node[solid] (node1) at (-.5,0) {};
      \node[solid] (node2) at (.5,0) {};
      \draw (main) to (node1);
      \draw (main) to (node2);
\end{tikzpicture}
\caption{$(2,0,0,\snode)$}
\label{fig:200s}
\end{subfigure}
\begin{subfigure}[l]{0.3\textwidth}
\centering
\begin{tikzpicture}
      \tikzset{solid/.style={draw, circle, inner sep=0pt, minimum size=.15cm, fill=black}}
      \tikzset{hollow/.style={draw, circle, inner sep=0pt, minimum size=.15cm, fill=white}}

      \node[hollow
      ] (main) at (0,0.25) {};
      \node[solid] (node) at (0,-1) {};
  \draw (main) to [out=-160,in=160] (node);
  \draw (main) to [out=-20,in=20] (node);
  \draw (main) to (node);
\end{tikzpicture}
\caption{$(3,0,3,\hnode)$}
\label{fig:303h}
\end{subfigure}
\begin{subfigure}[l]{0.3\textwidth}
\centering
\begin{tikzpicture}
      \tikzset{solid/.style={draw, circle, inner sep=0pt, minimum size=.15cm, fill=black}}
      \tikzset{hollow/.style={draw, circle, inner sep=0pt, minimum size=.15cm, fill=white}}

      \node[solid
      ] (main) at (0,0.25) {};
      \node[solid] (node) at (0,-1) {};
  \draw (main) to [out=-160,in=160] (node);
  \draw (main) to [out=-20,in=20] (node);
  \draw (main) to (node);
\end{tikzpicture}
\caption{$(3,0,3,\snode)$}
\label{fig:303s}
\end{subfigure}
\begin{subfigure}[l]{0.3\textwidth}
\centering
\begin{tikzpicture}
      \tikzset{solid/.style={draw, circle, inner sep=0pt, minimum size=.15cm, fill=black}}
      \tikzset{hollow/.style={draw, circle, inner sep=0pt, minimum size=.15cm, fill=white}}

      \node[hollow
      ] (main) at (0,0) {};
      \node[solid] (node1) at (.5,.5) {};
      \node[solid] (node2) at (-.75,-.75) {};
  \draw (main) to [out=180,in=90] (node2);
  \draw (main) to [out=-90,in=0] (node2);
  \draw (main) to (node1);
\end{tikzpicture}
\caption{$(3,0,1,\hnode)$}
\label{fig:301h}
\end{subfigure}
\begin{subfigure}[l]{0.3\textwidth}
\centering
\begin{tikzpicture}
      \tikzset{solid/.style={draw, circle, inner sep=0pt, minimum size=.15cm, fill=black}}
      \tikzset{hollow/.style={draw, circle, inner sep=0pt, minimum size=.15cm, fill=white}}

      \node[solid
      ] (main) at (0,0) {};
      \node[solid] (node1) at (.5,.5) {};
      \node[solid] (node2) at (-.75,-.75) {};
  \draw (main) to [out=180,in=90] (node2);
  \draw (main) to [out=-90,in=0] (node2);
  \draw (main) to (node1);
\end{tikzpicture}
\caption{$(3,0,1,\snode)$}
\label{fig:301s}
\end{subfigure}
\begin{subfigure}[l]{0.3\textwidth}
\centering
\begin{tikzpicture}
      \tikzset{solid/.style={draw, circle, inner sep=0pt, minimum size=.15cm, fill=black}}
      \tikzset{hollow/.style={draw, circle, inner sep=0pt, minimum size=.15cm, fill=white}}

      \node[hollow
      ] (main) at (0,0) {};
      \node[solid] (node1) at (.4,.3) {};
      \node[solid] (node2) at (0,-.5) {};
      \node[solid] (node3) at (-.4,.3) {};
  \draw (main) to (node1);
  \draw (main) to (node2);
  \draw (main) to (node3);
\end{tikzpicture}
\caption{$(3,0,0,\hnode)$}
\label{fig:300h}
\end{subfigure}
\begin{subfigure}[l]{0.3\textwidth}
\centering
\begin{tikzpicture}
      \tikzset{solid/.style={draw, circle, inner sep=0pt, minimum size=.15cm, fill=black}}
      \tikzset{hollow/.style={draw, circle, inner sep=0pt, minimum size=.15cm, fill=white}}

      \node[solid
      ] (main) at (0,0) {};
      \node[solid] (node1) at (.4,.3) {};
      \node[solid] (node2) at (0,-.5) {};
      \node[solid] (node3) at (-.4,.3) {};
  \draw (main) to (node1);
  \draw (main) to (node2);
  \draw (main) to (node3);
\end{tikzpicture}
\caption{$(3,0,0,\snode)$}
\label{fig:300s}
\end{subfigure}
\end{table}



\section{Characteristic submanifolds}\label{sec:charac}
Here we review Johannson's characteristic submanifold
theory \cite{Joh:79} (see also \cite{CanMcC:04}),
and introduce the characteristic diagram and annulus diagram. A completeness criteria needed in 
Section \ref{sec:classification} is also developed.

\subsection{Characteristic submanifold theory}
\begin{definition}
Given a compact $n$-manifold $M$,
a boundary-pattern $\bp m$ for $M$  
is a finite set of compact, connected $(n-1)$-submanifolds of $\partial M$
such that the intersection of 
any $i$ of them is either empty or an $(n-i)$-manifold.
\end{definition}

We denote by $\vert\bp m\vert$ the union of $\elem G$, 
$\elem G\in\bp m$.
An $i$-faced disk is a disk 
$D$ whose boundary-pattern $\bp d$ 
consists of $i$ elements
with $\vert \bp d\vert=\partial D$. 
When $i\leq 3$ (resp.\ $i=4$), 
$(D,\bp d)$ is called a small-faced disk (resp.\ a square). 
The empty boundary-pattern is denoted by $\bp \phi$, and 
the completion $\cbp m$ 
of a boundary-pattern $\bp m$ for $M$
is the boundary-pattern given by 
\[
\cbp m:=\{\elem G \in \bp m\}\cup \{\text{components of $\overline{M-\vert \bp m\vert}$}\}.
\]
 
Throughout the paper, 
an annulus (or arc) is assumed to carry the boundary-pattern $\cbp\phi$. 
Given a manifold $(M,\bp m)$ with boundary-pattern,
and a submanifold $N\subset M$ of positive codimension,
if $N\cap \partial M$ meets every intersection of elements of 
$\cbp m$ transversely, then $N$
inherits a natural boundary-pattern given by
\begin{equation}\label{eq:submfd_b_pattern}
\bp n:=\{\elem G\cap \partial N\mid
\forall  \elem G\in \bp m\}.
\end{equation} 
Similarly, $\bp n$ defines a boundary-pattern
for a codimension-zero submanifold $N$ of $M$, 
provided the intersection $\fr{N}{M}\cap \partial M$ 
meets every intersection of elements in $\cbp m$ transversely.
The boundary-pattern $\bp n$ for $N$ is called
the \emph{submanifold boundary-pattern};
when $N$ is of codimension zero, we call 
the completion $\cbp n$ 
the \emph{proper boundary-pattern} for $N$.
Throughout the paper, a submanifold $N\subset M$ is assumed to satisfy the transversality condition, and
unless otherwise specified, $N$ always carries 
the submanifold boundary-pattern $\bp n$. We drop $\bp n$
from the notation when there is no
risk of confusion, but specify in the notation the 
proper boundary-pattern $\cbp n$ whenever useful.
When $N$ is considered as 
the exterior
$\Compl W$ of some submanifold $W$ in $M$, 
the proper boundary-pattern is assumed and denoted by $\pbp m$.   


\begin{definition}\label{def:essentiality}
An \emph{arc} $\gamma$ in a surface $(S,\bp{s})$
with boundary-pattern is essential if 
no component of $(\Compl \gamma,\pbp s)$ 
is a small-faced disk.

A \emph{surface} $S$ in a $3$-manifold $(M,\bp m)$ 
with boundary-pattern is essential 
if no component $X$ of $(\Compl S,\pbp m)$ contains a small-faced disk that meets the frontier $\bm X$ in an essential arc in
$\bm X$.  
A \emph{codimension-zero submanifold} $N$ in $(M,\bp m)$
is essential if its frontier $\fr{N}{M}$ is essential
in $(M,\bp m)$.
\end{definition} 
In the case $\bp m=\bphi$, the definition is equivalent to the one
in terms of incompressibility and $\partial$-incompressibility.
A $3$-manifold $(M,\bp m)$ with boundary-pattern can be 
\emph{I-fibered} (resp.\ \emph{Seifert fibered})
if it admits an I-bundle (resp.\ Seifert bundle) structure
$X\xrightarrow{\pi} B$ with $B$ equipped with 
a boundary-pattern $\bp b$ such that 
\begin{equation}\label{eq:fiber_bundle_structure}
\bp m=\{ \pi^{-1}(\elem G)\mid \elem G\in \bp b\}\cup \{\text{components of $\overline{\partial M-\pi^{-1}(\partial B)}$}\}.
\end{equation}
If $(M,\bp m)$ is I-fibered over $(B,\bp b)$, a component 
of $\overline{\partial M-\pi^{-1}(\partial B)}$
is called a \emph{lid} of $(M,\bp m)$ (with respect to $\pi$),
and any other element in $\bp m$ is called a \emph{side} of
$(M,\bp m)$ (with respect to $\pi$).
$(M,\bp m)$ is called a \emph{cylindrical shell}
if it can be I-fibered over an annulus.
An annulus $A$ in $(M,\bp m)$ is parallel to an element $\elem A\in \bp m$ (resp.\ to another annulus $A'$ in $(M,\bp m)$) 
if a component of 
$(\Compl{A\cup \elem A},\pbp m)$ (resp.\ $(\Compl {A\cup A'},\pbp m)$) 
is a cylindrical shell 
meeting both the regular neighborhoods of $A$ and of $\elem A$ (resp.\ of $A'$).
The following is a corollary of the vertical-horizontal theorem 
\cite[Proposition $5.6$; Corollary $5.7$]{Joh:79}.
\begin{lemma}\label{lm:vertical}
Suppose $(M,\bp m)$ is I-fibered over $(B,\bp b)$
with $\chi(B)<0$. Let $A$ be 
an essential annulus in $(M,\bp m)$.
Then the boundary $\partial A$ 
is in the lid(s) $\Li\in \bp m$, and
there exists an isotopy $F_t:(A,\partial A)\rightarrow 
(M,\Li)$ with $F_0=\id$ and 
$F_1(A)$ the preimage of an 
essential loop in $B$.
\end{lemma}

\begin{definition}
An $\F$-manifold $W$ in $(M,\bp m)$ is a codimension-zero 
essential submanifold of $M$ such that each component of $W$ 
can be I- or Seifert fibered.
An $\F$-manifold $W$ in $M$ is full if 
there exists no component $Y$ of $\Compl W$ 
such that $Y\cup W$ is an $\F$-manifold in $(M,\bp m)$. 
\end{definition}

\begin{definition}\label{def:completeness}
An $\F$-manifold $W$ in $(M,\bp m)$ is complete if
for any component $Y$ of $(\Compl W,\pbp m)$
and any essential square, annulus or torus $S$ in $Y$,
one of the following holds.
\begin{align}
\label{cond:completeness1}\tag{C1} &\text{If $S\cap \fr{Y}{M}\neq\emptyset$, then
$Y$ can be fibered as a product I-bundle or $S^1$-bundle over $S$.}\\
\label{cond:completeness2}\tag{C2}
&\text{If $S\cap \fr{Y}{M}=\emptyset$, then $S$
is parallel to a component of $\fr{Y}{M}$
in $Y$.}
\end{align}
\end{definition}

\begin{definition}
A characteristic submanifold $W$ for $(M,\bp m)$ is
a full, complete $\F$-manifold in $(M,\bp m)$.
\end{definition}


\subsection{Characteristic submanifiolds of atoroidal manifolds}
Here we consider the case where $M$ 
is an orientable, irreducible,
$\partial$-irreducible, atoroidal $3$-manifold 
equipped with the boundary-pattern $\bphi$; $\bphi$ is dropped from the notation when no confusion may arise. 
The existence and uniqueness of characteristic submanifolds are guaranteed in this case. 
\begin{theorem}[{\cite[Proposition $9.4$; Corollary $10.9$]{Joh:79}}]\label{teo:unique_char_submfd}
There exists a characteristic submanifold $W$ for $M$,
and two characteristic submanifolds $W_1,W_2$ for $M$
are ambient isotopic.
\end{theorem} 

Furthermore, they have 
the engulfing property.
\begin{theorem}[{\cite[Proposition $10.8$]{Joh:79}}]\label{teo:engulfing}
Let $W$ be a characteristic submanifold for $M$.
Then, for every $\F$-manifold $X\subset M$, 
there exists an ambient isotopy $F_t$ 
such that $F_1(X)\subset W$. 
\end{theorem}  

The following, a direct consequence 
of \cite[Theorem $2.9.3$]{CanMcC:04}, gives an 
alternative description of characteristic submanifolds in
terms of simple manifolds. 
\begin{definition}
A manifold $(X,\bp x)$ with boundary-pattern
is simple if any component of 
a characteristic submanifold of $(X,\cbp x)$
is a regular neighborhood of 
a square, annulus or torus in $\bp x$.
\end{definition}
%
%
%
%
\begin{theorem}\label{teo:completeness_criterion_via_simpleness}
Given a full $\F$-manifold $W\subset M$, 
$W$ is a characteristic submanifold for $M$ 
if and only if, for every component $Y\subset (\Compl W,\pbp m)$, $Y$ either is simple or is 
a cylindrical shell.
\end{theorem}

We examine topological properties of 
submanifolds of $M$
that can be I- or Seifert fibered.

\begin{lemma}\label{lm:genus_one_component}
Let $X$ be an essential codimension-zero submanifold 
of $M$. 
Then $g(\partial X)=1$ if and only if
$(X,\cbp x)$ can be Seifert fibered over an
$n$-faced disk with at most one exceptional fiber, 
and $\bp x$ non-empty and containing disjoint elements;
additionally, it has 
exactly one exceptional fiber when $n=2$.
\end{lemma}
\begin{proof}
The direction ``$\Leftarrow$'' is clear. 
To see the direction ``$\Rightarrow$'', 
note first that by the essentiality of $X$
and the boundary-pattern $\bphi$ on $M$,
the intersection 
$X\cap \partial M$ is non-empty
and consists of disjoint annuli $A_1,\dots, A_m$ in $\partial X$.
This implies $X$ is a solid torus
by the atoroidality of $M$.
Since $M$ is $\partial$-irreducible,
$H_1(A_i)\rightarrow H_1(X)$ cannot be trivial, and
therefore, $(X,\cbp x)$ can be Seifert fibered over an
$n$-faced disk $(D,\bp d)$ with $n=2m>0$. 
In the case $n=2$, by the essentaility 
of $\bm X$, the Seifert fibering must 
contain at least one exceptional fiber.
\end{proof}

\begin{corollary}\label{cor:essen_annuli_in_seifert}
Let $X\subset M$ be an essential codimension-zero submanifold 
with $g(\partial X)=1$.  
Then $(X,\cbp{x})$ admits an essential annulus
meeting $\bm X$.
\end{corollary}
\begin{proof}
By Lemma \ref{lm:genus_one_component}, the frontier 
of a regular neighborhood of an element in $\bp x$
is an essential annulus meeting $\bm X$. 
\end{proof}

\begin{lemma}
Given an essential codimension-zero submanifold $X\subset M$, if $(X,\bp x)$ is I-fibered over $(B,\bp b)$, then $\bp b=\bp \phi$, that is, $\bp x$ consists of only lids.
\end{lemma}
\begin{proof}
By the definition \eqref{eq:fiber_bundle_structure}, 
the lid(s) of $(X,\bp x)$ is(are) element(s) in $\bp x$.
On the other hand, 
since the boundary pattern on $M$ is $\bphi$, 
the submanifold boundary-pattern $\bp x$ 
consists of disjoint elements. Thus $\bp x$ only 
contain the lid(s).
\end{proof}

\begin{lemma}\label{lm:disjoint_parallel_annuli_in_I_bundle}
Let $(X,\bp x)\xrightarrow{\pi} (B,\bphi)$ be an I-bundle and $g(\partial X)>1$.
Then every essential annulus in $(X,\bp x)$ 
disjoint from the sides of $(X,\bp x)$ 
is parallel to a side $\elem A\in \bp x$ 
if and only if $B$ is a pair of pants. 
\end{lemma}
\begin{proof}
The direction ``$\Leftarrow$''
follows from Lemma \ref{lm:vertical}.
We prove the direction ``$\Rightarrow$'' by contradiction. 
Observe first that since $g(\partial X)>1$, 
the Euler characteristic $\chi (B)\leq -1$ 
by the equality $2\chi(B)=2-2g(\partial X)$. 
In particular, $B$ is a closed surface $\hat B$ 
with $k$ open disks removed such that $k$ and the genus 
$g(\hat B)$ satisfy
$3-2g(\hat B)\leq k$ when $B$ is orientable and
$3-g(\hat B)\leq k$ otherwise.
Let $l$ be a non-separating loop in $B$
if $\hat B$ is neither a $2$-sphere nor a projective plane,
or a loop cutting a Mobius band off from $B$ if $\hat B$ 
is a projective space, or a loop
cutting a pair of pants off from $B$ if
$\hat B$ is a $2$-sphere.   
Then if $B$ is not a pair of pants,
the preimage of $l$
is an essential annulus in $X$
disjoint from the sides and not parallel
to any side of $(X,\bp x)$. 
\end{proof}

The following is a corollary of \cite[Proposition $4.6$]{Joh:79}.
\begin{lemma}\label{lm:essentiality_in_X_in_M}
Let $S\subset M$ be a surface consisting 
of essential annuli, and $X$ a component of $(\Compl S,\pbp m)$.
Then first $X$ is atoroidal, and secondly, 
given
an annulus $A\subset X$ disjoint from $\bm X$,
$A$ is essential in $X$ if and only if $A$ 
is essential in $M$.
\end{lemma}

\begin{theorem}[{\bf Completeness Criterion}]\label{teo:completeness_criterion}
Let $W\subset M$ be a full $\F$-manifold. Then
$W$ is complete if and only if, 
for every component $Y$ of $(\Compl W, \pbp m)$,
either $Y$ is a cylindrical shell or 
$g(\partial Y)>1$, $Y$ cannot be I-fibered 
over a pair or pants, and every essential annulus 
in $Y$ disjoint from $\bm Y$ is parallel
to a component of $\bm Y$.
\end{theorem}
\begin{proof}
``$\Rightarrow$'': Given a component $Y$ 
of $(\Compl W,\pbp m)$,
either $Y$ admits an essential square or annulus
that meets $\bm Y$ or it does not.
By \eqref{cond:completeness1} in Definition \ref{def:completeness}, 
$Y$ is a cylindrical shell if it is the former.
Suppose it is the latter. Then,
since $Y$ contains no essential square, it
cannot be I-fibered over a pair of pants, 
and by Corollary \ref{cor:essen_annuli_in_seifert}, 
$g(\partial Y)$ cannot be $1$. 
The rest follows directly from \eqref{cond:completeness2} of Definition \ref{def:completeness}.

``$\Leftarrow$'': It is clear that 
the conditions \eqref{cond:completeness1} and \eqref{cond:completeness2}
in Definition \eqref{def:completeness} are satisfied
if $Y$ is a cylindrical shell.
So, we suppose otherwise; 
by Theorem \ref{teo:completeness_criterion_via_simpleness},
it suffices to show that $Y$ is simple.
Let $W_y$ be the characteristic submanifold of $Y$;
note that since $Y$ is a component of $(\Compl W,\pbp m)$,
$Y\subset M$ is equipped with the proper boundary-pattern.
If $W_y=\emptyset$, then $Y$ is simple by the definition. 
If $W_y\neq\emptyset$ but $\fr{W_y}{Y}=\emptyset$,
then $Y=W_y$. Since $g(\partial Y)>1$, by 
Lemma \ref{lm:genus_one_component}, 
it cannot be Seifert fibered, so $Y$
admits an I-bundle structure, contradicting the assumption
by Lemma \ref{lm:disjoint_parallel_annuli_in_I_bundle}.

Suppose $\fr{W_y}{Y}\neq\emptyset$, and let 
$X_y$ be a component of $W_y$, and
$A$ be a component of 
the frontier $\fr{X_y}{Y}$.
Then $A$ is disjoint from 
$\bm Y$ since $W_y$ contains a regular neighborhood
of $\bm Y$ by Theorem \ref{teo:engulfing}. 
$A$ therefore cannot be a square by 
the boundary-pattern $\bphi$ on $M$;
neither can it be a torus because of Lemma \ref{lm:essentiality_in_X_in_M}.
$A$ is hence is an annulus. 
By the assumption, 
$A$ is parallel to a component $A'$ of $\bm Y$ in $Y$.
Let $P\subset Y$ be the cylindrical shell 
between $A$ and $A'$. 
Then by the fullness of
$W_y$, $P\supset X_y$ and $A'\subset X_y$.
On the other hand, the essentiality of $X_y$
implies $\partial_P X_y=\emptyset$,
so $P=X_y$. It follows that $W_y$ 
is the union of regular neighborhoods of components
in $\bm Y$, so $Y$ is simple.
\end{proof}

\subsection{Characteristic diagram}  
Let $M$ be as in the previous subsection. 
\begin{definition}[{\textbf{Characteristic Surfaces}}]
A characteristic surface $S$ of $M$ is a union of 
components of $\bm W$ 
such that 
\begin{itemize}
\item no two components of $S$ are parallel, and
\item every component of $\bm W$ is parallel to some component of $S$, 
\end{itemize}
where $W\subset M$ is a characteristic submanifold.
\end{definition}

The existence of a characteristic surface 
follows from the existence of 
a characteristic submanifold $W$ of $M$: 
for instance, a maximal subset of mutually 
non-parallel annuli in $\bm W$ 
is a characteristic surface. 
Characteristic surfaces of $M$ are unique, up to isotopy, by Theorem \ref{teo:unique_char_submfd}.
 
\begin{corollary}\label{cor:unique_char_surface}
Given two characteristic surfaces 
$S_1,S_2$ of $M$, there exists an ambient isotope
$F_t$ such that $F_1(S_1)=S_2$.
\end{corollary} 
Furthermore, by Theorem \ref{teo:completeness_criterion_via_simpleness}, 
every component of $\Compl S:=M-\openrnbhd{S}$ 
is either Seifert/I-fibered or simple.

\begin{definition}
Given a characteristic surface $S$ of $M$,
denote by $\Compl S$ the complement $M-\openrnbhd{S}$.
Then the associated characteristic diagram $\charM$
is a graph defined as follows:
\begin{itemize} 
\item Assign a solid node $\snode$ to each component
$\Compl S$ that can be I-or Seifert fibered.
\item Assign a hollow node $\hnode$ to each component
$\Compl S$ that is simple.
\item To each component of $\rnbhd{S}$,
assign an edge between node(s)
corresponding to component(s) of $\Compl S$ 
meeting the component of $\rnbhd{S}$.
\end{itemize}
\end{definition}
A node in $\charM$ or the component $X\subset \Compl S$ it represents is said to be 
of \emph{genus $g$} if $g(\partial X)=g$.
Two characteristic diagrams are \emph{isomorphic}
if there is a graph isomorphism between them 
sending solid (resp.\ hollow) nodes to solid (resp.\
hollow) nodes of the same genus.
By Corollary \ref{cor:unique_char_surface},  
the characteristic diagram $\charM$ of $M$ is determined by $M$, up to isomorphism.

We say an annulus $A\subset M$ 
is \emph{characteristic} if
it is isotopic to a component of a characteristic surface $S$ of $M$.
 
\subsection{Classification and annulus diagram}
Throughout the subsection,
$M$ is a $\partial$-irreducible, atoroidal $3$-submanifold in $\sphere$ with $g(\partial M)=2$, 
and $\charM$ is its characteristic diagram.  

\begin{proposition}\label{prop:properties_char_diagram_Ehk}\hfill
\begin{enumerate}[label=\textnormal{(\roman*)}]
\item\label{itm:one_label}
$\charM$ has exactly one genus two node, and 
all the other nodes are of genus one. 
\item\label{itm:unlabel_is_seifert} Genus one nodes in $\charM$ are all solid, and each corresponds
to a Seifert-fibered solid torus that is not a cylindrical shell.

\item\label{itm:no_solid_loops} No loop in $\charM$ contains a solid node.
\item\label{itm:adjacency} All edges in $\charM$ 
are adjacent to the genus two node.

\item\label{itm:I_bundles} If the genus two 
node in $\charM$ is solid, it corresponds to an I-bundle over a pair of pants or a Mobius band or Klein bottle with 
an open disk removed.
\item\label{itm:no_solid_bigon} 
If the genus two node in $\charM$ 
is solid, then $\charE$ cannot be a bigon.

\item\label{itm:trivalency} Every node in $\charM$ is at most trivalent.
\end{enumerate} 
\end{proposition}
\begin{proof}
Let $W$ be a characteristic submanifold of $M$
and $S$ a corresponding 
characteristic surface of $M$. 
Suppose the complement $\Compl S:=M-\openrnbhd{S}$ 
contains $n$ components
$X_1,\dots, X_n$.
Then the equality of Euler characteristic: 
\[-2=2-2g\big(\partial M\big)=\chi\big(\partial M\big)=
\sum_{i=1}^n\chi(\partial X_i)=
\sum_{i=1}^n \big(2-2g(\partial X_i)\big)\]
implies that
\[\sum_{i=1}^n\big(g(\partial X_i)-1\big)=1.\]
In particular, there exists exactly one genus $2$ 
component in $\Compl S$, and other components are 
of genus $1$ and hence Seifert-fibered by Lemma \ref{lm:genus_one_component} with none of them
a cylindrical shell by the definition of $S$. 
This proves \ref{itm:one_label} and \ref{itm:unlabel_is_seifert}.

We prove \ref{itm:no_solid_loops} by contradiction. 
Suppose there is a loop with a solid node in $\charM$, and
denote by $A$ the annulus corresponding to the loop,
and by $X\subset\Compl S$ the component corresponding to the solid node.
Then the union $X$ and $\rnbhd{A}$
is either Seifert-fibered or I-fibered, contradicting
the fullness of $W$.  

To see \ref{itm:adjacency}, it suffices to 
show there is no edge connecting two genus one solid nodes, given \ref{itm:no_solid_loops}.
Suppose such an edge exists, and 
$X_1,X_2\subset \Compl S$ be the Seifert components corresponding to the solid nodes. 
Let $A$ be the annulus corresponding 
to the edge. 
Then $X_1\cup\rnbhd{A}\cup X_2$
is Seifert fibered, contradicting 
the fullness of $W\subset M$.

For \ref{itm:I_bundles}, 
we observe first that the component $U$ in $\Compl S$ 
corresponding to a genus two solid node cannot be Seifert fibered by
Lemma \ref{lm:genus_one_component},
and hence is I-fibered. Since the lid(s)
of $U$ has Euler characteristic $-2$, the base
is either a pair of pants or 
a Mobius band, torus, or Klein bottle
with one open disk removed.
Suppose the base is 
a torus with one open disk removed.  
Then $\charM$ is 
\raisebox{0cm}{\includegraphics[scale=.12]{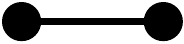}} 
by \ref{itm:adjacency}.
Denote by $A$ the annulus corresponding to the edge,
and let $V$ be the solid torus corresponding to the genus one node. Choose generators of $H_1(A)$ and $H_1(V)$ 
so that the homomorphism
$H_1(A)\simeq\mathbb{Z} \xrightarrow{m}\mathbb{Z}\simeq H_1(V)$ 
has $m\geq 0$. Since $A$ is essential, 
we have $m\neq 0,1$. 
The short exact sequence
\[0\rightarrow H_1(A)\xrightarrow{(m,0)}H_1(V)\oplus H_1(U)
\rightarrow H_1(M)\rightarrow 0\]
then implies $H_1(M)\simeq \Z\oplus \Z\oplus \Z_m$, contradicting $M\subset \sphere$.  

We prove \ref{itm:no_solid_bigon} 
by contradiction. Suppose $\charM$
is a bigon, and 
let $U$ (resp.\ $V$) be the components of $\Compl S$ 
corresponding to
the genus two (resp.\ genus one) 
node, and 
$A_1,A_2$ the annuli corresponding
to the edges.
Then by \ref{itm:I_bundles},
$U$ is an admissible I-bundle over
a Mobius band with one open disk removed.  
Choose generators of $H_1(A_i)$, $H_1(V), H_1(U)$ 
so that $H_1(A_i)\simeq \Z\xrightarrow{m} \Z\simeq H_1(V)$, 
$i=1,2$, 
and 
\[H_1(A_1)\simeq \Z \xrightarrow{\tiny \begin{pmatrix}
1\\0
\end{pmatrix}} \Z\oplus \Z\simeq H_1(U)
\quad\text{and}\quad 
H_1(A_2)\simeq \Z \xrightarrow{\tiny \begin{pmatrix}
\pm 1\\2
\end{pmatrix}} \Z\oplus \Z\simeq H_1(U).
\] 
Then by the exact sequence
\[0\rightarrow H_1(A_1\cup A_2)
\rightarrow H_1(V)\oplus H_1(U)
\rightarrow H_1(M)
\rightarrow \tilde{H}_1(A_1\cup A_2)
\rightarrow 0,
\]
either $(m,0,1)$ or $(0,0,1)$ in $H_1(V)\oplus H_1(U)$
induces an element of order $2$ in $H_1(M)$,
contradicting $M\subset \sphere$.

Lastly, in view of \ref{itm:adjacency},
to prove \ref{itm:trivalency}, it suffices to 
consider the genus two node. The case with a solid genus two node 
follows from \ref{itm:I_bundles}, so we assume 
the genus two node is hollow, and 
$Y\subset \Compl S$ is the corresponding genus two component.
Suppose $\bm Y$ has more than $3$ components.
Then there exists an annular component $A$ in $\overline{\partial Y- \bm Y}$.
Let $A_1,A_2$ be the components 
of $\bm Y$ that meet $\partial A$.
Suppose the frontier $A'$ 
of a regular neighborhood of 
$A_1\cup A\cup A_2$ in $Y$ is inessential, 
then there is an essential square in $Y$, contradicting the 
completeness of $W\subset M$; on the other hand, 
since $Y$ is of genus $2$, $A'$ is not parallel to
any component of $\bm Y$; thus by 
by the simpleness of $Y$, neither can   
$A'$ be essential.
 . 
\end{proof}

\begin{definition}
We say the characteristic diagram of $M$
is of type $(e,l,b,\square)$
if $\charM$ has $e$ edges, $l$ loops, and $b$ bigons, 
and $\square=\snode$ (resp.\ $\hnode$) if the genus 
two node in $\charM$ 
is solid (resp.\ hollow). 
\end{definition}

\begin{theorem}\label{teo:possible_shapes}
Characteristic diagrams of $M$ 
are classified, up to isomorphism, by
their types $(e,l,b,\square)$ 
into $13$ classes in Table \ref{tab:configurations}.
\end{theorem}
\begin{proof}
Note first characteristic diagrams of the same type are isomorphic. 
By Proposition \ref{prop:properties_char_diagram_Ehk}\ref{itm:adjacency}, \ref{itm:trivalency}, 
we have $1\leq e\leq 3$, $l=0$ or $1$, and $b=0,1,3$,
and $(1,1,0, \snode)$, $(2,1,0,\snode)$ are ruled out
by \ref{itm:no_solid_loops} and 
$(2,0,1,\snode)$ ruled out by \ref{itm:no_solid_bigon}
in Proposition \ref{prop:properties_char_diagram_Ehk}.  
\end{proof}
Recall that a handlebody-knot $\pair$ 
is \emph{irreducible} if $\ComplHK$ is $\partial$-irreducible.
\begin{lemma}\label{lm:tri_redu_toro} 
Suppose $\pair$ is reducible.
Then it is trivial if and only if it is atoroidal.
\end{lemma} 
\begin{proof}
Observe first that there exists a separating essential disk 
$D\subset \ComplHK$.  
The disk $D$ splits $\Compl\HK$ into two knot exteriors 
$\Compl {K_1}$, $\Compl{K_2}$, for some knots $K_1,K_2$ in $\sphere$.  
Then $\pair$ is trivial if and only if both $K_1$, $K_2$
are trivial and therefore if and only if $\pair$ is atoroidal.
\end{proof}

\begin{corollary}\label{cor:irre_atoro_vs_tri_atoro}
Suppose $\pair$ is atoroidal. Then $\pair$ is non-trivial if and only if
$\pair$ is irreducible. 
\end{corollary}
\begin{proof}
``$\Leftarrow$'' is straightforward, while
``$\Rightarrow$'' follows from Lemma \ref{lm:tri_redu_toro}.
\end{proof}
 
\begin{definition}[Annulus Diagram]
Let $\pair$ be a non-trivial, atoroidal handlebody-knot. Then the annulus diagram $\charhk$ of $\pair$
is the characteristic diagram $\charE$ 
of $\Compl\HK$ together with
a labeling, $\hopf_i,\knot_i,\link(r_1,r_2)$ or $\EM$, for each edge, based on the type of the annulus the edge represents. 
\end{definition}
%
%
 

\section{Classification}\label{sec:classification}
Throughout the section, 
$\pair$ is
a non-trivial atoroidal handlebody-knot. 
We examine here combinations of 
non-separating annuli of various types in $\ComplHK$.
$A\subset\ComplHK$
is a non-separating essential 
annulus, and $\HK_A:=\HK\cup\rnbhd{A}$. 
The frontier of $\rnbhd{A}$ in $\ComplHK$
are two annuli $A_+, A_-$, whose cores we denote by $l_+,l_-$, respectively.
We orient $l_+,l_-$ so that $[l_+]=[l_-]\in H_1(\rnbhd{A})$. In the case $A$ is of type $2$-$2$, one of $l_+,l_-$, say $l_-$, is separating in $\partial \HKA$.
We denote the components of $\partial A$ by $l_1,l_2$
if $A$ is of type $3$-$3$, and by $l_A,l$ 
if $A$ is of type $2$ with $l_A$ the one bounding a disk in $\HK$.  
In addition, by ``unique'', we understand ``unique, up to isotopy''. 
\subsection{Annulus configuration}
Recall first a result on type $4$-$1$ annuli
\cite[Lemma $3.7$]{FunKod:20}, \cite[Lemma $2.2$]{Wan:23p}. 
\begin{lemma}\label{lm:EM_excludes_nonseparating}
Let $\hat A\subset\ComplHK$ be a type $4$-$1$ annulus. 
Then no non-separating essential annulus
in $\ComplHK$ disjoint from $\hat A$ exists.
\end{lemma}

Given a type $3$-$3$ annulus $A$,
we fix an oriented disk $\disk_A\subset\HK$ disjoint from $\partial A$. Recall the definition of meridional basis from \cite{Wan:22}. 
\begin{definition}
Suppose $A$ is 
of type $3$-$3$ with a slope pair $(\frac{p}{q},pq)$.
Then a meridional basis of $H_1(\Compl\HKA)$ is a basis
given by the homology classes of the boundary of 
two oriented, disjoint,
non-parallel meridian disks
$D_1,D_2\subset\HKA$ disjoint from $\disk_A$
with $[\partial D_1]-[\partial D_2]=[\partial\disk_A]\in H_1(\ComplHKA)$. 
\end{definition}
\begin{lemma}\label{lm:l_pm_meridional_basis}
Suppose $A$ is of type $3$-$3$ with a slope pair 
$(\frac{p}{q},pq)$
and  $\{b_1,b_2\}$ a meridional basis of $H_1(\ComplHKA)$.
If $[l_+]=(p_1,p_2)$ in terms of $\{b_1,b_2\}$,
then $[l_-]=(p_1\mp 1,p_2\pm 1)$ and $p_1+p_2=\pm p$. 
\end{lemma}
\begin{proof}
Denote by $V_1,V_2$ the solid tori in 
$\HK-\openrnbhd{\disk_A}$,
and by $U$ the solid torus 
$V_1\cup V_2\cup \rnbhd{A}$. 
Then $l_+,l_-$ are two 
parallel curves in $\partial U$, and they 
separate the two disk components of the frontier
$\partial_\HK \rnbhd{\disk_A} \subset \partial U$,
so $[l_+]-[l_-]=\pm [\partial \disk_A]\in H_1(\ComplHKA)$
and therefore the first assertion.   
Consider the short exact sequence
\[
0\rightarrow \langle [\partial\disk_A]\rangle
\rightarrow H_1(\ComplHKA)
\rightarrow H_1(U)\simeq \langle b_1=b_2\rangle
\rightarrow 0,
\]
and note that the slopes of $l_+,l_-\subset \partial U$ are
$\frac{p}{q}$ with respect to $(\sphere,U)$. Hence $p_1+p_2=\pm p$. 
\end{proof}

\begin{lemma}\label{lm:lpm_typethreethree}
Suppose $A$ is of type $3$-$3$ with a boundary slop pair $(r_1,r_2)$. 

If $(r_1,r_2)=(\frac{p}{q},\frac{q}{p})$, $pq\neq 0$,
then $\{[l_+],[l_-]\}$ is a basis of $H_1(\ComplHKA)$.

If $(r_1,r_2)=(\frac{p}{q},pq)$, $pq\neq 0$,
then $\langle [l_+],[l_-]\rangle$ is a subgroup of $H_1(\ComplHKA)$ with index $\vert p\vert$.

If $(r_1,r_2)=(0,0)$,
then $\langle [l_+],[l_-]\rangle$ is a rank one subgroup of $H_1(\ComplHKA)$.
\end{lemma}
\begin{proof}
Denote by $V_1,V_2$ the solid tori in $\HK-\openrnbhd{\disk_A}$,
and by $U$ the union $V_1\cup V_2\cup \rnbhd{A}$.

Suppose $(r_1,r_2)=(\frac{p}{q},\frac{q}{p})$, $\vert p\vert,\vert q\vert>1$. 
Then $U$ is a Seifert 
fibered space with two exceptional fibers, and
therefore the exterior $W:=\Compl U$ of $U$ in $\sphere$
is a solid torus, whose core is a 
$(p,q)$-torus knot in $\sphere$.
Since $l_+,l_-$ are parallel to the core of $W$
in $W$ by \cite{Sei:33}, $[l_+]=[l_-]$ generates $H_1(W)$.
On the other hand, we have 
$\ComplHKA=W-\rnbhd{\disk_A}$,
that is, $\ComplHKA$ is obtained by removing 
a regular neighborhood of an arc in $W$ dual to $\disk_A$, so 
$H_1(W,\ComplHKA)=0$. This together with $H_2(W)=0$ implies
the short exact sequence
\[
0\rightarrow 
H_2(W,\ComplHKA)\rightarrow 
H_1(\ComplHKA)\rightarrow 
H_1(W)\rightarrow 0
\]
given by the inclusion 
$\ComplHKA\hookrightarrow W$. 
Since $[\disk_A]$
generates $H_2(W,\ComplHKA)$,
$\pm [\partial \disk_A]=[l_+]-[l_-]$,
and $[l_+]=[l_-]$ generates $H_1(W)$, 
we have $\{[l_+],[l_-]\}$ is a basis of $H_1(\ComplHKA)$.

Suppose $(r_1,r_2)=(\frac{p}{q},pq)$, $q\neq 0$.
Then by Lemma \ref{lm:l_pm_meridional_basis},  
$[l_+]=(p_1,p_2)$ and $[l_-]=(p_1\mp 1, p_2\pm 1)$
with $p_1+p_2=p$ in terms of a meridional basis
of $H_1(\ComplHKA)$, and hence the determinant
\[
\begin{vmatrix}
p_1& p_2\\
p_1\mp 1& p_2\pm 1. 
\end{vmatrix}=\pm (p_1+p_2)=\pm p
\] 
In other words, when $p\neq 0$, 
$\langle [l_+],[l_-]\rangle<H_1(\ComplHKA)$
is a subgroup of rank two with index $\vert p\vert$.
When $p=0$, since $[l_+]-[l_-]=\mp (1,-1)$, at least
one of $[l_+],[l_-]\in H_1(\ComplHKA)$ is non-trivial,
so $\langle [l_+],[l_-]\rangle$ is a subgroup 
isomorphic to $\Z$. 
\end{proof}

\begin{corollary}\label{cor:typethreethree_nontrivial_slope}
Suppose $A$ is of type $3$-$3$ with a non-trivial slope pair, and $A'$ is a non-separating annulus disjoint from $A$.
Then $\partial A, \partial A'$ are parallel in $\partial\HK$. 
In particular, $A'$ is of type $3$-$3$ with the same slope pair.
\end{corollary}
\begin{proof}
Let $P$ be the planar surface
$\bComplHKA-\mathring{A_+}\cup \mathring{A_-}$. 
Denote by $l_{1\pm},l_{2\pm}$
the components of $\partial A_\pm$ and by $l_1', l_2'$
the components of $\partial A'$.
%
Since $l_1',l_2'\subset P$, 
one of $l_1',l_2'$ is parallel to one of 
$l_{1\pm},l_{2\pm}$; it may be assumed that $l_1'$ is parallel to $l_{1+}$.
By Lemma \ref{lm:lpm_typethreethree}, 
$[l_+]\neq \pm [l_-]$ and none of $[l_+],[l_-]$
is trivial in $H_1(\ComplHKA)$.
These, together with $[l_1']=[l_2']\in H_1(\ComplHKA)$,
imply that $l_2'$ is parallel to either $l_{2+}$
or $l_{1+}$. The latter is impossible since 
$l_1',l_2'$ are not parallel 
in $\partial \HK$ and hence not parallel in $P$.
Therefore $\partial A'$ is parallel to $\partial A_+$
and hence to $\partial A$.
\end{proof}

There is an analog of Lemma \ref{lm:lpm_typethreethree}
for type $2$ annuli.    
\begin{lemma}\label{lm:l_pm_typetwo}
If $A$ is of type $2$-$1$, then 
$\{[l_+],[l_-]\}$ is a basis of $H_1(\ComplHKA)$.
If $A$ is of type $2$-$2$, then 
$[l_-]$ is trivial and the quotient $H_1(\ComplHKA)/\langle [l_+]\rangle \simeq \Z$. 
\end{lemma}
\begin{proof}
It follows from the fact that $l_+,l_-$
bound non-parallel, non-separating meridian disks in
$\HK$ if $A$ is of type $2$-$1$, and $l_-$ (resp.\ $l_+$) bounds a separating (resp.\ non-separating) disk in $\HK$ if $A$ is of type $2$-$2$.
\end{proof} 

\begin{lemma}\label{lm:typethreethree_trivial_slope}
Suppose $A$ is of type $3$-$3$ with a trivial slope pair, and $A'$ is a type $3$-$3$ annulus disjoint from $A$.
Then $A,A'$ are parallel $\ComplHK$.  
\end{lemma}
\begin{proof}
Suppose $\partial A$ and $\partial A'$ are parallel
in $\partial \HK$. Let $B_1,B_2\subset\partial \HK$ 
be the annuli cut off by $\partial A,\partial A'$.
Then $A\cup A'\cup B_1\cup B_2$ bounds
a solid torus $V$ in $\ComplHK$ by the atoroidality of 
$\pair$. Since $A$ has a trivial slope pair,
the linking number $\lk{l_1}{l_2}$ is $0$ and hence the core of $A$
is a preferred longitude with respect to $(\sphere,V)$;
this implies $H_1(A)\rightarrow H_1(V)$
is an isomorphism, so $A,A'$ are parallel
through $V$.

Suppose $\partial A$ and $\partial A'$ are not parallel.
Let $l_1',l_2'$ be the components of $\partial A'$.
Then since $\partial\HK-\partial A$
is a four-times punctured sphere,
it may be assumed that $l_1, l_1'$ 
are parallel in $\partial\HK$, and $l_2,l_2'$ are not.
Let $B_1\subset\partial \HK$ be the annulus cut off by $l_1.l_1'$. 
Then $B_1\cup A\cup A'$ induces an annulus $A''\subset\ComplHK$
disjoint from $A\cup A'$
with $\partial A''$ parallel to $l_2, l_2'$. 
Let $B_2,B_3\subset\partial\HK$ be the annuli cut off by 
$\partial A''$ and $l_2\cup l_2'$. 
Then the torus 
$B_1\cup B_2\cup B_3\cup A\cup A'\cup A''$
bounds a solid torus $W$ in $\ComplHK$
since $\pair$ is atoroidal. 

Let $P_1,P_2\subset \partial\HK$ 
be the pairs of pants cut off by $B_1\cup B_2\cup B_3$.
Then $P_1,P_2$ can be regarded as 
a planar surface in $\Compl W$.
By \cite[Lemma $3.5$]{KodOzaGor:15}, 
$P_1,P_2$ are inessential in $\Compl W$.

\textbf{Case $1$: $P_1$ is compressible.}
Let $D$ be a compressing disk of $P_1$
that minimizes 
\[
\#\{D \cap P_2\mid D\text{ a compressing disk of $P_1$}\}.
\] 

\textbf{Subcase $1.1$: $D\cap P_2=\emptyset$.}   
$D$ is either in $\HK$ or in $\ComplHK$.
Since $\partial D$ is essential in $P_1$, 
$\partial D$ is essential in $\partial\HK$,
so $D$ is a compressing disk of $\partial\HK$ 
in $\sphere$.
On the other hand, 
$\partial A\cup \partial A'\cup\partial A''$ contains three mutually non-parallel simple loops
in $\partial \HK$ that bound no disks in $\HK$, 
so every meridian disk in $\HK$
meets $\partial A\cup \partial A'\cup\partial A''$,
and hence $D\subset \ComplHK$, but this 
contradicts $\pair$ is irreducible.

\textbf{Subcase $1.2$: $D\cap P_2\neq\emptyset$.} 
Note first that $D\cap P_2$ only contains circles.
Let $D'\subset D$ be the disk cut off by a circle in $D\cap P_2$
innermost in $D$. By the minimality $\partial D'$
is essential in $P_2$; hence $D'$ is a compressing
disk of $\partial\HK$ in $\sphere$, 
a contradiction as in \textbf{Subcase $1.1$}. 

The same argument applies to the case where
$P_2$ is compressible.

\textbf{Case $2$: $P_1, P_2$ are incompressible.}
First observe that, since none of the components of 
$\partial A\cup\partial A'\cup \partial A''$ 
is separating in $\partial \HK$, $P_1$ (resp.\ $P_2$) meets 
$B_i$ for each $i$.
Let $D$ be a $\partial$-compressing
disk of $P_1$ that minimizes 
\[
\#\{D\cap P_2\mid \text{$D$ a $\partial$-compressing disk of $P_1$}\}.
\]
Then by the minimality and incompressibility of $P_2$,
$D\cap P_2$ is either empty or some arcs.

\textbf{Subcase $2.1$: $D\cap P_2=\emptyset$.}
Denote by $\gamma$ the arc $D\cap \Compl W$, and 
note that $\gamma \subset B_u:=B_1\cup B_2\cup B_3$ 
if $D\subset\HK$; otherwise $\gamma\subset A_u:=A\cup A'\cup A''$. In addition, $\gamma$ is inessential in 
in either case:
in the former, it follows from the fact that 
none of $B_i$, $i=1,2,3$, has 
two boundary components lying in $P_1$,
whereas in the latter, it results 
from the $\partial$-incompressibility of $A,A',A''$. 
 
Let $D'$ be the disk cut off
from $B_u$ (resp.\ $A_u$).
Then $D\cup D'$ induces a disk $D''$
disjoint from $B_u$ (resp.\ $A_u$).
Since $D$ is a $\partial$-compressing disk of $P_1$ in $\Compl W$,
$\partial D''$ is essential in $P_1$, contradicting
the incompressibility of $P_1$.

\textbf{Subcase $2.2$: $D\cap P_2\neq\emptyset$.}
Let $D'\subset D$ be a disk cut off by 
an arc in $D\cap P_2$ outermost in $D$.
Denote by $\gamma$ the arc $D'\cap \partial W$;
as with \textbf{Subcase $2.1$}, $\gamma$ 
is either in $A_u$ or in $B_u$, and inessential whichever way. 
Let $D''$ be the disk cut off by $\gamma$ from $A_u$ or $B_u$. Then 
$D'\cup D''$ induces a disk $D'''$
disjoint from $P_1$ with $\partial D''\subset P_2$.
By the minimality of $\# D\cap P_2$, 
$\partial D'''$ is essential in $P_2$, contradicting
the incompressibility of $P_2$.

\cout{An alternative proof:
Suppose $D\cap P_2\neq\emptyset$,
and $D'\subset D$ be a disk cut off by 
an arc in $D\cap P_2$ that is outermost in $D$.
If $D'\cap P_2$ is an essential arc in $P_2$,
then $D'$ is a $\partial$-compressing disk of $P_2$
disjoint from $P_1$; in other words, $D'$ is either in $\HK$
or in $\ComplHKA$, and as with Subcase $2.1$, this is impossible.

If $D'\cap P_2$ is an inessential arc $\alpha$ in $P_2$, 
which cuts a disk $D''$ from $P_2$. $D'\cup D''$ 
induces a disk $D'''\subset\Compl W$ 
disjoint from $P_1\cup P_2$.
By the minimality,
$D'''$ is essential in $\Compl W$,
and hence $\partial D'''$ is an essential loop
in either $B_1\cup B_2\cup B_3$ or $A\cup A'\cup A''$,
depending on whether $D'''$ is in $\HK$ or $\ComplHKA$.
The former contradicts the definition of a type $3$-$3$
annulus, whereas the latter contradicts the essentiality of 
$A,A'$ or $A''$.
}
\end{proof}

\begin{lemma}\label{lm:basis_uniqueness}
If $\{[l_+],[l_-]\}$ is a basis of $H_1(\ComplHKA)$,
then $A$ is the unique annulus in $\ComplHK$.
\end{lemma}
\begin{proof}
By Theorem \ref{teo:engulfing}, 
it suffices to show that 
$\rnbhd{A}\subset\ComplHK$
is a characteristic submanifold of $(\ComplHK,\bphi)$.
To see this, we employ Theorem \ref{teo:completeness_criterion}.
Since $\rnbhd{A}$ is a full $\F$-manifold of 
$(\ComplHK,\bphi)$, 
it amounts to showing that every essential annulus
$A'$ in $\ComplHKA$ disjoint from $A_+,A_-$
is parallel to $A_+,A_-$, 
where $\ComplHKA\subset (\ComplHK,\bphi)$ 
is endowed with the proper boundary pattern.
Denote by $l'$ a core of $A'$.

\textbf{Case $1$: $A'$ is non-separating in $\ComplHK$.}
Since $\{[l_+],[l_-]\}$ is a basis of 
$H_1(\ComplHKA)$, the argument for Corollary \ref{cor:typethreethree_nontrivial_slope} applies 
and thus $\partial A'$ is parallel to $\partial A_+$
or $\partial A_-$ in $\bComplHKA$; 
it may be assumed that it is the former, and denote by
$B_1,B_2$ the annuli
cut off by $\partial A_+,\partial A'$ from $\partial \ComplHKA$.
Then $A\cup A'\cup B_1\cup B_2$ bounds a solid torus 
$W$ in $\ComplHKA$ by Corollary \ref{cor:atoro_HK_HKA}.
Let $X$ be the closure of the complement $\ComplHKA-W$ and 
$l_w$ a core of $W$, and orient $l',l_w$ 
so that $[l']=[l_+]$ and 
$[l']=k[l_w]$ with $k>0$ in $H_1(W)$.
Consider the short exact sequence 
\[
0\rightarrow
H_1(A')\xrightarrow{(\iota_1,\iota_2)}
H_1(W)\oplus H_1(X)\xrightarrow{\iota_3-\iota_4}
H_1(\ComplHKA)\rightarrow
0,
\]
where $\iota_i$, $i=1,2,3,4$, are induced by the inclusions.
Note that 
$\iota_4$ sends $[l']$ to $[l_+]$ and 
$[l_-]$ to itself, and $\iota_1$ sends $[l']$
to $k[l_w]$. 
Since $\{[l_+],[l_-]\}$ is a basis of 
$H_1(\ComplHKA)$, the image of $[l_w]$ under 
$\iota_3$ is $m[l_+]+n[l_-]$, for some $m,n\in\mathbb{Z}$.
Then the identity $\iota_3\circ \iota_1=\iota_4\circ\iota_2$ gives us $km[l_+]+kn[l_-]=[l_+]$, and 
therefore $n=0,k=m=1$. This implies
$H_1(A')\xrightarrow{\iota_1} H_1(W)$ is an isomorphism, and
hence $A'$ is parallel to $A_+$ through $W$ in $\ComplHKA$.

\textbf{Case $2$: $A'$ is separating in $\ComplHK$.}
Since the components of
$\partial A'$ are parallel 
and do not separate the components of $\partial A$ 
in $\partial \HK$, the components of $\partial A'$ 
are 
parallel in $\partial \ComplHKA$.
Let $B\subset \partial\ComplHKA$ 
be the annulus cut off by $\partial A'$.
Then $B\cup A'$ bounds a solid torus $W$
in $\ComplHKA$ by Corollary \ref{cor:atoro_HK_HKA}. 
Set $X:=\overline{\ComplHKA-W}$, and consider 
the short exact sequence
\[
0\rightarrow
H_1(A')\xrightarrow{(\iota_1,\iota_2)}
H_1(W)\oplus H_1(X)\xrightarrow{\iota_3-\iota_4}
H_1(\ComplHKA)\rightarrow
0,
\]
where $\iota_i$, $i=1,2,3,4$, are induced by the inclusions.
Let $l_w$ be a core of $W$, and orient $l',l_w$ 
so that $[l']=k[l_w]$ with $k>0$. Note that
$k$ is necessarily larger than $1$ by the essentiality of 
$A'$.
Since $[l_+],[l_-]\in H_1(X)$ and $H_2(\ComplHKA,X)=0$,
$\iota_4:H_1(X)\rightarrow H_1(\ComplHKA)$ is an isomorphism.
Let the image of $[l_w]$ under $\iota_3$ be $m[l_+]+n[l_-]$, 
and the image of $[l']$ under $\iota_2$ be 
$m'[l_+]+n'[l_-]$,
for some $m,n,m',n'\in\mathbb{Z}$.
Then $x=([l_w],m[l_+]+n[l_-])\in H_1(W)\oplus H_1(X)$ 
is in the kernel of $\iota_3-\iota_4$, and therefore,
there exists $c\in\mathbb{Z}$
such that the image of $c[l']$
under $(\iota_1,\iota_2)$ is $x$; in other words, we have 
the equality 
\[
\big(kc[l_w],m'c[l_+]+n'c[l_-]\big)=\big([l_w],m[l_+]+n[l_-]\big)\in H_1(W)\oplus H_1(X), 
\]  
but this implies $k=c=1,m=m',n=n'$, contradicting $k>1$.
\end{proof}

\begin{lemma}\label{lm:l_pm_basis}
$\{[l_+],[l_-]\}$ is a basis of $H_1(\ComplHKA)$
if and only if $A$ is of type $2$-$1$ or 
of type $3$-$3$ with the slope pair
$(\frac{p}{q},\frac{q}{p})$, $pq\neq 0$.
\end{lemma}
\begin{proof}
``$\Leftarrow$'' follows from Lemmas \ref{lm:l_pm_typetwo}
and \ref{lm:lpm_typethreethree}.
``$\Rightarrow$'' also results from 
the same lemmas as 
$\{[l_+],[l_-]\}$ can form a basis of $H_1(\ComplHKA)$
only if $A$ is of type $2$ or of type $3$-$3$.
\end{proof}

Lemmas \ref{lm:basis_uniqueness} and \ref{lm:l_pm_basis} give us 
the following uniqueness result.
\begin{corollary}\label{cor:unique_annulus}
If $A$ is of type $2$-$1$ or 
of type $3$-$3$ with the slope pair
$(\frac{p}{q},\frac{q}{p})$, $pq\neq 0$, 
then $A$ is the unique annulus in $\ComplHK$.
\end{corollary}

\begin{lemma}\label{lm:typetwotwo_parallel_boundary}
Let $A,A'$ be two disjoint type $2$-$2$ annuli in $\ComplHK$.
If $\partial A,\partial A'$ are parallel in $\partial\HK$,
then $A,A'$ are parallel in $\ComplHK$.
\end{lemma}
\begin{proof}
Let $B_1,B_2\subset \partial\HK$ 
be the annuli cut off by $\partial A, \partial A'$.
Then $B_1\cup B_2\cup A\cup A'$ bounds 
a solid torus $W$ in $\ComplHK$ by the atoroidality of $\pair$.
Observe that 
$l_A$ is a longitude of $(\sphere, W)$ 
since it bounds a disk in $\HK$. This implies 
$H_1(A)\rightarrow H_1(W)$
is an isomorphism, and hence $A,A'$ are parallel through $W$ in $\ComplHK$.
\end{proof}
\begin{lemma}\label{lm:two_typetwotwo_one_typethreethree}
Suppose $A$ is of type $2$-$2$.
Then there exists another type $2$-$2$
annulus $A'$ disjoint from and non-parallel to $A$
if and only if there exists a type $3$-$3$
annulus $A''$ with a trivial slop pair disjoint from $A$. 
\end{lemma}
\begin{proof}
``$\Rightarrow$'': Let $l_{A'}\subset\partial A'$ be the 
component that bounds a disk in $\HK$ and $l'\subset \partial A'$
another component. Then $l_A,l_{A'}$ are parallel and bound an
annulus $B$ in $\partial\HK$, and 
$l,l'$ are non-parallel in $\partial\HK$ by Lemma \ref{lm:typetwotwo_parallel_boundary}.
The union $A\cup A'\cup B$ induces a type $3$-$3$
annulus, which has a trivial slope pair since
$\lk{l}{l'}=\lk{l_A}{l_{A'}}=0$.

``$\Leftarrow$'': 
Let $l_1'',l_2''$ be components of $\partial A''$. Then 
one of them, say $l_1''$, is parallel to $l$ in $\partial\HK$. 
Let $B\subset\partial\HK$ be the annulus cut off by $l,l_1''$.
Then the union $A\cup B\cup A''$ induces an annulus
whose boundary components parallel to $l_A, l_2''$,
so it is of type $2$-$2$.
\end{proof}
\begin{corollary}\label{cor:three_typetwotwo}
Let $A,A',A''$ be three disjoint type $2$-$2$ annuli in $\ComplHK$. Then at least two of them are parallel in $\ComplHK$.
\end{corollary}
\begin{proof}
Let $l'\subset \partial A',l''\subset \partial A''$ be the components that do not bound a disk in $\HK$,
and $l_{A'}\subset \partial A', l_{A''}\subset\partial A''$
the other components. Then $l_A,l_{A'},l_{A''}$
are parallel in $\partial\HK$ by the definition of 
a type $2$-$2$ annulus.
 
Suppose $A,A'$ are not parallel in $\ComplHK$. 
Then $l,l'$ are longitudes of the solid tori
$V,V'$ in $\HK-\mathring{U}$, where $U\subset \HK$
is the $3$-ball cut off by the disks bounded by $l_A,l_{A'}$. 
In particular, $l''$
is parallel to either $l$ or $l'$, so 
by Lemma \ref{lm:typetwotwo_parallel_boundary}, $A''$ is parallel
to $A$ or $A'$. 
\end{proof}

\subsection{Classification theorems}
Let $\charE$ be 
the characteristic diagram of $\ComplHK$,
and $\charhk$ the annulus diagram of $\pair$.
\begin{theorem}[\textbf{$\theta$-shape characteristic diagram}]
\label{teo:char_three_bigons}
If $\charE$ is \raisebox{-.2 cm}{\includegraphics[scale=.12]{theta_char}},
then $\charhk$ is 
\raisebox{-.2 cm}{\includegraphics[scale=.12]{typetwotwo_ann_3}}, $\square=\snode$ or $\hnode$, 
and the Seifert
fibered solid torus has no exceptional fiber.
\end{theorem}
\begin{proof}
Let $A,A',A''$ be the non-separating 
annuli corresponding to the edges of 
$\charE$. None of them is of type $2$-$1$ by Corollary \ref{cor:unique_annulus} or of type $3$-$3$
with a non-trivial slope pair by Corollaries \ref{cor:unique_annulus} and 
\ref{cor:typethreethree_nontrivial_slope}
since no two of them separate $\ComplHK$. 
Therefore, $A,A',A''$ are of type $2$-$2$ or of type $3$-$3$
with a trivial slope.
By Lemma \ref{lm:typethreethree_trivial_slope}, 
at most one of them is of type $3$-$3$, whereas by
Corollary \ref{cor:three_typetwotwo}, 
at most two of them is of type $2$-$2$, so $\charhk$ is 
\raisebox{-.25 cm}{\includegraphics[scale=.1]{typetwotwo_ann_3}}.

Let $W$ be the component corresponding to the genus one node,
and $A$ the type $3$-$3$ annulus.
If a core of $A$ 
is a $(p,q)$-curve with respect to 
$(\sphere,W)$, then the  
the linking number of the components of $\partial A$ in $\sphere$ is $\pm pq$. Since $A$ has a trivial slope pair, $pq=0$,
and by the essentiality of $A$, $q\neq 0$ and therefore $(p,q)=(0,\pm 1)$. Thus $W$ has no exceptional fiber.
\end{proof}

\begin{lemma}\label{lm:nonchar_annuli}
$\ComplHK$ contains a non-characteristic, 
non-separating annulus $A$ if and only if
$\charE$ is 
\raisebox{.05 cm}{\includegraphics[scale=.1]{stick_char}}. 
In addition, 
$A$ is of type $3$-$3$ with a boundary 
slope pair $(\frac{p}{q},pq)$,
$pq\neq 0$, and is the unique non-separating annulus in $\ComplHK$.
\end{lemma}
\begin{proof}
``$\Leftarrow$'': 
Let $X$ be the component corresponding to
the genus two node. By Proposition \ref{prop:properties_char_diagram_Ehk},
$X$ is I-fibered over a Klein bottle $B$ with
an open disk removed. 
Any non-separating simple loop $l$ in $B$ 
induces an essential annulus $A$ in
$X$ and hence in $\ComplHK$ by Lemma \ref{lm:essentiality_in_X_in_M}.
Since $l$ cannot be isotoped away from
essential separating loops 
that are not parallel to $\partial B$ in $B$ by \cite[Theorem $3.3$]{Gom:17},
$A$ is not characteristic.

``$\Rightarrow$'': 
By Theorem \ref{teo:engulfing} and Lemma \ref{lm:essentiality_in_X_in_M}, 
we may assume $A$ is an essential annulus 
in a component $X$ of a 
characteristic submanifold of $\ComplHK$ 
with $A$ non-parallel to any component of $\bEhk X$.
By Proposition \ref{prop:properties_char_diagram_Ehk}, $X$
is either an I-bundle with $\chi(\partial X)<0$
or a Seifert fibered solid torus.
The latter is impossible because  
$\#\bEhk X\leq 3$ by Theorem \ref{teo:possible_shapes} and $X$ has no exceptional fiber
by Theorem \ref{teo:char_three_bigons}
when $\#\bEhk X=3$. 

Therefore, $X$ is an I-bundle over a Mobius band
or Klein bottle with an open disk removed; 
in particular, $\charE$ is \raisebox{-.15 cm}{\includegraphics[scale=.1]{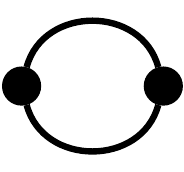}} or \raisebox{.05 cm}{\includegraphics[scale=.1]{stick_char}}. 
The former is ruled out by Proposition \ref{prop:properties_char_diagram_Ehk}\ref{itm:no_solid_bigon}, so
$X$ is an I-bundle
over a Klein bottle with an opened disk removed $B$, and $\charE$ 
is \raisebox{.05 cm}{\includegraphics[scale=.1]{stick_char}}.

By \cite[Theorem $3.3$]{Gom:17}, 
every two non-separating simple loops in
a Klein bottle with an opened disk removed
are isotopic, 
so $A$ is the unique non-separating annulus in $\ComplHK$.  
Now, to determine the type of $A$,
first note that 
the annulus $A':=\bEhk X\subset \ComplHK$
is an annulus non-isotopic to $A$, 
so $A$ is not of type $2$-$1$ or of type $3$-$3$ with a slope pair $(\frac{p}{q},\frac{q}{p})$, $pq\neq 0$, by Corollary \ref{cor:unique_annulus}.   
Denote by $X'$ the solid torus $\overline{\ComplHK-X}$ 
and observe that, by the essentiality of $A'=X\cap X'\subset\ComplHK$, 
the homomorphism
\[H_1(A')\simeq \Z \xrightarrow{k} \Z \simeq H_1(X')\]
induced by the inclusion
neither is trivial nor is an isomorphism, namely $k\neq 0, \pm 1$. On the other hand, 
the decomposition $\ComplHKA=(X-\openrnbhd{A})\cup X'$ 
gives us the isomorphism: 
\begin{equation}\label{eq:Klein_seifert_homology_HKA}
H_1(\ComplHKA)\simeq \langle v_+,v_-,u\rangle/(v_+ +  v_- =\pm ku), 
\end{equation}
where $u$ is a generator of $H_1(X')$, 
$v_\pm=[l_\pm]$, and $l_\pm$ are the cores of 
the frontier $\bEhk\rnbhd{A}$. 
If $A$ is of type $2$-$2$, then 
$v_-$ is trivial in $H_1(\ComplHKA)$ by Lemma \ref{lm:l_pm_typetwo}, so 
$H_1(\ComplHKA)\simeq \Z$, a contradiction. 
If $A$ is of type $3$-$3$ with a trivial slope pair, then at least one of $v_+,v_-$
is not a generator by Lemma \ref{lm:l_pm_meridional_basis},
contradicting \eqref{eq:Klein_seifert_homology_HKA}, as both $\{v_+, u\}$
and $\{v_-,u\}$ form a basis of $H_1(\ComplHKA)$.
Therefore $A$ is of type $3$-$3$ with a slope pair $(\frac{p}{q},pq)$, $pq\neq 0$.
\end{proof}

\begin{corollary}\label{cor:char_annuli}
If $A$ is of type $2$ or of type $3$-$3$ with a trivial slope pair or a slope pair $(\frac{p}{q},\frac{q}{p})$, $pq\neq 0$, then 
$A$ is characteristic.
\end{corollary}

\begin{corollary}\label{cor:disjoint_non_separating_annuli}
Up to isotopy, non-separating annuli in $\Compl\HK$
are mutually disjoint.
\end{corollary}

\begin{theorem}[\textbf{Classification Theorem}]\label{teo:classification_char_diagram_type_two}
\hfill
\begin{enumerate}[label=\textnormal{(\roman*)}]
\item\label{itm:charhk_type2_1}
If $A$ is of type $2$-$1$, then  
$\charhk$ is \raisebox{-.2 cm}{\includegraphics[scale=.1]{typetwoone_ann}}.
\item\label{itm:charhk_type2_2}
If $A$ is of type $2$-$2$, then $\charhk$ is one of the following:
\begin{figure}[H]
\begin{subfigure}{0.32\linewidth}
\centering
\includegraphics[scale=.12]{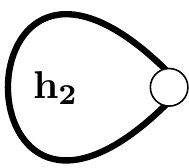}
\end{subfigure}
\begin{subfigure}{0.32\linewidth}
\centering
\includegraphics[scale=.12]{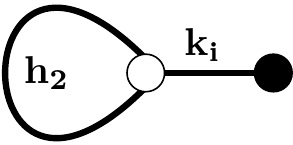} 
\quad 
\raisebox{.2 cm}{\footnotesize{$i=1$ or $2$}}
\end{subfigure}
\begin{subfigure}{0.32\linewidth}
\centering
\includegraphics[scale=.12]{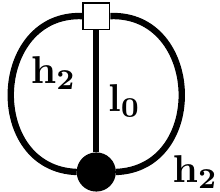} 
\quad \raisebox{.4 cm}{\footnotesize{$\square=\snode$ or $\hnode$}}.
\end{subfigure}
\end{figure} 
\end{enumerate}
\end{theorem} 
\begin{proof} 
\ref{itm:charhk_type2_1} follows from 
Corollary \ref{cor:unique_annulus}.
To see \ref{itm:charhk_type2_2}, let $S$ 
be a characteristic surface of $\ComplHK$.
By Theorem \ref{teo:possible_shapes}, 
$S$ consists of at most three annuli, 
one of which is $A$ by Corollary \ref{cor:char_annuli}. 

\textbf{Case $1$: $\# S=1$.} 
This implies $\charhk$ is 
\raisebox{-.2 cm}{\includegraphics[scale=.11]{typetwotwo_ann_1}}.
 
\textbf{Case $2$: $\# S=2$.} 
Let $A'\in S$ be
the other annulus. Then
by Corollaries \ref{cor:unique_annulus} and \ref{cor:typethreethree_nontrivial_slope},
it is not of type $2$-$1$ or of type $3$-$3$ with a non-trivial slope pair.
By Lemma \ref{lm:two_typetwotwo_one_typethreethree}
and Corollary \ref{cor:char_annuli}, it is not 
of type $2$-$2$ or of type $3$-$3$ with a trivial slope pair since $\# S=2$.
Therefore $A'$ is separating, and by Lemma \ref{lm:EM_excludes_nonseparating}, 
it is not of type $4$-$1$, so $\charhk$ is 
\raisebox{-.2 cm}{\includegraphics[scale=.12]{typetwotwo_ann_2}}, $i=1$ or $2$.

\textbf{Case $3$: $\# S=3$.}
Let $A',A''$ be the other two annuli.
Then at least one of them, say $A'$, is non-separating by Theorem \ref{teo:possible_shapes}. 
On the other hand, $A'$ cannot be of type $2$-$1$
or of type $3$-$3$ with a non-trivial slope 
by Corollaries \ref{cor:unique_annulus} and \ref{cor:typethreethree_nontrivial_slope},
so $A'$ is of type $2$-$2$
or of type $3$-$3$ with a trivial slope pair;
this implies that
$A''$ is of type $3$-$3$ with a trivial slope pair 
or of type $2$-$2$, respectively, by Lemma \ref{lm:two_typetwotwo_one_typethreethree}
and Corollary \ref{cor:char_annuli}. 
Therefore $\charhk$ is 
\raisebox{-.2 cm}{\includegraphics[scale=.12]{typetwotwo_ann_3}}, $\square=\snode$ or $\hnode$. 
\end{proof}
%
%
We now give a characterization of $\pairfourone$
in terms of characteristic diagrams.
%
%
%
\begin{figure}[b]
\centering
\def\svgwidth{.4\columnwidth}
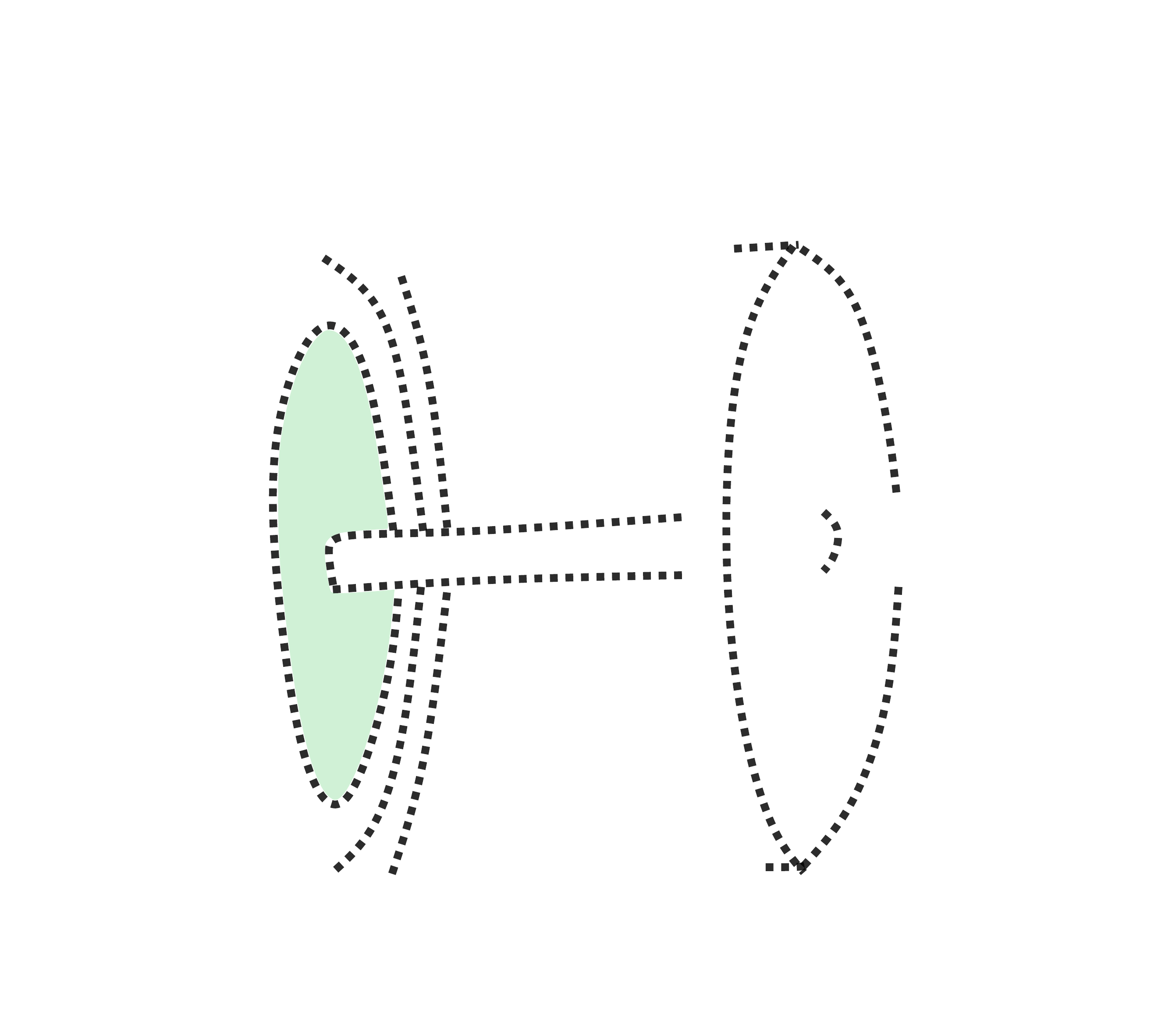
\caption{Decompose $\Compl\HK$ and $\HK$.}
\label{fig:decomposition_303s_hk}
\end{figure}
\begin{lemma}\label{lm:rigidity_solid_three_bigons}
Suppose the annulus diagrams of 
the handlebody-knots $\pair,\pairtilde$ both are 
\raisebox{-.2 cm}{\includegraphics[scale=.12]{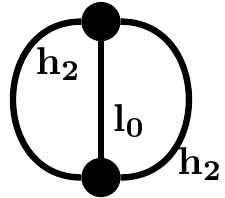}}. 
Then $\pair$ and $\pairtilde$ are equivalent.
\end{lemma}
\begin{proof}
Let $A$ (resp.\ $\til A$) and $A_0,A_1$
(resp.\ $\til A_0, \til A_1$) be the type $3$-$3$
annulus and the two type $2$-$2$ annuli
in $\ComplHK$ (resp.\ $\Compl\tildeHK$), respectively,
and denote by 
$l_{A_0},l_{A_1}$ (resp.\ $\til l_{A_0},\til l_{A_1}$ ) 
the boundary components of $A_0,A_1$ (resp.\ $\til A_0, \til A_1$) that bound disks $\disk_{A_0},\disk_{A_1}$
(resp.\ $\disk_{\til A_0},\disk_{\til A_1}$) in $\HK$ (resp.\ $\tildeHK$), respectively.
Also, let $U\subset \ComplHK,\tilde{U}\subset\Compl\tildeHK$ be
the I-bundles and 
$W, \tilde{W}$ their exteriors in $\ComplHK,\Compl\tildeHK$, respectively.
Note that $W$ (resp.\ $\til W$) are Seifert fibered solid torus
with frontier in $\ComplHK$
$A\cup A_0\cup A_1$ (resp.\ $\til A\cup \til A_0\cup \til A_1$),
and 
$l_{A_0},l_{A_1}$ (resp.\ $\til l_{A_0},\til l_{A_1}$) 
lie in different lids of $U$ (resp.\ $\til U$); 
see Fig.\ \ref{fig:decomposition_303s_hk}. 

To show $\pair,\pairtilde$ are equivalent,
we first construct a homeomorphism 
\[f_0:(U,A,A_0,A_1,l_{A_0},l_{A_1})\rightarrow (\til U,\til A,\til A_0, \til A_1, l_{\til A_0}, l_{\til A_1}).\]
To do this, we identify 
$U,\tilde{U}$ with $P\times I,\tilde{P}\times I$,
respectively, where $P,\tilde{P}$ are pairs of pants.
Let $C,C_0,C_1$ (resp.\ $\til C,\til C_0,\til C_1$)
be the components of 
$\partial P$ (resp.\ $\partial \til P$),
and identify $(C_0\times I, C_0\times 0)$ and $(C_1\times I, C_1\times 1)$
with $(A_0,l_{A_0})$ and $(A_1,l_{A_1})$ (resp.\ $(\til C_0\times I, \til C_0\times 0)$ and $(\til C_1\times I, \til C_1\times 1)$
with $(\til A_0, \til l_{A_0})$ and $(\til A_1, \til l_{A_1})$),
respectively.
 
It is not difficult to see there exist homeomorphisms 
$g_i:P\times i\rightarrow\til P\times i$ 
that map $(C\times i,C_0\times i, C_1\times i)$
to $(\til C\times i, \til C_0\times i, \til C_1\times i)$,
$i=0,1$. 
On the other hand, since the mapping class group of 
a three-times punctured sphere is given by 
the permutation group
on the punctures,  
$g_0, g_1$ can be extended to $f_0$. 
 
Now, let $V, V_0,V_1\subset \HK$ (resp.\ $\til V, \til V_0, \til V_1\subset\tildeHK$)
be the $3$-ball and two solid tori 
cut off by $\disk_{A_0},\disk_{A_1}$
(resp.\ $\disk_{\til A_0},\disk_{\til A_1}$)
such that $\disk_{A_i}, P\times i\subset \partial V_i$
(resp.\ $\disk_{\til A_i}, \til P\times i\subset \partial \til V_i$), $i=0,1$. Then the exterior $\Compl {V\cup W}$ 
(resp.\ $\Compl{\til V\cup \til W}$)
of $V\cup W$ 
(resp.\ $\til V\cup \til W$)
in $\sphere$ is $U\cup V_0\cup V_1$ 
(resp.\ $\til U\cup \til V_0\cup \til V_1$);
see Fig.\ \ref{fig:decomposition_303s_hk}, and
$f_0$ can be extended to a homeomorphism
\[f_1:(\Compl{V\cup W},U,V_0,V_1)
\rightarrow 
(\Compl{\til V\cup \til W},\til U,\til V_0,\til V_1)
\] 
as follows. 
Extend first the restriction $f_0\vert_{P\times i}$
to a homeomorphism 
\[
\bar{f_0}:\partial (V_0\cup V_1)\rightarrow \partial (\til V_0\cup\til V_1)\] 
that sends a meridian
of $V_i$ to a meridian of $\til V_i$, $i=0,1$; 
this can be done because 
$\partial V_i-\mathring{P}\times i$ 
consists of an annulus and the disk $\disk_{A_i}$.
Then extend $\bar{f_0}$ to a homeomorphism from 
$V_0\cup V_1$ to $\til V_0\cup \til V_1$, which, together with $f_0$, induces $f_1$.

Observe that $\Compl{V\cup W}$
(resp.\ $\Compl{\til V\cup \til W}$)
meets $W$ (resp.\ $\til W$) at an annulus $A^\flat$ (resp.\ $\til A^\flat$) 
Thus we can extend the restriction
$f_1\vert_{A^\flat}$ to a homeomorphism 
\[\bar{f}_1:(W, A^\flat)\rightarrow (\til W, \til A^\flat).\]  
Gluing $\bar{f}_1$ and $f_1$ together yields 
a homeomorphism
\[
f_2:(\Compl{V},U,V_1,V_2,W)
\rightarrow 
(\Compl{\til V},\til U,\til V_1,\til V_2,\til W).
\]
Since $V\subset \HK,\til V\subset\tildeHK$ are $3$-balls,
by the Alexander trick,  
$f_2\vert_{\partial V}$ can be extended to a homeomorphism 
\[\bar{f}_2:(V,\partial V)\rightarrow (\til V,\partial \til V).\]
Gluing $\bar{f}_2$ and $f_2$ together yields a homeomorphism  
\[(\sphere, U, W, V_1,V_2,V)
\rightarrow 
(\sphere, \til U, \til W, \til V_1, \til V_2, \til V),\]
and hence an equivalence between
$\pair$ and $\pairtilde$.
\end{proof}
\begin{figure}[h]
\begin{subfigure}{0.47\textwidth}
\centering
\includegraphics[scale=.13]{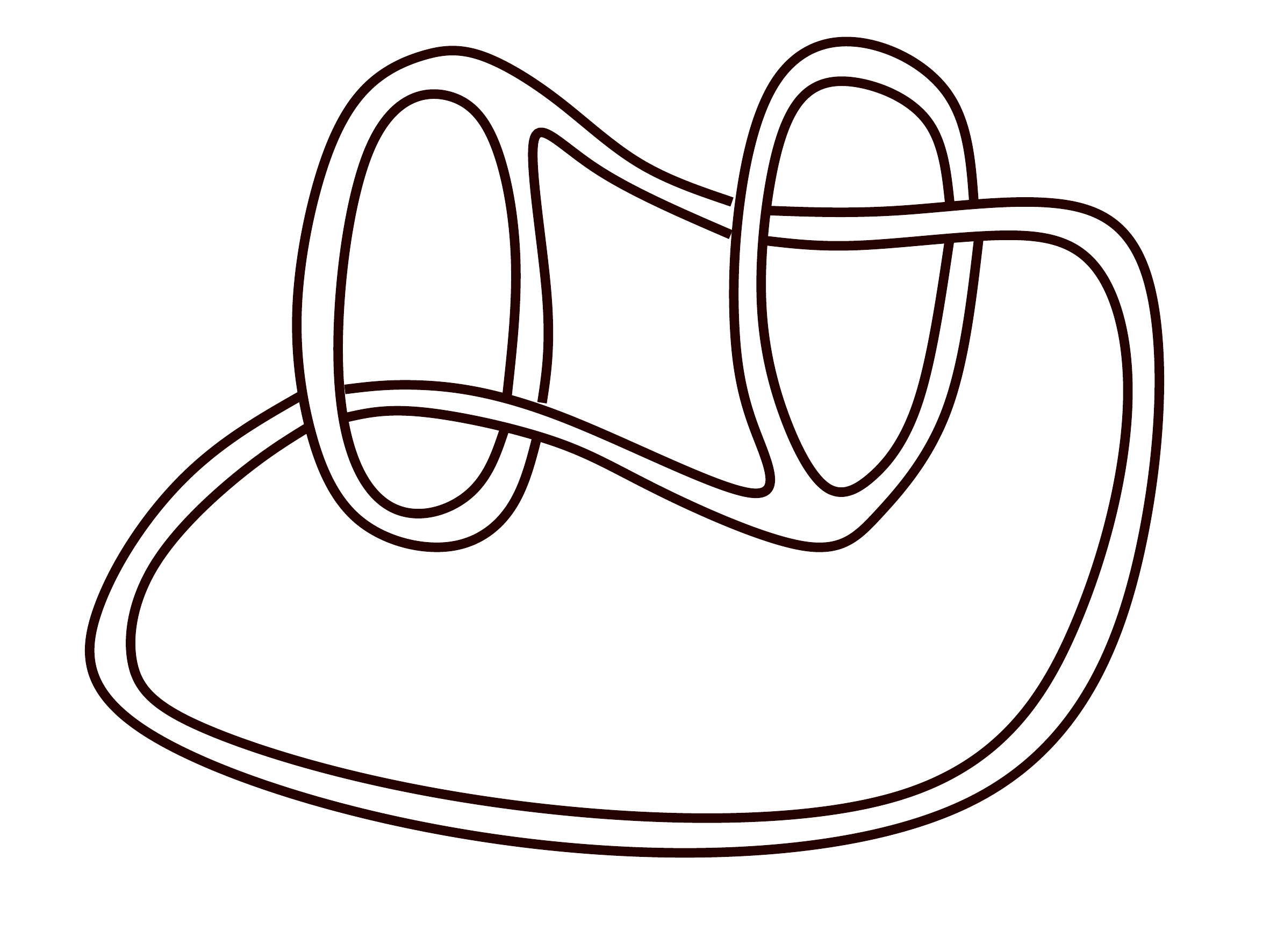}
\caption{$(\sphere,4_1)$.}
\label{fig:hk4_1}
\end{subfigure}
\begin{subfigure}{0.47\textwidth}
\centering
\def\svgwidth{.9 \columnwidth}
\begingroup%
  \makeatletter%
  \providecommand\color[2][]{%
    \errmessage{(Inkscape) Color is used for the text in Inkscape, but the package 'color.sty' is not loaded}%
    \renewcommand\color[2][]{}%
  }%
  \providecommand\transparent[1]{%
    \errmessage{(Inkscape) Transparency is used (non-zero) for the text in Inkscape, but the package 'transparent.sty' is not loaded}%
    \renewcommand\transparent[1]{}%
  }%
  \providecommand\rotatebox[2]{#2}%
  \newcommand*\fsize{\dimexpr\f@size pt\relax}%
  \newcommand*\lineheight[1]{\fontsize{\fsize}{#1\fsize}\selectfont}%
  \ifx\svgwidth\undefined%
    \setlength{\unitlength}{1133.85826772bp}%
    \ifx\svgscale\undefined%
      \relax%
    \else%
      \setlength{\unitlength}{\unitlength * \real{\svgscale}}%
    \fi%
  \else%
    \setlength{\unitlength}{\svgwidth}%
  \fi%
  \global\let\svgwidth\undefined%
  \global\let\svgscale\undefined%
  \makeatother%
  \begin{picture}(1,0.75)%
    \lineheight{1}%
    \setlength\tabcolsep{0pt}%
    \put(0.08836546,0.57861219){\color[rgb]{0,0,0}\makebox(0,0)[lt]{\lineheight{1.25}\smash{\begin{tabular}[t]{l}{\footnotesize $W$}\end{tabular}}}}%
    \put(0,0){\includegraphics[width=\unitlength,page=1]{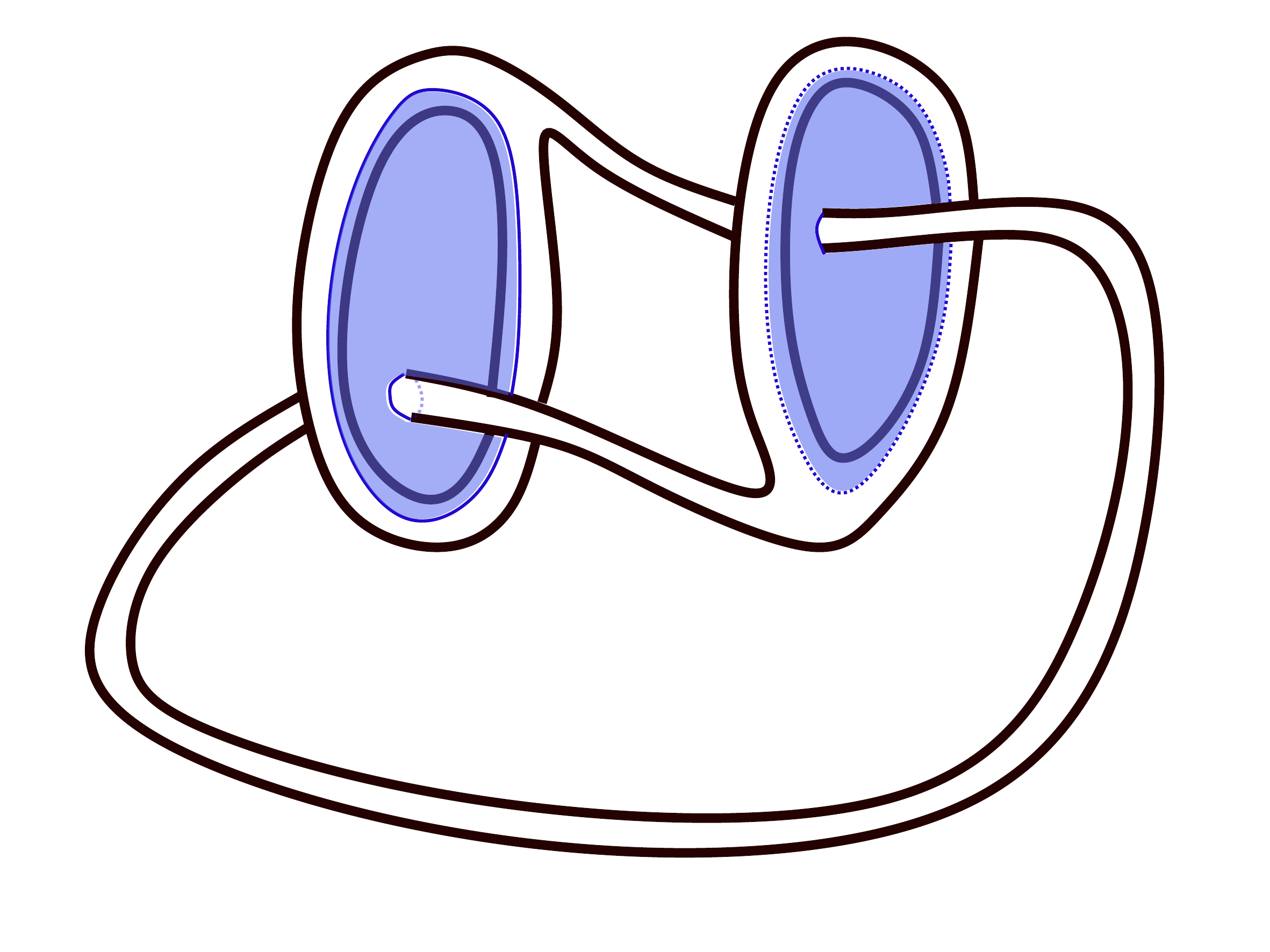}}%
    \put(0.30262144,0.53503469){\color[rgb]{0,0,0}\makebox(0,0)[lt]{\lineheight{1.25}\smash{\begin{tabular}[t]{l}{\footnotesize $A'$}\end{tabular}}}}%
    \put(0.63550499,0.60524287){\color[rgb]{0,0,0}\makebox(0,0)[lt]{\lineheight{1.25}\smash{\begin{tabular}[t]{l}{\footnotesize $A''$}\end{tabular}}}}%
    \put(0,0){\includegraphics[width=\unitlength,page=2]{hk4_1_annuli.pdf}}%
    \put(0.0326831,0.46119507){\color[rgb]{0,0,0}\makebox(0,0)[lt]{\lineheight{1.25}\smash{\begin{tabular}[t]{l}{\footnotesize $A$}\end{tabular}}}}%
  \end{picture}%
\endgroup%

\caption{Annuli in $\Compl {4_1}$.}
\label{fig:hk4_1_annuli}
\end{subfigure}
%
\begin{subfigure}{0.47\textwidth}
\centering
\def\svgwidth{.9 \columnwidth}
\begingroup%
  \makeatletter%
  \providecommand\color[2][]{%
    \errmessage{(Inkscape) Color is used for the text in Inkscape, but the package 'color.sty' is not loaded}%
    \renewcommand\color[2][]{}%
  }%
  \providecommand\transparent[1]{%
    \errmessage{(Inkscape) Transparency is used (non-zero) for the text in Inkscape, but the package 'transparent.sty' is not loaded}%
    \renewcommand\transparent[1]{}%
  }%
  \providecommand\rotatebox[2]{#2}%
  \newcommand*\fsize{\dimexpr\f@size pt\relax}%
  \newcommand*\lineheight[1]{\fontsize{\fsize}{#1\fsize}\selectfont}%
  \ifx\svgwidth\undefined%
    \setlength{\unitlength}{1133.85826772bp}%
    \ifx\svgscale\undefined%
      \relax%
    \else%
      \setlength{\unitlength}{\unitlength * \real{\svgscale}}%
    \fi%
  \else%
    \setlength{\unitlength}{\svgwidth}%
  \fi%
  \global\let\svgwidth\undefined%
  \global\let\svgscale\undefined%
  \makeatother%
  \begin{picture}(1,0.75)%
    \lineheight{1}%
    \setlength\tabcolsep{0pt}%
    \put(0,0){\includegraphics[width=\unitlength,page=1]{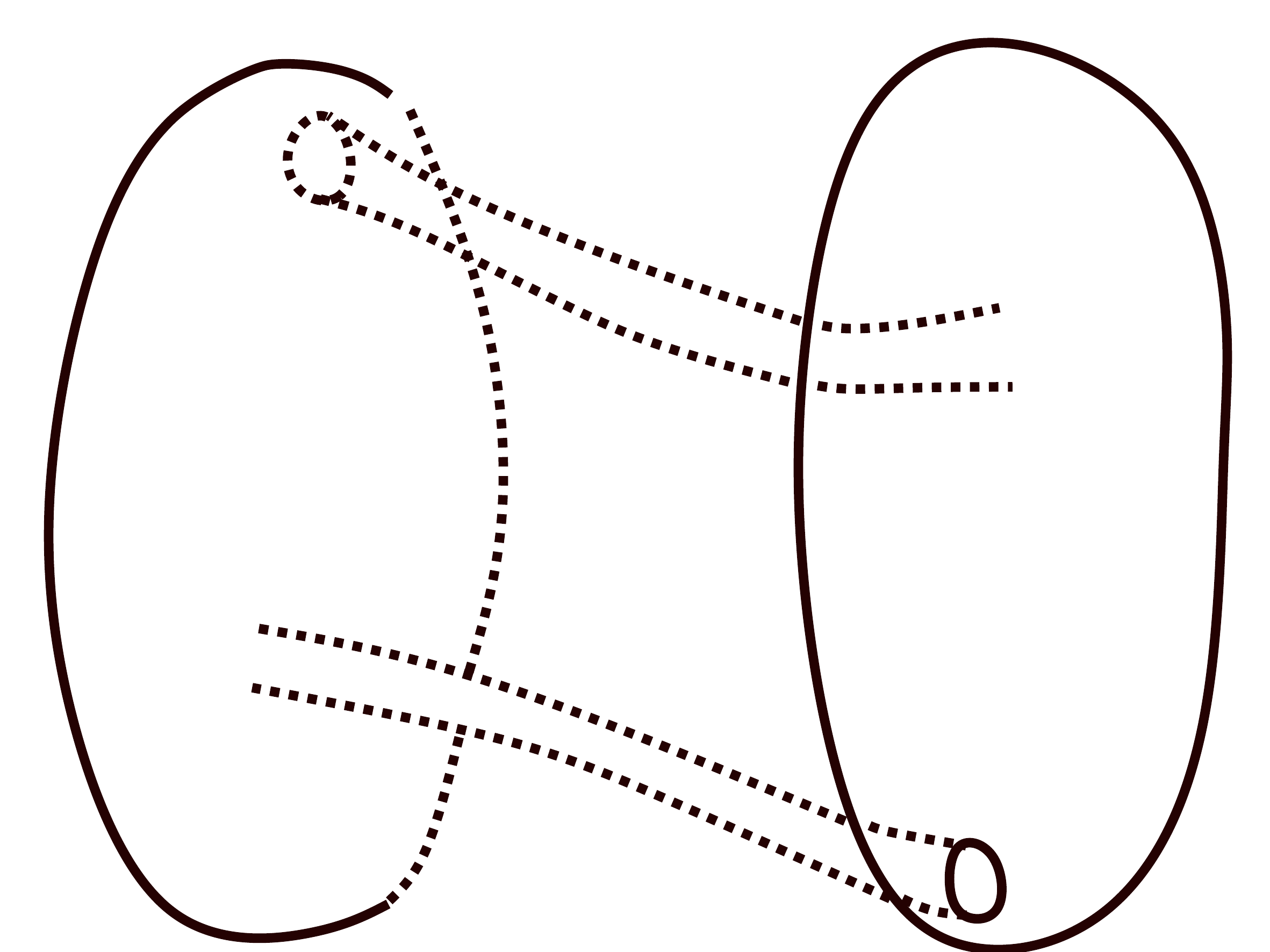}}%
    \put(0.11486233,0.32184701){\color[rgb]{0,0,0}\makebox(0,0)[lt]{\lineheight{1.25}\smash{\begin{tabular}[t]{l}{\footnotesize $A'$}\end{tabular}}}}%
    \put(0.82478418,0.37587362){\color[rgb]{0,0,0}\makebox(0,0)[lt]{\lineheight{1.25}\smash{\begin{tabular}[t]{l}{\footnotesize $A''$}\end{tabular}}}}%
    \put(0,0){\includegraphics[width=\unitlength,page=2]{I_bundle_exterior_hk4_1i.pdf}}%
    \put(0.48535145,0.62510868){\color[rgb]{0,0,0}\makebox(0,0)[lt]{\lineheight{1.25}\smash{\begin{tabular}[t]{l}{\footnotesize $A$}\end{tabular}}}}%
    \put(0,0){\includegraphics[width=\unitlength,page=3]{I_bundle_exterior_hk4_1i.pdf}}%
  \end{picture}%
\endgroup%

\caption{$\Compl {4_1}-\mathring{W}$.}
\label{fig:exterior_minus_seifert_solid_torus}
\end{subfigure}
\begin{subfigure}{0.47\textwidth}
\centering
\def\svgwidth{.9 \columnwidth}
\begingroup%
  \makeatletter%
  \providecommand\color[2][]{%
    \errmessage{(Inkscape) Color is used for the text in Inkscape, but the package 'color.sty' is not loaded}%
    \renewcommand\color[2][]{}%
  }%
  \providecommand\transparent[1]{%
    \errmessage{(Inkscape) Transparency is used (non-zero) for the text in Inkscape, but the package 'transparent.sty' is not loaded}%
    \renewcommand\transparent[1]{}%
  }%
  \providecommand\rotatebox[2]{#2}%
  \newcommand*\fsize{\dimexpr\f@size pt\relax}%
  \newcommand*\lineheight[1]{\fontsize{\fsize}{#1\fsize}\selectfont}%
  \ifx\svgwidth\undefined%
    \setlength{\unitlength}{1133.85826772bp}%
    \ifx\svgscale\undefined%
      \relax%
    \else%
      \setlength{\unitlength}{\unitlength * \real{\svgscale}}%
    \fi%
  \else%
    \setlength{\unitlength}{\svgwidth}%
  \fi%
  \global\let\svgwidth\undefined%
  \global\let\svgscale\undefined%
  \makeatother%
  \begin{picture}(1,0.75)%
    \lineheight{1}%
    \setlength\tabcolsep{0pt}%
    \put(0.2606692,0.16746327){\color[rgb]{0,0,0}\makebox(0,0)[lt]{\lineheight{1.25}\smash{\begin{tabular}[t]{l}{\footnotesize $A'$}\end{tabular}}}}%
    \put(0.65038776,0.45497145){\color[rgb]{0,0,0}\makebox(0,0)[lt]{\lineheight{1.25}\smash{\begin{tabular}[t]{l}{\footnotesize $A''$}\end{tabular}}}}%
    \put(0,0){\includegraphics[width=\unitlength,page=1]{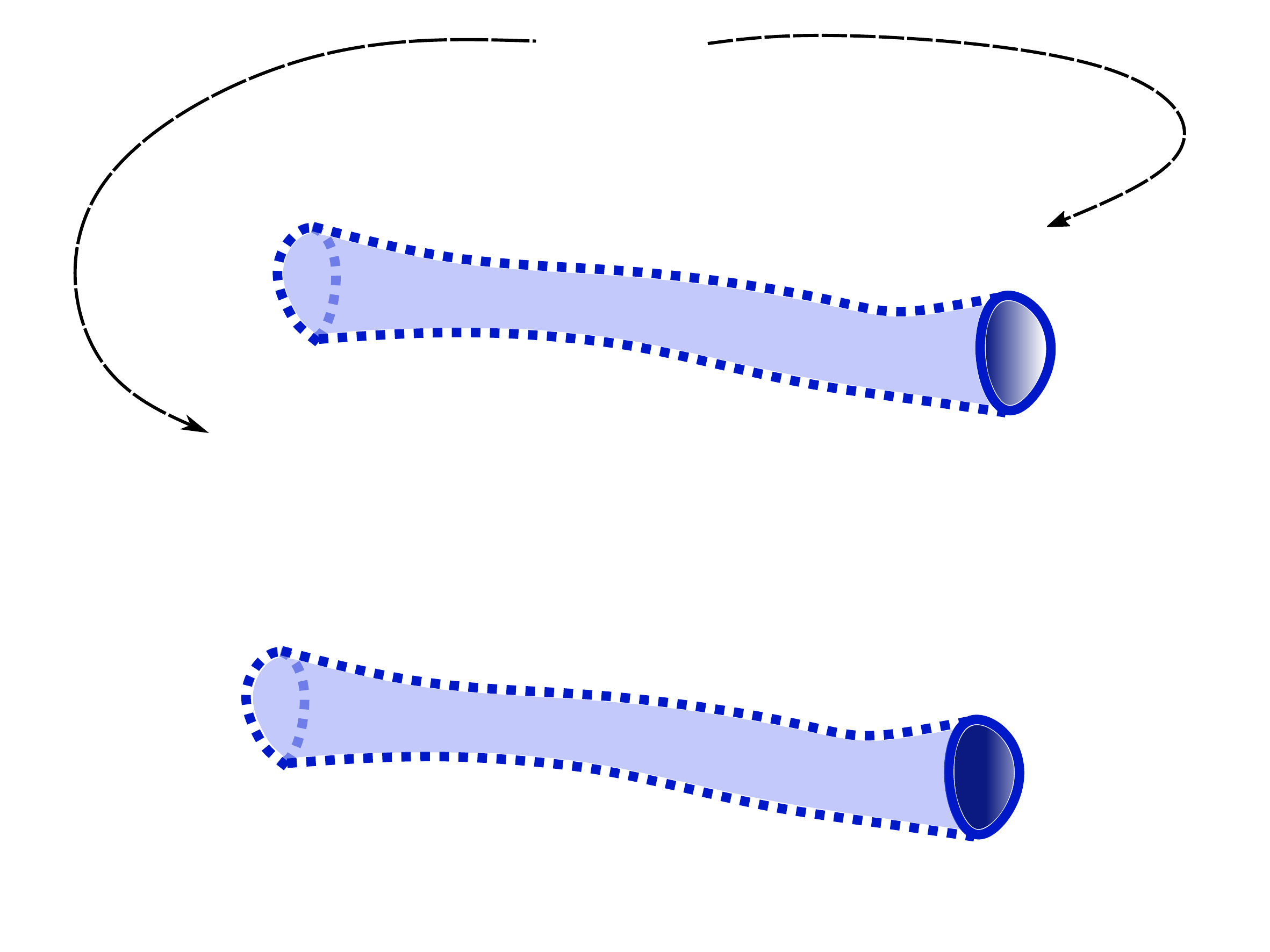}}%
    \put(0.43170265,0.7023507){\color[rgb]{0,0,0}\makebox(0,0)[lt]{\lineheight{1.25}\smash{\begin{tabular}[t]{l}{\footnotesize Lids}\end{tabular}}}}%
    \put(0.46152679,0.58222431){\color[rgb]{0,0,0}\makebox(0,0)[lt]{\lineheight{1.25}\smash{\begin{tabular}[t]{l}{\footnotesize $A$}\end{tabular}}}}%
    \put(0,0){\includegraphics[width=\unitlength,page=2]{I_bundle_exterior_hk4_1ii.pdf}}%
  \end{picture}%
\endgroup%

\caption{I-bundle over a twice-punctured disk.}
\label{fig:I_bundle}
\end{subfigure}
\caption{Annulus diagram of $\pairfourone$.}
\end{figure}
\begin{lemma}\label{lm:ann_diagram_hk4_1}
The annulus diagram of $(\sphere,4_1)$ 
is
\raisebox{-.2 cm}{\includegraphics[scale=.12]{hkfourone_ann}}.
\end{lemma}
\begin{proof}
Recall that $\pairfourone$ is equivalent to 
the handlebody-knot in Fig.\ \ref{fig:hk4_1},
and its exterior admits three annuli 
$A,A',A''$ as depicted in Fig.\ \ref{fig:hk4_1_annuli},
where $A$ is of type $3$-$3$,
and $A',A''$ are of type $2$-$2$.
By Corollary \ref{cor:char_annuli}, 
they are characteristic and hence 
the characteristic diagram of $\Compl{4_1}$ is 
\raisebox{-.3 cm}{\includegraphics[scale=.1]{theta_char}}, $\square=\hnode$ or $\snode$.
Let $W\subset\ComplHK$ 
be the Seifert fibered solid torus cut off by $A\cup A'\cup A''$ (Fig.\ \ref{fig:hk4_1_annuli}).
Then as shown in Fig.\ \ref{fig:exterior_minus_seifert_solid_torus} and
\ref{fig:I_bundle}, the exterior
of $W\subset \ComplHK$ together with $A\cup A'\cup A''$
is an I-bundle over a pair of pants, and hence the assertion.
\end{proof}


\begin{theorem}\label{teo:characterization_fourone}
$\charE$ is 
\raisebox{-.35 cm}{\includegraphics[scale=.1]{hkfourone_char}}
if and only if $\pair$ is equivalent to $(\sphere,4_1)$. 
\end{theorem}
\begin{proof}
This follows from Theorem \ref{teo:char_three_bigons} and
Lemmas \ref{lm:rigidity_solid_three_bigons} and \ref{lm:ann_diagram_hk4_1}.
\end{proof}

\section{Handlebody-knot symmetries}\label{sec:symmetry}

In this section, we compute the symmetry groups of handlebody-knots whose exteriors contain a type $2$ annulus, based on
the classification in Theorem \ref{teo:classification_char_diagram_type_two}.

\subsection{Mapping class group}
Here we 
recall some properties of mapping class groups.
Given subpolyhedra $X_1,\dots,X_n$ of a manifold $M$, 
the space of self-homeomorphisms of $M$ preserving 
$X_i$, $i=1,\dots, n$, 
setwise (resp.\ pointwise) is denoted by 
\[\Aut{M,X_1,\dots,X_n} \qquad\text{(resp.\ $\Aut{M,\rel X_1,\dots ,X_n}$ )},\]
and the mapping class group of $(M,X_1,\dots,X_n)$ is defined as 
\begin{multline*}
\MCG{M,X_1,\dots,X_n}:=\pi_0\big(\Aut{M,X_1,\dots ,X_n}\big)\\
\text{\Big(resp.\ $\MCG{M,\rel X_1,\dots,X_n}:=\pi_0\big(\Aut{M,\rel X_1,\dots ,X_n}\big)$\Big )}.
\end{multline*} 
The ``+'' subscript is added when considering 
only orientation-preserving
homeomorphisms:
\begin{align*}
\pAut{M,X_1,\dots, X_n}&\quad \text{\big(resp.\ $\pAut{M,\rel X_1,\dots, X_n}$\big)},\\
\pMCG{M, X_1,\dots,X_n}&\quad\text{\big(resp.\ $\pMCG{M,\rel X_1,\dots , X_n}$\big)}.
\end{align*}
Given $f\in\Aut{M,X_1,\cdots,X_n}$, 
$[f]$ denotes the mapping class it represents.
If $M=\sphere$, we call the mapping class group
the \emph{symmetry group} of $(M,X_1,\dots,X_n)$. 
\begin{lemma}[Cutting Homomorphism, {\cite[Proposition $3.20$]{FarMar:11}}]\label{lm:cut_homo}
Let $\Sigma$ be a closed surface and 
$\alpha_1,\dots,\alpha_n$ 
mutually disjoint and non-homotopic simple loops in $\Sigma$. Then there is a well-defined
homomorphism
\[\cut:\pMCG{\Sigma,[\alpha_1],\dots,[\alpha_n]}
\rightarrow
\pMCG{\Sigma-\rnbhd{\alpha_1\cup\cdots\cup\alpha_n}}
\]
whose kernel is generated by the Dehn twists
about $\alpha_1,\dots,\alpha_n$, where
the group
\[\pMCG{\Sigma,[\alpha_1],\dots,[\alpha_n]}\]
is the subgroup of $\pMCG{\Sigma}$ given by 
homeomorphisms that preserve the isotopy classes of $\alpha_1,\dots,\alpha_n$,
respectively.
\end{lemma}

Then next two lemmas are proved in \cite{ChoKod:13}
and \cite{FunKod:20} (see also \cite[Remark $2.1$]{Wan:21}).
\begin{lemma}[{\cite[Lemma $2.3$]{ChoKod:13}}]\label{lm:cho_koda_lemma}
If $\pair$ is atoroidal, then
\[\pMCG{\ComplHK,\rel \bCompl\HK}\simeq \{1\}.\]
\end{lemma}

\begin{lemma}[\cite{FunKod:20}]\label{lm:fun_koda_lemma}
$\pair$ is non-trivial and atoroidal if and only if
$\Sym\HK$ is finite.
\end{lemma}

\begin{lemma}\label{lm:solid_torus_w_bp_MCG}
Let $(W,\bp w)$ be a solid torus with  
boundary pattern, where 
$\bp w=\{\elem G_1,\elem G_2,\dots, \elem G_n\}$,
and $\elem G_i$, $i=1,\dots,n$, are all annuli, and $\vert \bp w\vert=\partial W$. 
Suppose $f\in \pAut{W,\elem G_1,\dots, \elem G_n}$ 
does not swap the components of $\partial \elem G_1$---which holds automatically when $n>2$.
Then $f$ is isotopic to $\id$ in $\pAut{W,\elem G_1,\dots,\elem G_n}$.
\end{lemma}
\begin{proof}
Without loss of generality,
it may be assumed that 
$\elem G_i\cap \elem G_j\neq\emptyset$ 
if and only if $\vert i-j\vert\leq 1$.
Denote by $\elem U_k$ the union 
$\elem G_1\cup \cdots\cup \elem G_k$ and set $\elem U_0=\emptyset$.
Observe that, if 
$f\vert_{\elem U_{k-1}}=\id, 1\leq k\leq n$,
then $f$ can be isotoped in   
\begin{equation}\label{eq:isotopy_U_k}
\pAut{W,\elem G_k,\dots,\elem G_n,\rel \elem U_{k-1}},
\end{equation} 
so that $f\vert_{U_k}=\id$.
To see this, we first isotope 
$f\vert_{\elem U_k}$ to $\id$ in
$\Aut{\elem U_k,\rel \elem U_{k-1}}$ as follows: 
In the case $k=1$, it results from
the assumption that
$f$ does not swap components of $\partial \elem G_1$,
whereas if $1<k<n$, it follows from the fact that 
$\MCG{\elem U_k,\rel \elem U_{k-1}}=\{1\}$. 
If $k=n$, it is a consequence of $f$ sending 
meridian disks of $W$ to themselves.
Via a regular neighborhood of $\elem U_k$
in $W$, the isotopy of $f\vert_{\elem U_k}$
can be extended to an isotopy in \eqref{eq:isotopy_U_k} that isotopes $f$ so that $f\vert_{\elem U_k}=\id$. 
Hence by induction, we may assume 
$f\in \Aut{W,\rel \partial W}$, and 
the assertion follows since 
$\MCG{W,\rel W}\simeq \{1\}$.
\end{proof}

\begin{lemma}\label{lm:one_annulus_MCG}
Let $W$ be a solid torus and $A\subset \partial W$
an annulus with $H_1(A)\rightarrow H_1(W)$
non-trivial and not an isomorphism.
Then $\MCG{W,A}\simeq\pMCG{W,A}\simeq \mathbb{Z}_2$. 
\end{lemma}
\begin{proof}
Identify $W$ with the subspace of $\mathbb{C}^2$
\[\{(z_1,z_2)\in\mathbb{C}^2\mid z_1=re^{i\theta},z_2=\sqrt{2-r^2}e^{i\phi},-\pi\leq \theta,\phi\leq  \pi, 0\leq r\leq 1\},\]
and $A\subset W$ with
\[\{z_1=e^{ip (t+s)},z_2=e^{iq t},-\pi\leq t \leq \pi, -\epsilon \leq s\leq \epsilon\},\]
where $p,q$ are coprime integers
with $q>0$ and $0<\epsilon <\frac{\pi}{q}$. 
Let $c_A,c_W$ be the cores of $A,W$ given by $s=0,r=0$, respectively,
and orient them so that $[c_A]=q [c_W]\in H_1(W)$
By the assumption, $q\neq 0$ or $\pm 1$, and therefore
the linking number 
$\lk{c_A}{c_W}$ is non-trivial and 
every homeomorphism $f$
of $(W,A)$ either preserves the orientations
of both $c_A,c_W$ or reverses them.
This implies $f$ is orientation-preserving.

On the other hand,  the conjugation 
\[\conj:(z_1,z_2)\in W\mapsto (\bar{z}_1 ,\bar{z}_2)\in W\]
preserves $A$ but swaps its boundary components, so
it is is non-trivial in 
$\pMCG{W,A}$. By Lemma \ref{lm:solid_torus_w_bp_MCG},
it generates the entire group  
since any homeomorphism $f\in\Aut{W,A}$ 
either swap components of $\partial A$ or preserve them. 
  
\end{proof}

\begin{lemma}\label{lm:two_annuli_MCG}
Let $W$ be a solid torus, and $A_1,A_2\subset \partial W$ 
two disjoint annuli 
with $H_1(A_i)\rightarrow H_1(W)$, $i=1,2$, isomorphisms. 
Then 
$\pMCG{W,A_1\cup A_2}\simeq \mathbb{Z}_2\times \mathbb{Z}_2$ and  
$\MCG{W,A_1\cup A_2}\simeq \mathbb{Z}_2\times \mathbb{Z}_2\times \mathbb{Z}_2$.
\end{lemma}
\begin{proof}
Identify $W$ with $\Q\times S^1\subset \mathbb{R}^2\times \mathbb{C}$, where $S^1$ is the unit circle $\{z=e^{i\theta}\}$
and $\Q$ is the square given by  
\[\{(x,y)\mid -1\leq x,y\leq 1\}.\]
Identify $A_1,A_2$ with the annuli given by $y=\pm 1$,
and their cores $c_1,c_2$ the loops given by 
$x=0$, and 
denote by $B_1,B_2$ the annuli in
the closure of $\partial W-A$.

Consider $\rot_i \in \pAut{W,A_1\cup A_2}, i=1,2$, defined by
the assignments:
\begin{align*}
\Q\times S^1&\rightarrow \Q\times S^1\\
(x,y,z)&\mapsto (-x,-y,z),\\
(x,y,z)&\mapsto (-x,y,\bar{z})
\end{align*}
respectively. Note that $\rot_1,\rot_2$
both are of order $2$ and commute with each other.
In addition, $\rot_1$
swaps $A_1,A_2$ and also $B_1,B_2$,
whereas $\rot_2$ swaps $A_1,A_2$ but preserves $B_1,B_2$, so
their composition $\rot_1\circ \rot_2$ swaps $B_1,B_2$ but preserves $A_1,A_2$. This implies they represent
distinct mapping classes.
Since every $f\in \Aut{W,A_1\cup A_2}$ 
either swaps $A_1,A_2$ (resp.\ $B_1,B_2$)
or preserves them, by Lemma \ref{lm:solid_torus_w_bp_MCG}, 
$\{[\rot_1],[\rot_2]\}$ generates $\pMCG{W,A_1\cup A_2}$.

To see $\MCG{W,A_1\cup A_2}\simeq \Z_2\times \Z_2\times \Z_2$,
consider $\mir\in \Aut{W,A_1\cup A_2}$ defined by the assignment
\begin{align*}
\Q\times S^1&\rightarrow \Q\times S^1\\
(x,y,z)&\mapsto (-x,y,z), 
\end{align*}
which is orientation-reversing, commutes with $\rot_i,i=1,2$, and together with $\rot_i,i=1,2$, generates
$\MCG{W,A_1\cup A_2}$.
\end{proof}

\begin{lemma}\label{lm:three_annuli_MCG}
Let $W$ be a solid torus and $A_1,A_2, A_3\subset \partial W$ 
three disjoint annuli with 
$H_1(A_i)\rightarrow H_1(W)$, $i=1,2,3$, isomorphisms. 
Then  
$\pMCG{W,A_1, A_2\cup A_3}\simeq \mathbb{Z}_2$ and 
$\MCG{W,A_1, A_2\cup A_3}\simeq \mathbb{Z}_2\times \mathbb{Z}_2$.
\end{lemma}
\begin{proof}
Identify $W$ with $\Hex\times S^1\subset\mathbb{C}\times \mathbb{C}$, where $S^1\subset \mathbb{C}$ is the unit circle,
and $\Hex\subset\mathbb{C}$ the regular hexagon 
with center at origin and vertices $v_k=e^{\frac{2\pi k}{6}}$, 
$k=1,\dots,6$. Identify $A_k$ with the product of $S^1$ and the edge $e_k$ connecting $v_{2k-1},v_{2k}$, $k=1,2,3$.
Denote by $\rot\in \pAut{W,A_1,A_2\cup A_3}$
the homeomorphism given by
\begin{align*}
\Hex\times S^1&\rightarrow \Hex\times S^1\\
(z_1,z_2)&\mapsto (-\bar{z}_1,\bar{z}_2);
\end{align*}
$\rot$ swaps $A_2,A_3$ and hence represents a non-trivial 
mapping class in $\pMCG{W,A_1,A_2\cup A_3}$.
Since every $f\in\pAut{W,A_1,A_2\cup A_3}$ 
either swaps $A_2,A_3$ or preserves them,
by Lemma \ref{lm:solid_torus_w_bp_MCG}, 
either $[f]=[\rot]$ or $[f]$ is trivial, so $\pMCG{W,A_1,A_2\cup A_3}\simeq \Z_2$.
On the other hand, there is an orientation-reversing homeomorphism 
$\mir\in \Aut{W,A_1,A_2\cup A_3}$ defined by
\begin{align*}
\Hex\times S^1&\rightarrow \Hex\times S^1\\
(z_1,z_2)&\mapsto (z_1,\bar{z}_2),
\end{align*}
which is of order $2$ and commutes with $\rot$, and 
$\{[\rot],[\mir]\}$ generates $\MCG{W,A_1,A_2\cup A_3}$.
\end{proof}

The next lemma follows from 
\cite[Section $2$]{Hat:76} (see also 
\cite[Theorem $1$]{Hat:99}).
\begin{lemma}\label{lm:NA_preserving_homo}
Given a handlebody-knot $\pair$ and an essential surface $S$ in $\Compl\HK$, 
the natural homomorphisms
\begin{align*}
\Sym{\HK,S}&\rightarrow \Sym\HK,\\
\Sym{\HK,\rnbhd S }&\rightarrow \Sym\HK
\end{align*} 
are injective.
\end{lemma}

\subsection{Symmetry groups of handlebody-knots} 
Here $\pair$ is an atoroidal handlebody-knot, 
and $A\subset\Compl\HK$ a type $2$ essential annulus. 
The symbols $l,l_A,\HK_A, A_+,A_-$ are as in Section \ref{sec:classification}. In addition, we identify 
the intersection $\rnbhd{A}\cap \partial\HK$
with $\rnbhd{l\cup l_A}=\rnbhd{l}\cup \rnbhd{l_A}$.


\begin{theorem}\label{teo:symmetry_type_2_1}
If $A$ is of type $2$-$1$,
then $\pSym\HK<\Z_2$ and $\Sym\HK<\Z_2\times \Z_2$. 
\end{theorem}
\begin{proof}
Note first that the injection
$\pSym{\HK,\rnbhd A}\rightarrow \pSym\HK$
in Lemma \ref{lm:NA_preserving_homo}
is an isomorphism since $A$ is unique by 
Theorem \ref{teo:classification_char_diagram_type_two},
composing its inverse with the homomorphism 
$\pSym{\HK,\rnbhd A}\xrightarrow{\Phi} \pMCG{\rnbhd A, A_+\cup A_-}$
given by restriction to $\rnbhd{A}$ yields the homomorphism
\[\pSym\HK\simeq \pSym{\HK,\rnbhd A}\rightarrow \pMCG{\rnbhd A, A_+\cup A_-}.\] 
By Lemma \ref{lm:two_annuli_MCG},
it then suffices to show the injectivity of $\Phi$
as it entails the injectivity of 
\[\Sym{\HK,\rnbhd A}\rightarrow \MCG{\rnbhd{A},A_+\cup A_-}.\]

To see $\Phi$ is injective, let $[f]\in \pSym{\HK,\rnbhd{A}}$ with 
$\Phi([f])=1$.
This implies $f\vert_{\partial\HK-\rnbhd{l\cup l_A}}$
does not permute punctures of 
the four-times punctured sphere
$\partial\HK-\rnbhd{l\cup l_A}$, and
thus $[f\vert_{\partial\HK-\rnbhd{l  \cup l_A}}]=1
\in \pMCG{\partial\HK-\rnbhd{l\cup l_A}}$
since $[f\vert_{\partial\HK-\rnbhd{l \cup l_A}}]$ is of finite order by Lemma \ref{lm:fun_koda_lemma}.
Again by Lemma \ref{lm:fun_koda_lemma}, $[f\vert_{\partial\HK}]$
is of finite order in $\pMCG{\partial\HK,[l],[l_A]}$; hence by Lemma \ref{lm:cut_homo}, it is the identity.
Because $f\vert_{\partial\HK}$ is isotopic to $\id$,
$f$ can be isotoped in $\Aut{\sphere,\HK}$   
so that $f\vert_{\partial\HK}=\id$. Applying Lemma \ref{lm:cho_koda_lemma}, one can further isotope $f$ to $\id$ in 
$\Aut{\sphere,\rel \partial \HK}$.
\end{proof}

\begin{theorem}\label{teo:symmetry_one_type_2_2}
If $A\subset\ComplHK$ is the unique type $2$-$2$ annulus, 
then $\pSym\HK\simeq\{1\}$ and $\Sym\HK<\Z_2$. 
If in addition $\ComplHK$ admits an annulus $A'$ 
of another type, then $\Sym\HK\simeq \pSym\HK\simeq\{1\}$.
\end{theorem}
\begin{proof}
As in the previous case, the uniqueness of $A$ gives us
the homomorphism
\[\pSym\HK\simeq \pSym{\HK,\rnbhd A}\xrightarrow{\Phi} \pMCG{\rnbhd A, A_+\cup A_-}.\]
The first assertion follows once we show 
the injectivity of $\Phi$ because, given
any $f\in\pAut{\sphere,\HK,\rnbhd A}$, it
can neither swap $A_+, A_-$ nor swap $\rnbhd{l},\rnbhd{l_A}$
by the definition of a type $2$-$2$
annulus. On the other hand, The second assertion 
can be derived from the first as follows: by Theorem \ref{teo:classification_char_diagram_type_two}, $A'$ is the unique type $3$-$2$ annulus in $\ComplHK$. 
Let $W\subset \ComplHK$ be the solid torus cut off by $A'$. 
Then by the essentiality of $A$, $H_1(A)\rightarrow H_1(W)$
is non-trivial and not an isomorphism.
On the other hand, by Lemma \ref{lm:NA_preserving_homo}, 
there is a homomorphism 
\[\Sym\HK\simeq \Sym{\HK,A}=\Sym{\HK,W}\rightarrow \MCG{W,A}.\]
Now, if $\Sym\HK$ is non-trivial, then by the first assertion, $\MCG{W,A}$
contains a mapping class represented by an orientation-reversing
homeomorphism, contradicting Lemma \ref{lm:one_annulus_MCG}.

We now prove the injectivity of $\Phi$. 
Let $[f]\in\pSym\HK$ with $\Phi([f])=1\in \pMCG{\rnbhd{A},A_+\cup A_-}$. 
We can isotope $g:=f\vert_{\partial\HK}$
in $\pAut{\partial\HK,\rnbhd{l\cup l_A}}$
so that $g\vert_{\rnbhd{l\cup l_A}}=\id$.
Let $D$ be the meridian disk disjoint from $l_A$
and dual to $l$. Then one can further isotope
$g$ in $\pAut{\partial\HK,\rel \rnbhd{l\cup l_A}}$
so that $g\vert_{\rnbhd{\partial D}}=\id$. 
In other words, $f\vert_{\partial\HK}$
represents a mapping class in $\pMCG{\partial\HK, \rel \rnbhd{\partial D\cup l}}$.
Now, the homomorphism induced by the inclusion
\[
\pMCG{\partial\HK,\rel \rnbhd{\partial D\cup l}}
\rightarrow \pMCG{\partial\HK}
\]
is injective by \cite[Theorem $3.18$]{FarMar:11}, and 
by Lemma \ref{lm:fun_koda_lemma},
$[f\vert_{\partial\HK}]\in 
\pMCG{\partial \HK}$ is of finite order, 
so $[f\vert_{\partial\HK}]\in\pMCG{\partial \HK,\rel \rnbhd{\partial D\cup l}}$ 
is also of finite order.
The group $\pMCG{\partial\HK, \rel \rnbhd{\partial D\cup l}}$ is,
however, torsion free, and hence $f\vert_{\partial\HK}$
is isotopic to $\id$ in $\pAut{\partial\HK}$.
We may thence isotope $f$ in $\pAut{\sphere,\HK}$
so that $f\vert_{\partial\HK}=\id$. 
By Lemma \ref{lm:cho_koda_lemma}, $f$
can be further isotoped to $\id$ in $\pAut{\sphere,\rel \partial \HK}$.
\end{proof}

\begin{theorem}\label{teo:symmetry_two_type_2_2}
If $A\subset\ComplHK$ 
is of type $2$-$2$ but not the unique type $2$-$2$ annulus,
then $\pSym\HK<\Z_2$ and $\Sym\HK<\Z_2\times \Z_2$. 
\end{theorem}
\begin{proof}
By Theorem \ref{teo:classification_char_diagram_type_two},
$\ComplHK$ admits a unique type $3$-$3$ annulus $A_0$, and
exactly two non-isotopic annuli $A,A'$, which
cut off a solid torus $W\subset\ComplHK$ and form a characteristic surface of $\ComplHK$;
this together with Lemma \ref{lm:NA_preserving_homo}
gives us the homomorphism
\begin{multline*}
\pSym\HK\simeq \pSym{\HK,A_0,A_1\cup A_2}\\
=\pSym{\HK,W}\xrightarrow{\Phi} \pMCG{W,A_0,A_1\cup A_2}.
\end{multline*}
It suffices to prove that $\Phi$ is injective,
in view of Lemma \ref{lm:three_annuli_MCG}.

Let $[f]\in\pSym{\HK,W}$ with $\Phi([f])=1\in\pMCG{W,A_0,A_1\cup A_2}$.
Note that $\partial\HK\cap W$ consists of three annuli
$B_0,B_1,B_2$; denote by $c_0,c_1,c_2$
their cores, respectively.
Since $\Phi([f])=1$, 
$f\vert_{\partial\HK-(B_0\cup B_1\cup B_2)}$
does not permute punctures of $\partial\HK-(B_0\cup B_1\cup B_2)$, which is two copies of the three-times punctured sphere, 
and therefore $[f\vert_{\partial\HK-(B_0\cup B_1\cup B_2)}]=1\in\pMCG{\partial\HK-(B_0\cup B_1\cup B_2)}$. On the other hand by Lemma \ref{lm:fun_koda_lemma},
$[f\vert_{\partial\HK}]$ is of finite order in $\pMCG{\partial\HK,[c_0],[c_1],[c_2]}$, and hence trivial therein by Lemma \ref{lm:cut_homo}; in particular, 
$f\vert_{\partial\HK}$ is isotopic to $\id$ in 
$\pAut{\partial\HK}$.
We then isotope $f$ in $\pAut{\sphere,\HK}$
so that $f\vert_{\partial\HK}=\id$; by 
Lemma \ref{lm:cho_koda_lemma}, we can further isotope $f$
to $\id$ in $\pAut{\sphere,\rel\partial\HK}$.
\end{proof}

\section{Irreducibility and atoroidality}\label{sec:irre_atoro}
 
Let $\pair$ be a handlebody-knot, not necessarily atoroidal, and $A\subset\Compl\HK$ 
a type $2$ annulus, not necessarily essential.
The symbols $l_A,l\subset \partial A$, 
$\HK_A$, and $A_+, A_-, l_+,l_-\subset \partial \HK_A$ are as in Section \ref{sec:classification}.
We say $A$ is \emph{unknotting} if $\pairA$ 
is trivial.

\subsection{Essentiality, irreducibility and triviality}
\begin{lemma}\label{lm:esse_redu_type2_1}
If $A$ is of type $2$-$1$, then the following are equivalent:
\begin{enumerate}[label=(\roman*)]
\item $\pair$ is reducible.\label{lm:esse_redu_type2_1:pair}
\item $A$ is inessential.\label{lm:esse_redu_type2_1:A}
\item $\pairA$ is reducible and there exists 
a disk $D$ meeting $l_+\cup l_-$ at one point.
\label{lm:esse_redu_type2_1:pairA}
\end{enumerate}
\end{lemma}
\begin{proof}
Note first that by the definition 
$A$ is incompressible.
 
\ref{lm:esse_redu_type2_1:pair} 
$\Rightarrow$\ref{lm:esse_redu_type2_1:A}: 
Let $D\subset\ComplHK$ be a compressing disk 
of $\bComplHK$. Minimize $\# D\cap A$
in the isotopy class of $A$.
If $D\cap A=\emptyset$,
then $\partial D\subset\bComplHK$ 
is separating, and hence $D\subset\ComplHK$ is separating. 
Since $\partial A\subset\partial \HK$ is non-parallel and 
non-separating, 
components of $\partial A$
lie in different components of $\bComplHK-\partial D$, 
contradicting that $A$ is connected. 
If $D\cap A\neq\emptyset$, then, 
since $A$ is incompressible, 
any outermost disk in $D$ cut off by $D\cap A$
is a $\partial$-compressing disk of $A$ by the minimality.

\ref{lm:esse_redu_type2_1:A}$\Rightarrow$\ref{lm:esse_redu_type2_1:pairA} \&
\ref{lm:esse_redu_type2_1:A} 
$\Rightarrow$\ref{lm:esse_redu_type2_1:pair}:
Since $A$ is incompressible, it is $\partial$-compressible.
Let $D$ be a $\partial$-compressing disk of $A$.
Then $D$ induces a disk in $\Compl\HKA$
meeting $l_+\cup l_-$ at one point,
and hence $\pairA$ is reducible. 
On the other hand, the the frontier of a regular neighborhood of $A\cup D\subset \ComplHK$ 
is a $\partial$-compressing disk of $\partial \Compl\HK$, so $\pair$ is reducible.

\ref{lm:esse_redu_type2_1:pairA}$\Rightarrow$\ref{lm:esse_redu_type2_1:A}: The disk $D$
induces a $\partial$-compressing disk of $A$. 
\end{proof}
\begin{remark}
The union $l_+\cup l_-$ in 
Lemma \ref{lm:esse_redu_type2_1}
\ref{lm:esse_redu_type2_1:pairA}   
is necessary; there are irreducible handlebody-knots $\pair$ with $\pairA$ trivial, and one of $l_+,l_-$ is primitive, for instance, handlebody-knots in Fig.\ \ref{fig:looping_n_2_torus_link}.
\end{remark}
\begin{lemma}\label{lm:triv_type2_1}
Let $A$ be of type $2$-$1$.
Then $\pair$ is trivial if and only if 
$\pairA$ is trivial and
$\{l_+,l_-\}$ is primitive.
\end{lemma}
\begin{proof}
``$\Rightarrow$'':
By Lemma \ref{lm:esse_redu_type2_1} 
there exists a disk $D$ meeting $l_+\cup l_-$ 
at one point, say  
$D\cap l_+\neq\emptyset$.
Then the frontier of a regular neighborhood
$\rnbhd{A_+\cup D}$ of $A_+\cup D\subset \ComplHKA-l_-$
is an essential separating disk $D'\subset \ComplHKA$, 
which splits $\ComplHKA$
into two parts: a solid torus where $l_+$ lies and $D$ is a meridian disk and the exterior $\Compl K$ of a knot $(\sphere, K)$ 
where $l_-\subset \bCompl K$ is a meridian of $(\sphere, K)$.
If $\pairA$ is non-trivial, then
$(\sphere,K)$ is non-trivial and $\bCompl K$ induces 
an incompressible torus $T$
in $\ComplHKA$. 
$T$ is also incompressible in $\ComplHK$,
for given any compressing disk $D$ of $T$,
one can always isotope $A$ away from $D$,
given the incompressibility of $A$, contradicting 
$\pair$ is trivial.
So $(\sphere, K)$ is trivial, and 
$\Compl K$ is a solid torus with $l_-$
primitive in $\Compl K$, and hence the assertion.

``$\Leftarrow$'':
By \cite{Zie:70} (see also \cite{Gor:87}),
there exists a basis $\{ x_+,x_-\}$ of 
$\pi_1(\ComplHKA)$ with $x_\pm$
in the conjugate classes determined by $l_\pm$, respectively.
Since $\pi_1(\ComplHK)$ is the HNN extension of $\pi_1(\ComplHKA)$ 
with respect to
$\pi_1(A)$, $\pi_1(\ComplHK)$ is free, so $\pair$ is trivial.
\end{proof}

\begin{lemma}\label{lm:esse_redu_type2_2}
If $A$ is of type $2$-$2$, then the following are equivalent:
\begin{enumerate}[label=(\roman*)]
\item $\pair$ is reducible.\label{lm:esse_redu_type2_2:pair}
\item $A$ is inessential.\label{lm:esse_redu_type2_2:A}
\item $\pairA$ is reducible and $l_-$ is homotopically
trivial in $\ComplHKA$.
\label{lm:esse_redu_type2_2:pairA}
\end{enumerate}
\end{lemma}
\begin{proof}
\ref{lm:esse_redu_type2_2:pair}$\Rightarrow$\ref{lm:esse_redu_type2_2:A}:
Let $D$ be
an essential disk in $\ComplHK$.
Minimize $\# D\cap A$ in the isotopy class of $A$.
Suppose $D\cap A=\emptyset$.  
Then $\partial D$ lies in the once-punctured torus $T$ in $\partial\HK_A-l_+\cup l_-$. 
If $\partial D$ is separating,
then $\partial D$ is parallel to $l_-$, and so $A$ is compressible. 
If $\partial D$ is non-separating,  
then there is a loop $l$ in $T$
meeting $\partial D$ once. The frontier of a regular neighborhood
of $D\cup l$ in $\ComplHKA-l_-$ is an essential separating disk 
disjoint from $A$, and therefore, as in the previous case, $A$ compressible.
Suppose $D\cap A$ contains a circle, then 
any innermost disk in $D$ cut off by $D\cap A$
is a compressing disk of $A$.
If $D\cap A$ contains only arcs, then 
an outermost disk $D'$ in $D$ cut off by $D\cap A$ 
either is a $\partial$-compressing disk of $A$
or induces an essential disk $D''$ disjoint from $A$ in $\ComplHK$; 
either way implies $A$ is inessential.

\ref{lm:esse_redu_type2_2:A} $\Rightarrow$ \ref{lm:esse_redu_type2_2:pairA} \& 
\ref{lm:esse_redu_type2_2:A} $\Rightarrow$ \ref{lm:esse_redu_type2_2:pair}:
Consider first the case $A$ is compressible. Then any compressing disk 
$D$ induces a disk $D'\subset\ComplHKA$ with $\partial D'=l_-$
and a disk $D''\subset\ComplHK$ with $\partial D''=l_A$,
and therefore \ref{lm:esse_redu_type2_2:pairA} and 
\ref{lm:esse_redu_type2_2:pair}.
Now if $A$ is $\partial$-compressible, and $D$
is a $\partial$-compressing disk of $A$, then
$D$ induces a disk $D'\subset\ComplHKA$ with 
$D'\cap l_+$ a point and $D'\cap A_-=\emptyset$;
the frontier of a regular neighborhood
$\rnbhd{A_+\cup D'}$ in $\ComplHKA-A_-$ 
is a separating disk $D''$ with $\partial D''$
parallel to $l_-$; this implies 
$A$ is compressible, that is, the previous case.

\ref{lm:esse_redu_type2_2:pairA} $\Rightarrow$ \ref{lm:esse_redu_type2_2:pair} \&
\ref{lm:esse_redu_type2_2:pairA} $\Rightarrow$ \ref{lm:esse_redu_type2_2:A} follow from Dehn's lemma. 
\end{proof}

\begin{lemma}\label{lm:triv_type2_2}
If $\pairA$ is trivial and $l_-$ 
is homotopically trivial, then $\pair$ is trivial.
\end{lemma}
\begin{proof}
Denote by $D\subset\ComplHKA$ a disk bounded by $l_-$. 
Then $D$ splits $\ComplHKA$
into two solid tori, in one of which $l_+$ is primitive. 
Therefore $\pi_1(\ComplHKA)$ has a basis $\{x,y\}$
with $x$ in the conjugacy class determined by $l_+$. 
The assertion then follows from the fact that
$\pi_1(\ComplHK)$ is the HNN extension of $\pi_1(\ComplHKA)$ with respect to
$\pi_1(A)$.
\end{proof}
The converse of Lemma \ref{lm:triv_type2_2} 
is not true in general.
%
As a corollary of Corollary \ref{cor:irre_atoro_vs_tri_atoro} and Lemmas \ref{lm:esse_redu_type2_1} and \ref{lm:esse_redu_type2_2},
we have the following.
\begin{corollary}\label{cor:ntri_atoro_esse_A}
If $\pair$ is non-trivial and atoroidal,
then $A$ is essential.
\end{corollary}

\subsection{Non-triviality and atoroidality}\label{subsec:tri_atoro_criteria}
We present here criteria for $\pair$
to be non-trivial and atoroidal in terms of 
$\pairA$ and $l_+,l_-$.  
Recall first two results on atoroidality:
\begin{corollary}\label{cor:atoro_HK_HKA}
If $\pair$ is non-trivial and atoroidal, then $\pairA$ is atoroidal.
\end{corollary}
\begin{proof}
This follows from \cite[Lemma $4.1$]{Wan:22},
but can also be deduced from Lemmas \ref{lm:esse_redu_type2_1}
and \ref{lm:esse_redu_type2_2}: since $A$ is essential by Corollary \ref{cor:ntri_atoro_esse_A}, if there exists an incompressible torus $T\subset \ComplHKA$, then any compressing disk of $T$ can be isotoped away from $A$, contradicting the atoroidality of $\pair$.
\end{proof}

\begin{corollary}{\cite[Lemma $4.9$]{Wan:22}}\label{cor:atoro_HKA_HK}
Suppose $\pairA$ is atoroidal, and $l_-\subset\ComplHKA$ 
is not homotopically trivial if $A$ is of type $2$-$2$. Then $\pair$ is atoroidal.
\end{corollary}

\begin{proposition}\label{prop:irre_atoro_type2_1}
Suppose $A$ is of type $2$-$1$.
Then $\pair$ is atoroidal and $A$ is essential 
if and only if $\pairA$ either is trivial 
with $\{l_+,l_-\}$ not primitive in $\ComplHKA$ or is non-trivial and atoroidal. 
\end{proposition}
\begin{proof}
``$\Rightarrow$'':
By Corollary \ref{cor:atoro_HK_HKA}, 
$\pairA$ is atoroidal. 
If $\pairA$ is trivial, then, 
since $A$ is essential, 
$\{l_+,l_-\}$ cannot be primitive in $\ComplHKA$ by Lemma \ref{lm:esse_redu_type2_1}.

``$\Leftarrow$'': $\pair$ is atoroidal by Corollary \ref{cor:atoro_HKA_HK}. $\pair$ cannot be trivial
by Lemma \ref{lm:triv_type2_1}, so 
$A$ is essential by Corollary \ref{cor:ntri_atoro_esse_A}. 
\end{proof}

\begin{proposition}\label{prop:irre_atoro_type2_2}
Suppose $A$ is of type $2$-$2$.
Then $\pair$ is atoroidal and $A$ is essential 
if and only if $\pairA$ either is trivial with $l_-\subset\ComplHKA$ not homotopically trivial 
or is non-trivial and atoroidal. 
\end{proposition}
\begin{proof}
``$\Rightarrow$'':
By Corollary \ref{cor:atoro_HK_HKA},
$\pairA$ is atoroidal. If $\pairA$ 
is trivial, then by Lemma \ref{lm:esse_redu_type2_2},
$l_-$ cannot be homotopically trivial since $A$ is essential.
 
``$\Leftarrow$':' By Lemma \ref{cor:atoro_HKA_HK} 
$\pair$ is atoroidal. If $A$ is inessential, then
by Lemma \ref{lm:esse_redu_type2_2}, $\pairA$ 
is reducible with $l_-\subset\ComplHKA$ homotopically trivial, contradicting the assumption and 
Lemma \ref{lm:tri_redu_toro}.
\end{proof}

\section{Examples}\label{sec:examples}
Here we construct atoroidal handlebody-knots
that admit a type $2$ essential annulus, and show
that annulus diagrams in Theorem \ref{teo:classification_char_diagram_type_two}
can all be realized by such handlebody-knots.

\subsection{Looping trivalent spatial graphs}
\begin{figure}[h]
\begin{subfigure}{0.47\textwidth}
\centering
\def\svgwidth{.5\columnwidth}
\begingroup%
  \makeatletter%
  \providecommand\color[2][]{%
    \errmessage{(Inkscape) Color is used for the text in Inkscape, but the package 'color.sty' is not loaded}%
    \renewcommand\color[2][]{}%
  }%
  \providecommand\transparent[1]{%
    \errmessage{(Inkscape) Transparency is used (non-zero) for the text in Inkscape, but the package 'transparent.sty' is not loaded}%
    \renewcommand\transparent[1]{}%
  }%
  \providecommand\rotatebox[2]{#2}%
  \newcommand*\fsize{\dimexpr\f@size pt\relax}%
  \newcommand*\lineheight[1]{\fontsize{\fsize}{#1\fsize}\selectfont}%
  \ifx\svgwidth\undefined%
    \setlength{\unitlength}{850.39370079bp}%
    \ifx\svgscale\undefined%
      \relax%
    \else%
      \setlength{\unitlength}{\unitlength * \real{\svgscale}}%
    \fi%
  \else%
    \setlength{\unitlength}{\svgwidth}%
  \fi%
  \global\let\svgwidth\undefined%
  \global\let\svgscale\undefined%
  \makeatother%
  \begin{picture}(1,0.83333333)%
    \lineheight{1}%
    \setlength\tabcolsep{0pt}%
    \put(0,0){\includegraphics[width=\unitlength,page=1]{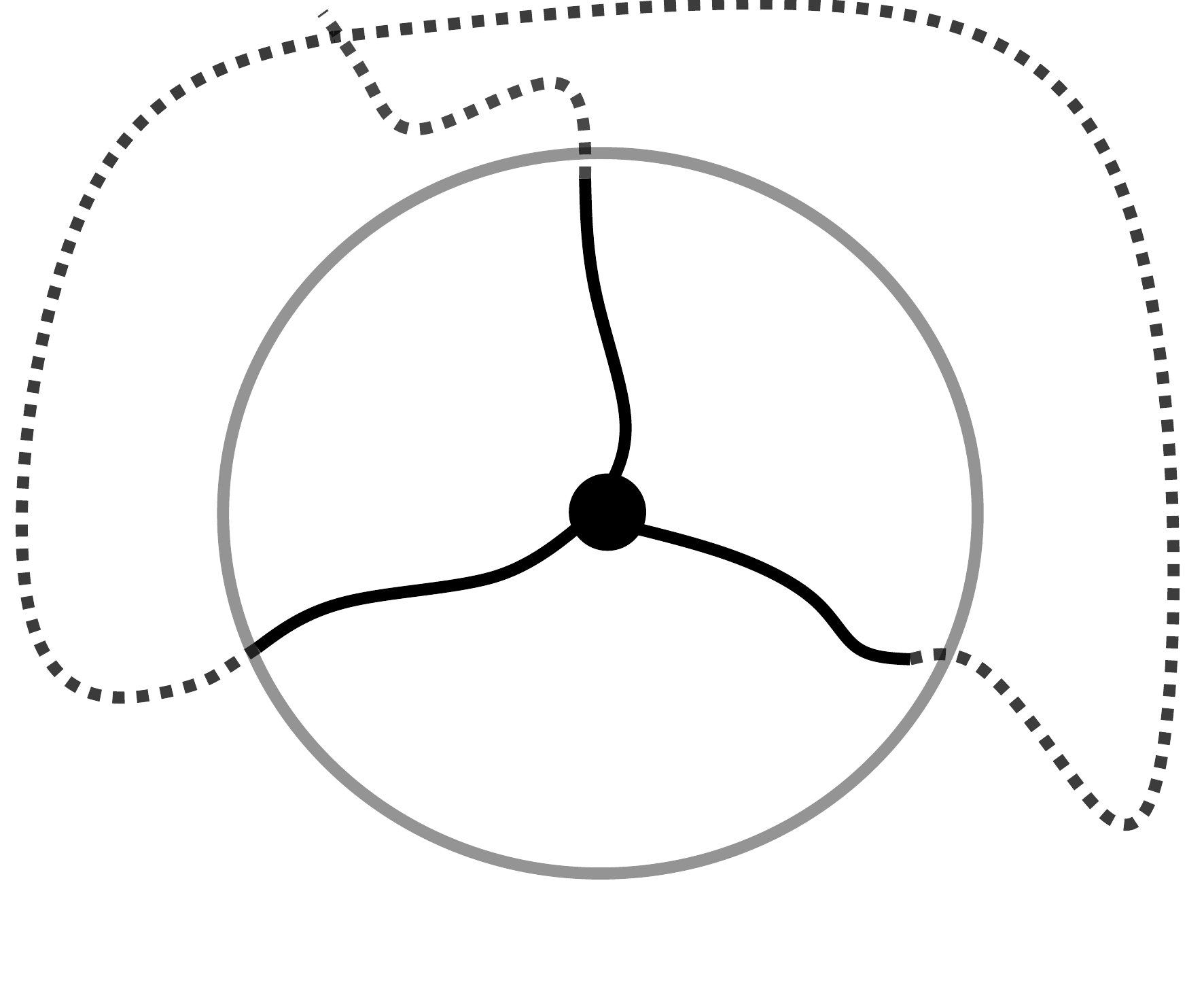}}%
    \put(0.43519694,0.31342651){\color[rgb]{0,0,0}\makebox(0,0)[lt]{\lineheight{1.25}\smash{\begin{tabular}[t]{l}{\footnotesize $v$}\end{tabular}}}}%
    \put(0.06236734,0.18185852){\color[rgb]{0,0,0}\makebox(0,0)[lt]{\lineheight{1.25}\smash{\begin{tabular}[t]{l}{\footnotesize $e_1$}\end{tabular}}}}%
    \put(0.49121104,0.74623805){\color[rgb]{0,0,0}\makebox(0,0)[lt]{\lineheight{1.25}\smash{\begin{tabular}[t]{l}{\footnotesize $e_2$}\end{tabular}}}}%
    \put(0.86429006,0.25280304){\color[rgb]{0,0,0}\makebox(0,0)[lt]{\lineheight{1.25}\smash{\begin{tabular}[t]{l}{\footnotesize $e_3$}\end{tabular}}}}%
    \put(0,0){\includegraphics[width=\unitlength,page=2]{looping_theta1.pdf}}%
  \end{picture}%
\endgroup%

\caption{Neighborhood of a trivalent node $v\in\SG$.}
\label{fig:neightborhood_trivalent_node}
\end{subfigure}
\raisebox{.4cm}{$\Longrightarrow$}
\begin{subfigure}{0.47\textwidth}
\centering
\def\svgwidth{.5 \columnwidth}
\begingroup%
  \makeatletter%
  \providecommand\color[2][]{%
    \errmessage{(Inkscape) Color is used for the text in Inkscape, but the package 'color.sty' is not loaded}%
    \renewcommand\color[2][]{}%
  }%
  \providecommand\transparent[1]{%
    \errmessage{(Inkscape) Transparency is used (non-zero) for the text in Inkscape, but the package 'transparent.sty' is not loaded}%
    \renewcommand\transparent[1]{}%
  }%
  \providecommand\rotatebox[2]{#2}%
  \newcommand*\fsize{\dimexpr\f@size pt\relax}%
  \newcommand*\lineheight[1]{\fontsize{\fsize}{#1\fsize}\selectfont}%
  \ifx\svgwidth\undefined%
    \setlength{\unitlength}{850.39370079bp}%
    \ifx\svgscale\undefined%
      \relax%
    \else%
      \setlength{\unitlength}{\unitlength * \real{\svgscale}}%
    \fi%
  \else%
    \setlength{\unitlength}{\svgwidth}%
  \fi%
  \global\let\svgwidth\undefined%
  \global\let\svgscale\undefined%
  \makeatother%
  \begin{picture}(1,0.83333333)%
    \lineheight{1}%
    \setlength\tabcolsep{0pt}%
    \put(0,0){\includegraphics[width=\unitlength,page=1]{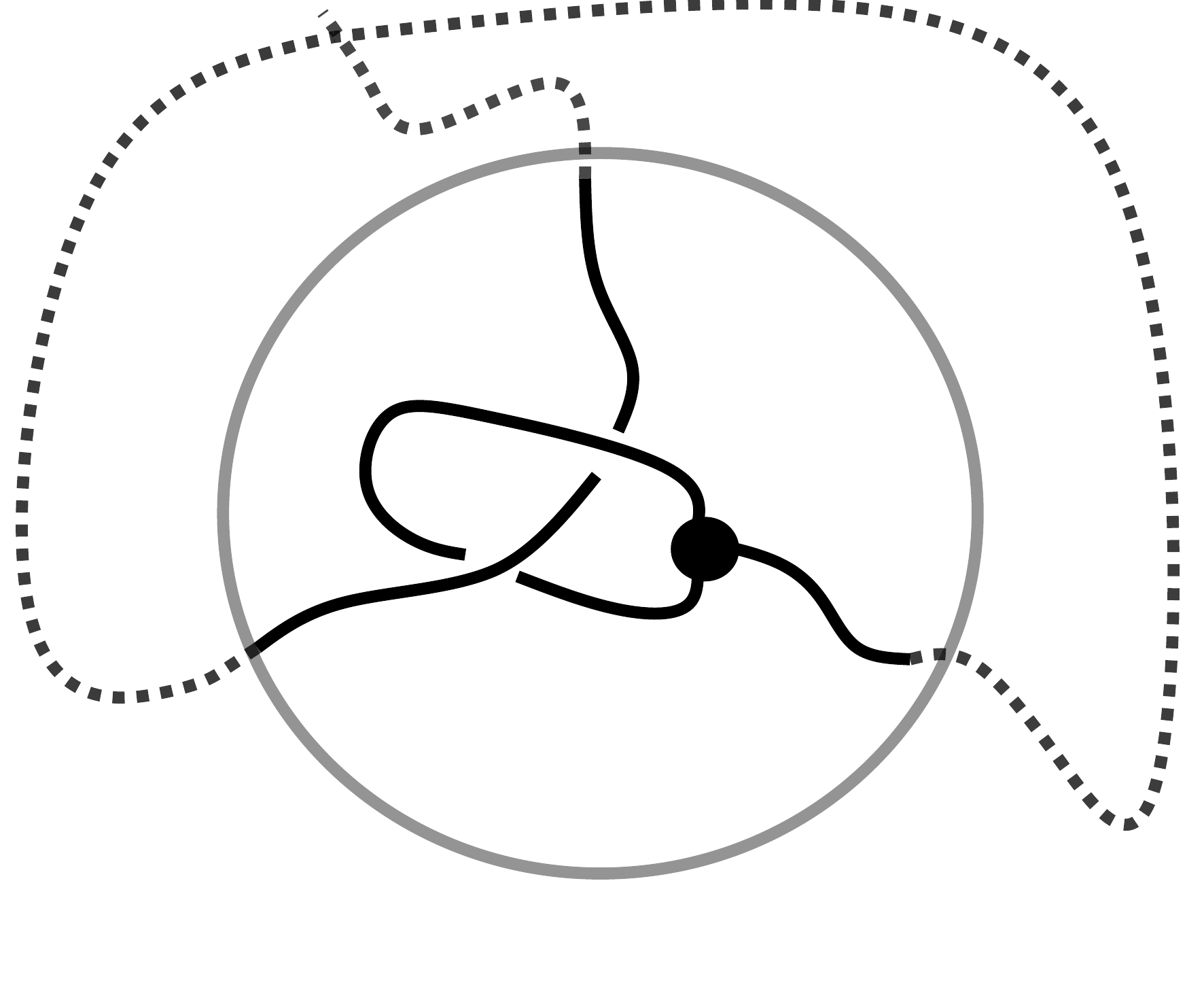}}%
    \put(0.06236734,0.18185842){\color[rgb]{0,0,0}\makebox(0,0)[lt]{\lineheight{1.25}\smash{\begin{tabular}[t]{l}{\footnotesize $e_1$}\end{tabular}}}}%
    \put(0.49121099,0.74623805){\color[rgb]{0,0,0}\makebox(0,0)[lt]{\lineheight{1.25}\smash{\begin{tabular}[t]{l}{\footnotesize $e_2$}\end{tabular}}}}%
    \put(0.86429006,0.25280299){\color[rgb]{0,0,0}\makebox(0,0)[lt]{\lineheight{1.25}\smash{\begin{tabular}[t]{l}{\footnotesize $e_3$}\end{tabular}}}}%
    \put(0,0){\includegraphics[width=\unitlength,page=2]{looping_theta2.pdf}}%
  \end{picture}%
\endgroup%

\caption{Replacing $v$ with a loop.}
\label{fig:looping_e_1_e_2_at_v}
\end{subfigure}
\caption{Looping of a spatial $\theta$-graph.}
\label{fig:looping_theta}
\end{figure}

Let $\pairSG$ be a spatial graph with $\SG$ 
either a $\theta$-graph or a handcuff graph. 
Then we can produce a new spatial graph $\pairLSG$
by replacing a small neighborhood of a 
trivalent node\footnote{A neighborhood
$\rnbhd{v}\in\SG$ of the trivalent node $v$ is 
a regular neighborhood of $v\subset\sphere$ such that 
$(\rnbhd{v},\rnbhd{v}\cap \SG)$ is homeomorphic to 
a unit $3$-ball with three non-negative axes.}
in $\SG$ with a loop 
as shown in Fig.\ \ref{fig:looping_theta}.

\begin{figure}[h]
\begin{subfigure}{0.47\textwidth}
\centering
\def\svgwidth{.5\columnwidth}
\begingroup%
  \makeatletter%
  \providecommand\color[2][]{%
    \errmessage{(Inkscape) Color is used for the text in Inkscape, but the package 'color.sty' is not loaded}%
    \renewcommand\color[2][]{}%
  }%
  \providecommand\transparent[1]{%
    \errmessage{(Inkscape) Transparency is used (non-zero) for the text in Inkscape, but the package 'transparent.sty' is not loaded}%
    \renewcommand\transparent[1]{}%
  }%
  \providecommand\rotatebox[2]{#2}%
  \newcommand*\fsize{\dimexpr\f@size pt\relax}%
  \newcommand*\lineheight[1]{\fontsize{\fsize}{#1\fsize}\selectfont}%
  \ifx\svgwidth\undefined%
    \setlength{\unitlength}{850.39370079bp}%
    \ifx\svgscale\undefined%
      \relax%
    \else%
      \setlength{\unitlength}{\unitlength * \real{\svgscale}}%
    \fi%
  \else%
    \setlength{\unitlength}{\svgwidth}%
  \fi%
  \global\let\svgwidth\undefined%
  \global\let\svgscale\undefined%
  \makeatother%
  \begin{picture}(1,0.83333333)%
    \lineheight{1}%
    \setlength\tabcolsep{0pt}%
    \put(0,0){\includegraphics[width=\unitlength,page=1]{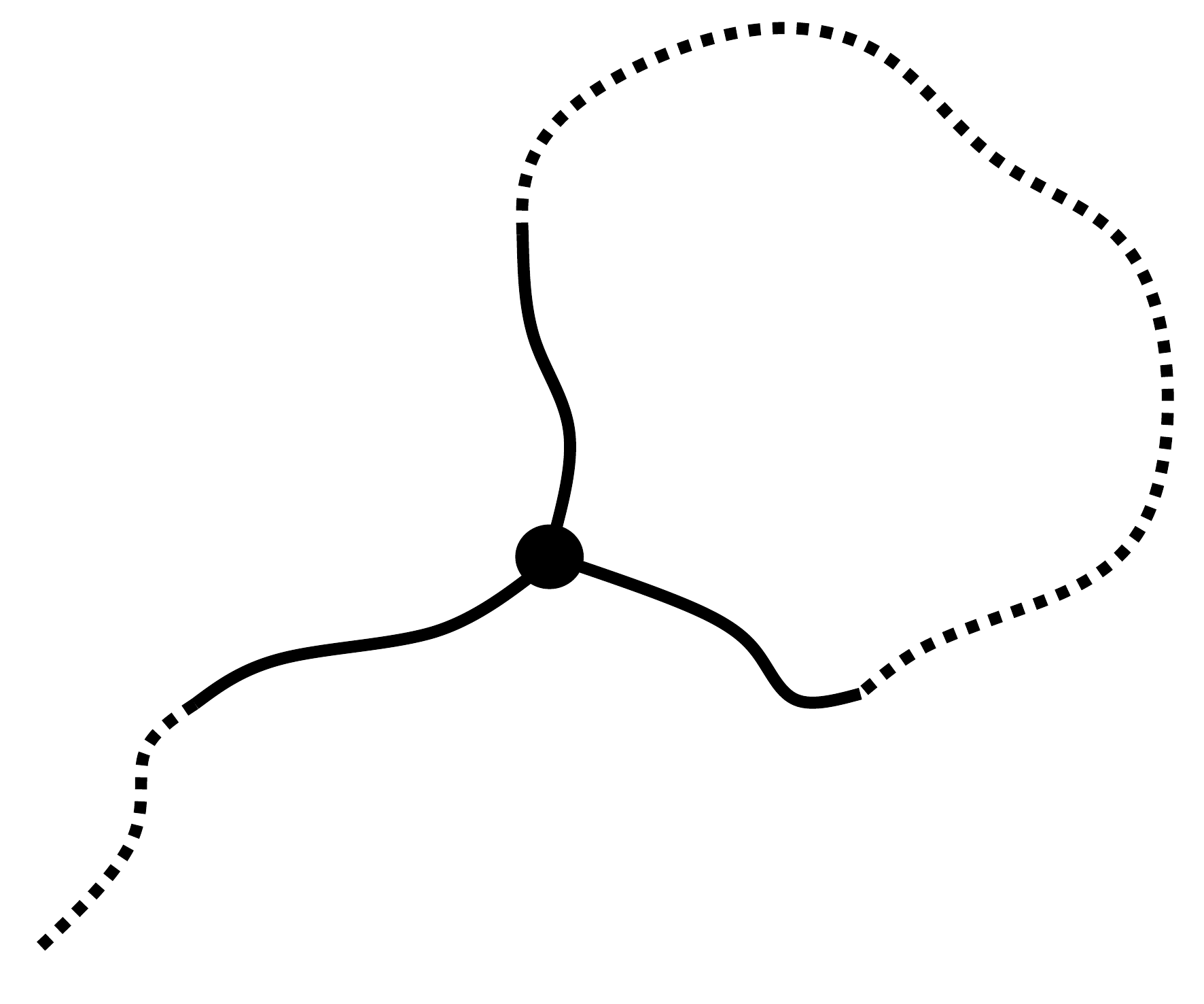}}%
    \put(0.40766624,0.2869547){\color[rgb]{0,0,0}\makebox(0,0)[lt]{\lineheight{1.25}\smash{\begin{tabular}[t]{l}{\footnotesize $v$}\end{tabular}}}}%
    \put(0,0){\includegraphics[width=\unitlength,page=2]{looping_handcuff1.pdf}}%
  \end{picture}%
\endgroup%

\caption{Neighborhood of a trivalent node $v\in\SG$.}
\end{subfigure}
\raisebox{.4cm}{$\Longrightarrow$}
\begin{subfigure}{0.47\textwidth}
\centering
\def\svgwidth{.5 \columnwidth}
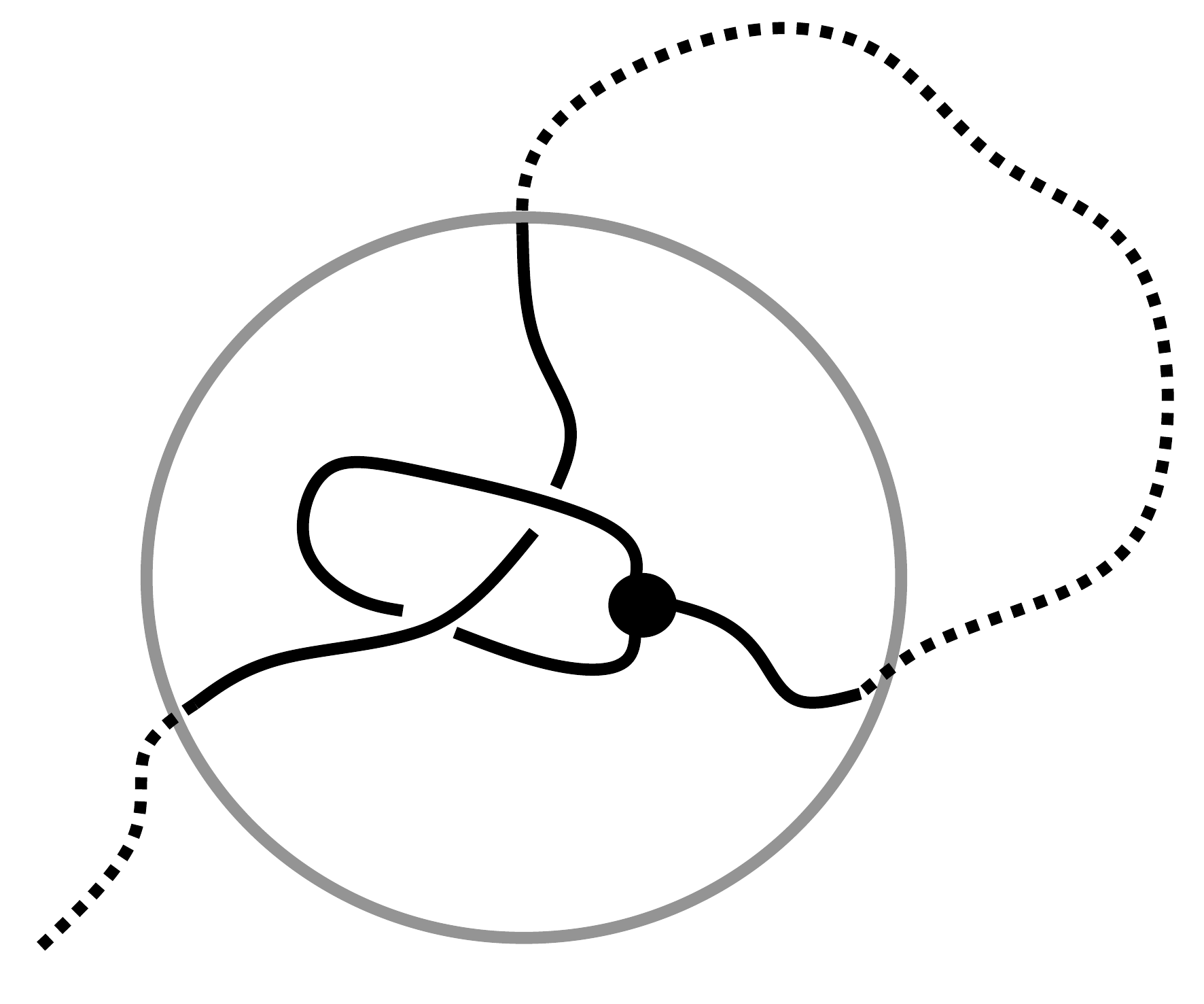
\caption{Replacing $v$ with a loop.}
\end{subfigure}
\caption{Looping of a spatial handcuff graph.}
\label{fig:looping_handcuff}
\end{figure}

Label the trivalent node with $v$ and its three adjacent edges $e_1,e_2,e_3$ as in Fig.\ \ref{fig:looping_theta}. 
Then the new spatial graph $\pairLSG$ in Fig.\ \ref{fig:looping_e_1_e_2_at_v}
is said to be obtained by \emph{looping $e_1e_2$ at $v$}. $\pairLSG$ is called a \emph{looping} of $\pairSG$,
provided the resulting spatial graph is connected (see Fig.\ \ref{fig:looping_handcuff}); there are six possible loopings for a spatial $\theta$-graph, and
four for a spatial handcuff graph. 
 
A double looping $\pairDLSG$ of $\pairSG$ 
is the spatial graph obtained by looping
at both trivalent nodes of $\SG$. 
Taking a regular neighborhood
of a looping $\LSG$ (resp.\ double looping $\DLSG$)
in $\sphere$ gives us a handlebody-knot, denoted by $\pairLHK$ (resp.\ $\pairDLHK$), whose exterior contains 
a canonical type $2$ annulus $\LA$ induced 
by the created loop
in $\pairLSG$.

A spatial graph $\pairSG$ is said to be \emph{nontrivially atoroidal} if the induced handlebody-knot $\pairNSG$
is non-trivial and atoroidal.  
\begin{lemma}\label{lm:looping}
If $\pairSG$ is nontrivially atoroidal,
then 
$\pairLHK$ induced by 
a looping of $\pairSG$
is atoroidal, and $\LA\subset \Compl\LHK$ 
is essential. Furthermore $\LA$ is 
of type $2$-$1$ and is the unique annulus 
if $\SG$ is a $\theta$-graph,
and is of type $2$-$2$ if $\SG$ 
is a handcuff graph. 
\end{lemma}
\begin{proof}
The disk bounded by a component of $\partial \LA$ in $\LHK$
is dual to the two edges being looped, so
$\LA$ is of type $2$-$1$ if $\Gamma$ is a $\theta$-graph
and is of type $2$-$2$ otherwise.
The essentiality of $\LA$ and atoroidality of 
$\pairLHK$ follow from Propositions \ref{prop:irre_atoro_type2_1}
and \ref{prop:irre_atoro_type2_2}. 
\end{proof}

\begin{corollary}\label{cor:double_looping}
If $\pairSG$ is nontrivially atoroidal,
then any handlebody-knot 
$\pairDLHK$ obtained by 
a double looping of $\pairSG$
is atoroidal, and its exterior
contains two non-isotopic type $2$-$2$ essential annuli.
\end{corollary}
\begin{proof}
The two canonical annuli are of type $2$-$2$ 
since any looping $\pairLSG$ 
is a spatial handcuff graph. The rest follows from Lemma \ref{lm:looping}.
\end{proof}
 
\begin{figure}[h]
\begin{subfigure}{0.32\textwidth}
\centering
\includegraphics[scale=.08]{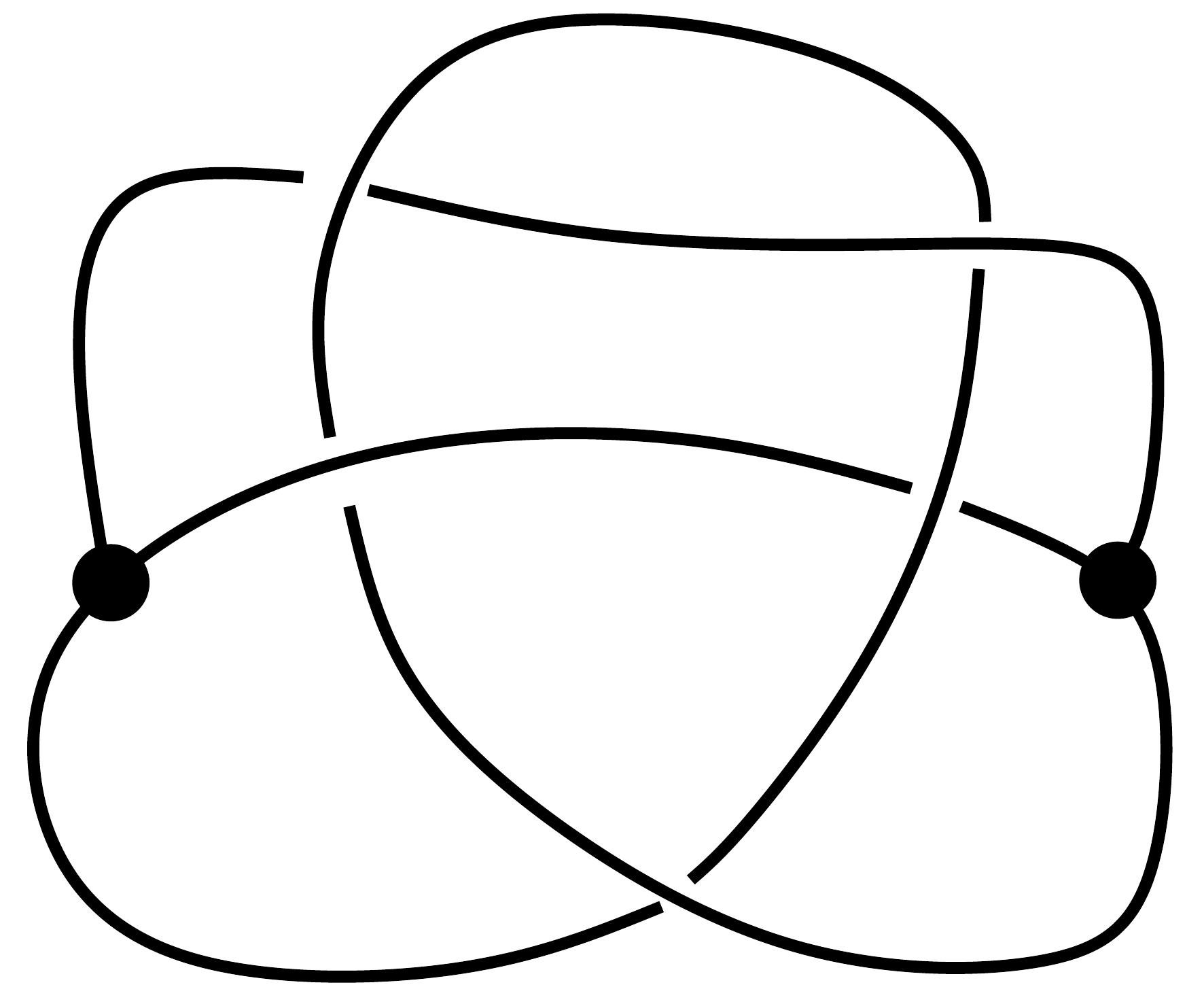}
\caption{Spine of $(\sphere,5_2)$.}
\label{fig:hk5_2}
\end{subfigure} 
\begin{subfigure}{0.32\textwidth}
\centering
\includegraphics[scale=.08]{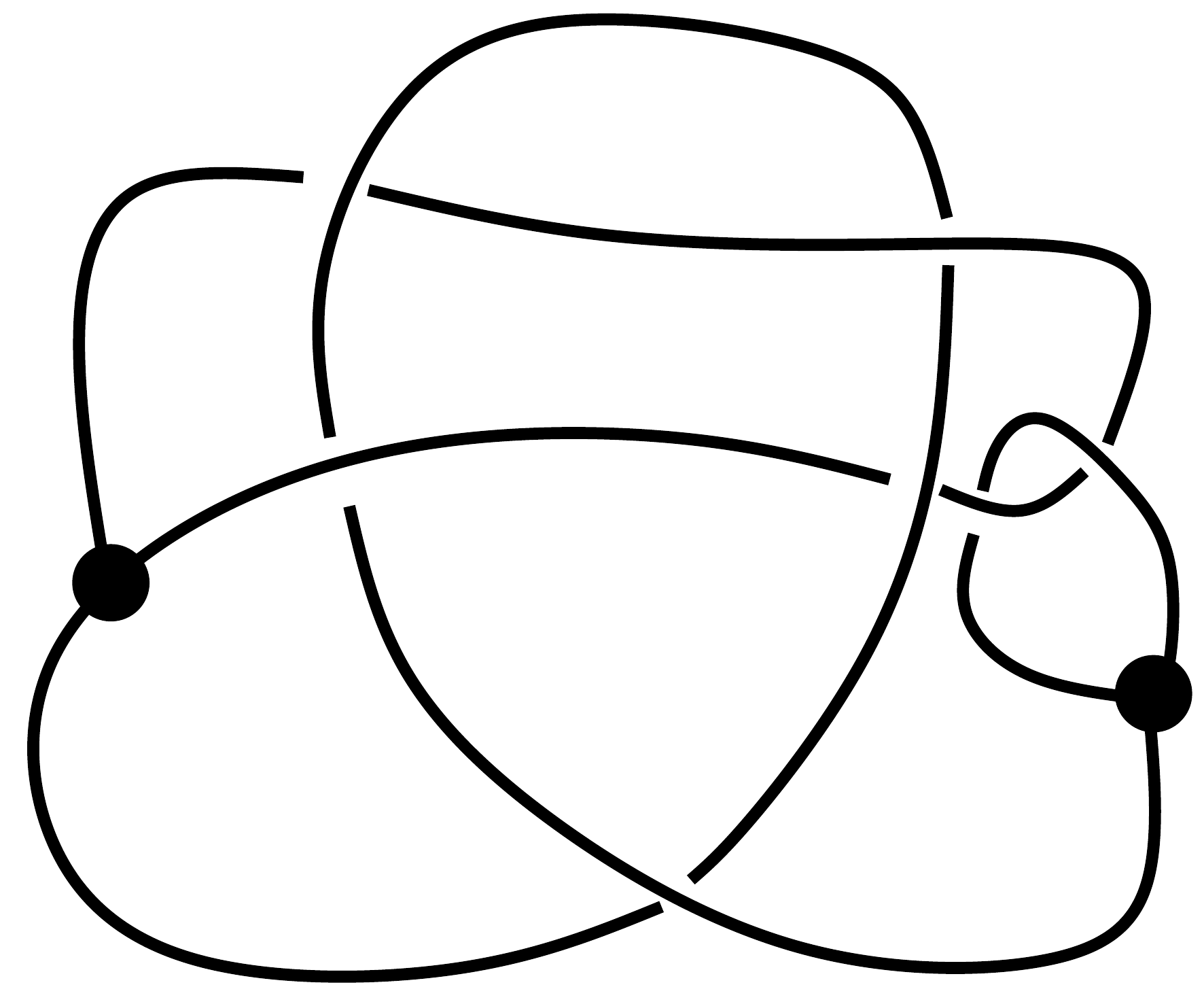}
\caption{Looping.}
\label{fig:looping_hk5_2}
\end{subfigure}
\begin{subfigure}{0.32\textwidth}
\centering
\includegraphics[scale=.08]{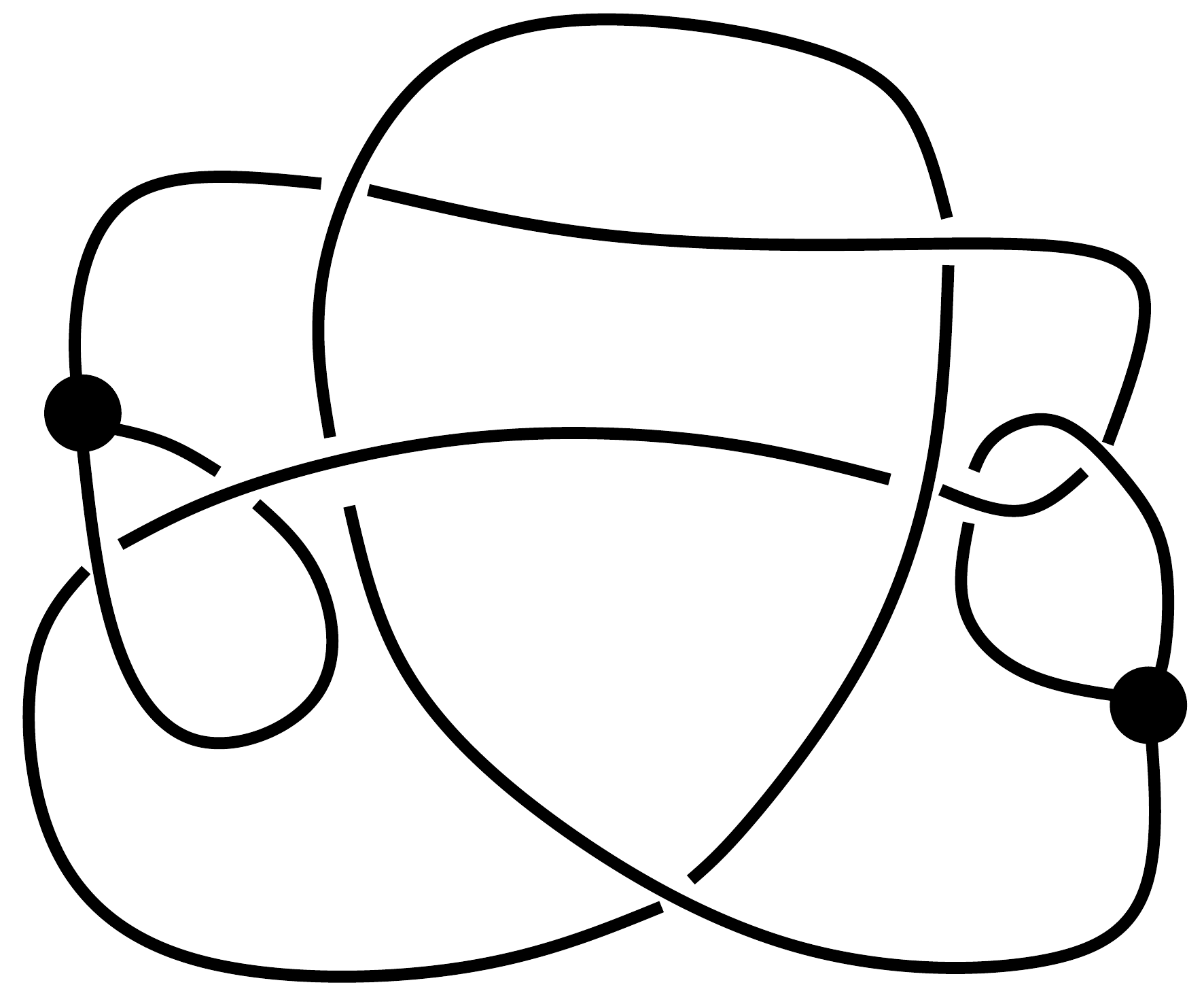}
\caption{Double looping.}
\label{fig:double_looping_hk5_2}
\end{subfigure}
\caption{Handlebody-knots with a type $2$ annulus.}
\end{figure}

As an application of Lemma \ref{lm:looping} and Corollary \ref{cor:double_looping}, 
we consider the spine $\pairSG$ of $(\sphere,5_2)$ in \cite{IshKisMorSuz:12} as shown in Fig.\ \ref{fig:hk5_2}. Then Fig.\ \ref{fig:looping_hk5_2} is 
a looping of $\pairSG$, 
whose associated handlebody-knot has the annulus diagram 
\raisebox{-.25cm}{\includegraphics[scale=.13]{typetwoone_ann}}.
On the other hand,  the double looping of $\pairSG$ in Fig.\ \ref{fig:double_looping_hk5_2}  
induces a handlebody-knot whose annulus diagram is 
\raisebox{-.2cm}{\includegraphics[scale=.12]{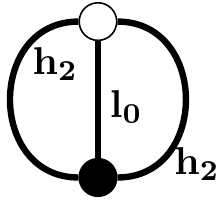}}.

\subsection{Unknotting annuli of type $2$}
As opposed to Lemma \ref{lm:looping} and Corollary \ref{cor:double_looping}, 
here we present a looping operation that yields
atoroidal handlebody-knots that admit an essential \emph{unknotting} type $2$ annulus.

Let $\pairSG$ be a spatial $\theta$-graph 
that is a union of a non-trivial knot $\pairK$ and a tunnel $\tau$
of $\pairK$. Let $\kappa_1,\kappa_2$ be the arcs of $K$
cut off by $\tau$.
Then a \emph{tunnel} looping of $\pairKtau$
is a looping obtained by looping $k_i\tau$
at a trivalent node of $\SG=K\cup\tau$, $i=1$ or $2$.
 

\begin{lemma}\label{lm:looping_nontrivial_knot_tunnel}
The handlebody-knot $\pairLHK$ induced by 
a tunnel looping of $\pairSG$ 
is atoroidal, and 
$\LA$ is an unknotting essential type $2$-$1$ annulus. 
\end{lemma}
\begin{proof}
It follows from the ``only if '' part of  
Proposition \ref{prop:irre_atoro_type2_1}
since $\pairK$ is non-trivial.
\end{proof}

Now, let 
$\pairSG$ be the union of a non-split link
$\pairL$ and a tunnel $\tau$ of $\pairL$.  
\begin{lemma}\label{lm:looping_nonsplit_link_tunnel}
The handlebody-knot $\pairLHK$ induced by 
a looping of $\pairSG$ 
is atoroidal, and
$\LA$ is an unknotting essential type $2$-$2$ annulus. 
\end{lemma}
\begin{proof}
Use $\pairL$ being non-split and apply the ``only if'' part of Proposition \ref{prop:irre_atoro_type2_2}. 
\end{proof}

\begin{figure}[h]
\begin{subfigure}{0.47\textwidth}
\centering
\def\svgwidth{.6\columnwidth}
\begingroup%
  \makeatletter%
  \providecommand\color[2][]{%
    \errmessage{(Inkscape) Color is used for the text in Inkscape, but the package 'color.sty' is not loaded}%
    \renewcommand\color[2][]{}%
  }%
  \providecommand\transparent[1]{%
    \errmessage{(Inkscape) Transparency is used (non-zero) for the text in Inkscape, but the package 'transparent.sty' is not loaded}%
    \renewcommand\transparent[1]{}%
  }%
  \providecommand\rotatebox[2]{#2}%
  \newcommand*\fsize{\dimexpr\f@size pt\relax}%
  \newcommand*\lineheight[1]{\fontsize{\fsize}{#1\fsize}\selectfont}%
  \ifx\svgwidth\undefined%
    \setlength{\unitlength}{850.39370079bp}%
    \ifx\svgscale\undefined%
      \relax%
    \else%
      \setlength{\unitlength}{\unitlength * \real{\svgscale}}%
    \fi%
  \else%
    \setlength{\unitlength}{\svgwidth}%
  \fi%
  \global\let\svgwidth\undefined%
  \global\let\svgscale\undefined%
  \makeatother%
  \begin{picture}(1,0.83333333)%
    \lineheight{1}%
    \setlength\tabcolsep{0pt}%
    \put(0,0){\includegraphics[width=\unitlength,page=1]{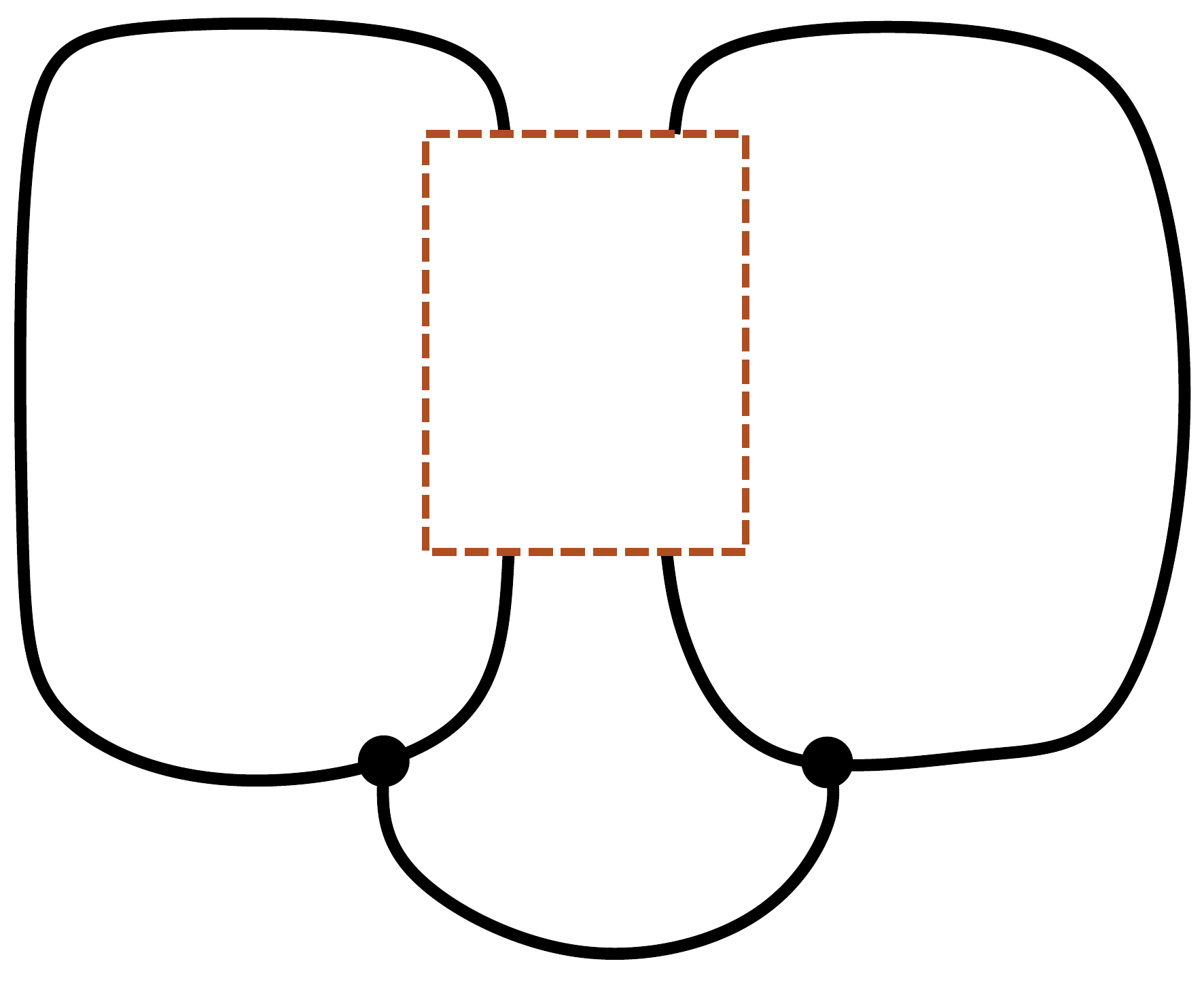}}%
    \put(0.42143155,0.59508682){\color[rgb]{0,0,0}\makebox(0,0)[lt]{\lineheight{1.25}\smash{\begin{tabular}[t]{l}{\footnotesize$n$}\end{tabular}}}}%
    \put(0.40167704,0.5344634){\color[rgb]{0,0,0}\makebox(0,0)[lt]{\lineheight{1.25}\smash{\begin{tabular}[t]{l}{\footnotesize half}\end{tabular}}}}%
    \put(0.38579389,0.47516647){\color[rgb]{0,0,0}\makebox(0,0)[lt]{\lineheight{1.25}\smash{\begin{tabular}[t]{l}{\footnotesize twists}\end{tabular}}}}%
    \put(0.09423973,0.71368065){\color[rgb]{0,0,0}\makebox(0,0)[lt]{\lineheight{1.25}\smash{\begin{tabular}[t]{l}{\footnotesize $l_1$}\end{tabular}}}}%
    \put(0.68576497,0.7204184){\color[rgb]{0,0,0}\makebox(0,0)[lt]{\lineheight{1.25}\smash{\begin{tabular}[t]{l}{\footnotesize $l_2$}\end{tabular}}}}%
    \put(0.46340159,0.06709351){\color[rgb]{0,0,0}\makebox(0,0)[lt]{\lineheight{1.25}\smash{\begin{tabular}[t]{l}{\footnotesize $\tau$}\end{tabular}}}}%
  \end{picture}%
\endgroup%

\caption{$(n,2)$-torus link $(\sphere,L_n)$ 
with a tunnel $\tau$.}
\label{fig:n_2_torus_link_tunnel}
\end{subfigure}
\raisebox{.4cm}{$\Longrightarrow$}
\begin{subfigure}{0.47\textwidth}
\centering
\def\svgwidth{.6 \columnwidth}
\begingroup%
  \makeatletter%
  \providecommand\color[2][]{%
    \errmessage{(Inkscape) Color is used for the text in Inkscape, but the package 'color.sty' is not loaded}%
    \renewcommand\color[2][]{}%
  }%
  \providecommand\transparent[1]{%
    \errmessage{(Inkscape) Transparency is used (non-zero) for the text in Inkscape, but the package 'transparent.sty' is not loaded}%
    \renewcommand\transparent[1]{}%
  }%
  \providecommand\rotatebox[2]{#2}%
  \newcommand*\fsize{\dimexpr\f@size pt\relax}%
  \newcommand*\lineheight[1]{\fontsize{\fsize}{#1\fsize}\selectfont}%
  \ifx\svgwidth\undefined%
    \setlength{\unitlength}{850.39370079bp}%
    \ifx\svgscale\undefined%
      \relax%
    \else%
      \setlength{\unitlength}{\unitlength * \real{\svgscale}}%
    \fi%
  \else%
    \setlength{\unitlength}{\svgwidth}%
  \fi%
  \global\let\svgwidth\undefined%
  \global\let\svgscale\undefined%
  \makeatother%
  \begin{picture}(1,0.83333333)%
    \lineheight{1}%
    \setlength\tabcolsep{0pt}%
    \put(0,0){\includegraphics[width=\unitlength,page=1]{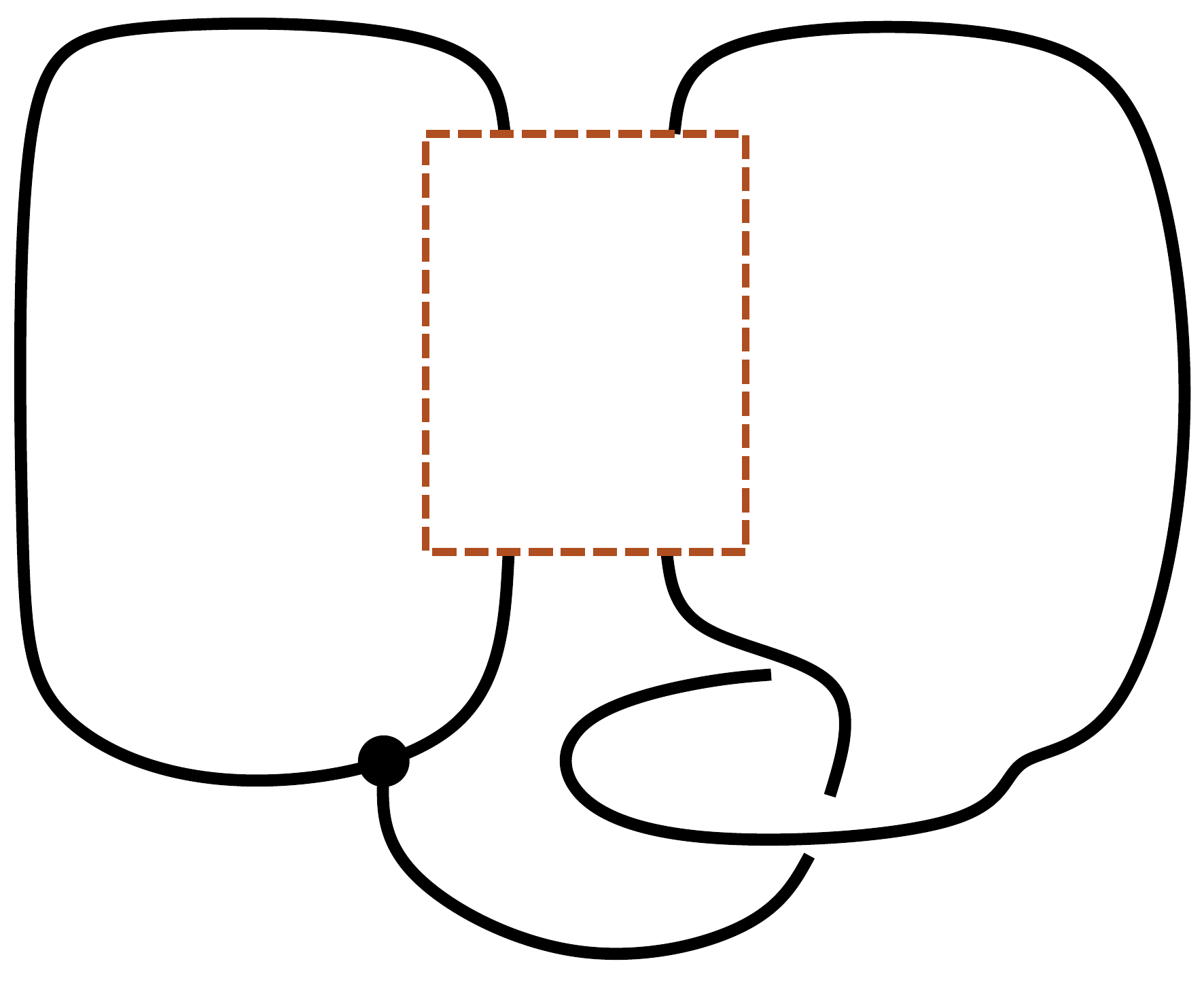}}%
    \put(0.42143155,0.59508682){\color[rgb]{0,0,0}\makebox(0,0)[lt]{\lineheight{1.25}\smash{\begin{tabular}[t]{l}{\footnotesize $n$}\end{tabular}}}}%
    \put(0.39744151,0.53552226){\color[rgb]{0,0,0}\makebox(0,0)[lt]{\lineheight{1.25}\smash{\begin{tabular}[t]{l}{\footnotesize half}\end{tabular}}}}%
    \put(0.37838181,0.47516647){\color[rgb]{0,0,0}\makebox(0,0)[lt]{\lineheight{1.25}\smash{\begin{tabular}[t]{l}{\footnotesize twists}\end{tabular}}}}%
    \put(0,0){\includegraphics[width=\unitlength,page=2]{looping_n_2_torus_link.pdf}}%
    \put(0.51990677,0.1747141){\color[rgb]{0,0,0}\makebox(0,0)[lt]{\lineheight{1.25}\smash{\begin{tabular}[t]{l}{\footnotesize $l$}\end{tabular}}}}%
  \end{picture}%
\endgroup%

\caption{(Tunnel) looping of $(\sphere,L_n\cup\tau)$.}
\label{fig:looping_n_2_torus_link}
\end{subfigure}
\caption{Construction of Koda's handlebody-knot family.}
\end{figure}

To show that all annulus diagrams in Theorem \ref{teo:classification_char_diagram_type_two} can
be realized by some atoroidal handlebody-knots, 
we consider the union of an $(n,2)$-torus link $(\sphere,L_n=l_1\cup l_2)$,
$n\in\mathbb{Z}$, with a tunnel $\tau$ as
depicted in Fig.\ \ref{fig:n_2_torus_link_tunnel}.
Denote by $\pairn$ the handlebody-knot induced by 
the looping of $\pairLntau$ in Fig.\ \ref{fig:looping_n_2_torus_link}.
Note that $(\sphere,\HK_2)$ is equivalent to 
$\pairfourone$, while $\{(\sphere,\HK_n)\}_{n>2}$ is Koda's handlebody-knot family in
\cite[Example $3$]{Kod:15}; Lemmas \ref{lm:looping_nontrivial_knot_tunnel} and \ref{lm:looping_nonsplit_link_tunnel} 
give an alternative way to see they are irreducible, in view of Corollary \ref{cor:irre_atoro_vs_tri_atoro}.

Observe that if $n>2$ and is even,
the handlebody-knot exterior $\Compl{\HK_n}$
contains a type $3$-$2$ annulus $A$ given as follows:
let $A_c$ be a cabling annulus in $\Compl {L_n}:=\sphere-\openrnbhd{L_n}$
with $\tau\cap \Compl{L_n}\subset A_c$. 
Let 
$\rnbhd{l_i}$ be the component of $\rnbhd{L_n}$ containing
$l_i$, $i=1,2$, and perform the looping construction 
entirely in
$\openrnbhd{l_2}$. Then the frontier of $\rnbhd{l_2}\cup \rnbhd{A_c}$ in $\Compl{l_1}:=\sphere-\openrnbhd{l_1}$ is an essential annulus $A\subset\Compl{\HK_n}$ of type $3$-$2$ii as $A$ is $\partial$-compressible in $\Compl{l_1}$. 

%
%
%
\begin{corollary}
Suppose $n>2$ and is even.
Then the annulus diagram of 
the handlebody-knot $(\sphere,\HK_n)$ 
obtained by the looping of $\pairLntau$ in Fig.\ \ref{fig:looping_n_2_torus_link} is 
\raisebox{-.2cm}{\includegraphics[scale=.14]{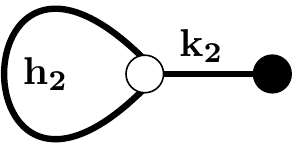}}.
\end{corollary}
\begin{figure}[h]
\begin{subfigure}{0.47\textwidth}
\centering
\def\svgwidth{.6\columnwidth}
\begingroup%
  \makeatletter%
  \providecommand\color[2][]{%
    \errmessage{(Inkscape) Color is used for the text in Inkscape, but the package 'color.sty' is not loaded}%
    \renewcommand\color[2][]{}%
  }%
  \providecommand\transparent[1]{%
    \errmessage{(Inkscape) Transparency is used (non-zero) for the text in Inkscape, but the package 'transparent.sty' is not loaded}%
    \renewcommand\transparent[1]{}%
  }%
  \providecommand\rotatebox[2]{#2}%
  \newcommand*\fsize{\dimexpr\f@size pt\relax}%
  \newcommand*\lineheight[1]{\fontsize{\fsize}{#1\fsize}\selectfont}%
  \ifx\svgwidth\undefined%
    \setlength{\unitlength}{850.39370079bp}%
    \ifx\svgscale\undefined%
      \relax%
    \else%
      \setlength{\unitlength}{\unitlength * \real{\svgscale}}%
    \fi%
  \else%
    \setlength{\unitlength}{\svgwidth}%
  \fi%
  \global\let\svgwidth\undefined%
  \global\let\svgscale\undefined%
  \makeatother%
  \begin{picture}(1,0.83333333)%
    \lineheight{1}%
    \setlength\tabcolsep{0pt}%
    \put(0,0){\includegraphics[width=\unitlength,page=1]{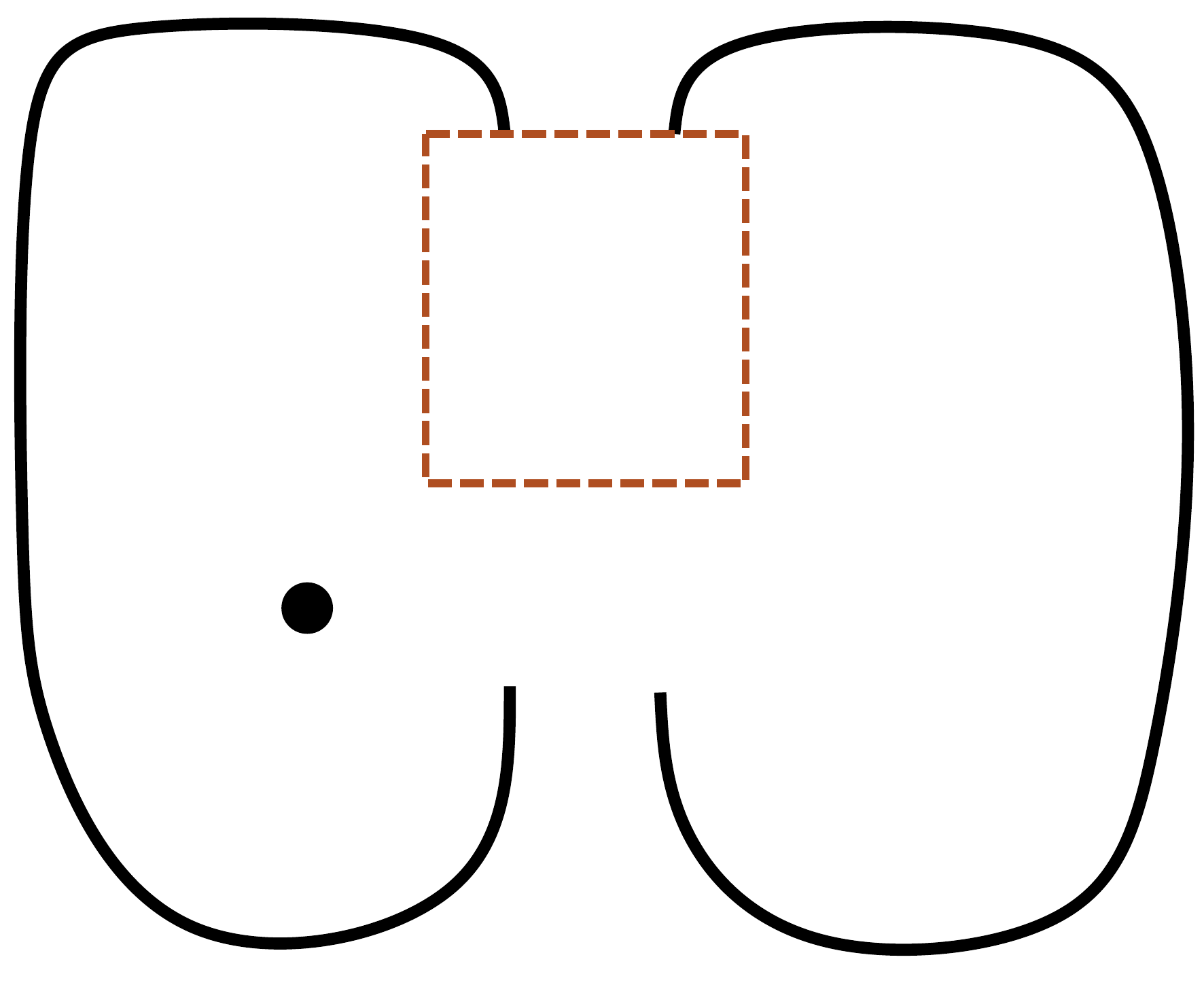}}%
    \put(0.41296056,0.63532405){\color[rgb]{0,0,0}\makebox(0,0)[lt]{\lineheight{1.25}\smash{\begin{tabular}[t]{l}{\footnotesize$n$}\end{tabular}}}}%
    \put(0.39320602,0.5747006){\color[rgb]{0,0,0}\makebox(0,0)[lt]{\lineheight{1.25}\smash{\begin{tabular}[t]{l}{\footnotesize half}\end{tabular}}}}%
    \put(0.3773229,0.51540365){\color[rgb]{0,0,0}\makebox(0,0)[lt]{\lineheight{1.25}\smash{\begin{tabular}[t]{l}{\footnotesize twists}\end{tabular}}}}%
    \put(0.11715,0.40699188){\color[rgb]{0,0,0}\makebox(0,0)[lt]{\lineheight{1.25}\smash{\begin{tabular}[t]{l}{\footnotesize $\tau$}\end{tabular}}}}%
    \put(0,0){\includegraphics[width=\unitlength,page=2]{ringed_n_2_torus_knot_tunnel.pdf}}%
  \end{picture}%
\endgroup%

\caption{Union of $(\sphere,L_n)$ with $n$ odd and a tunnel $\tau$.}
\label{fig:ringed_n_2_torus_knot_tunnel}
\end{subfigure}
\raisebox{.4cm}{$\Longrightarrow$}
\begin{subfigure}{0.47\textwidth}
\centering
\def\svgwidth{.6 \columnwidth}
\begingroup%
  \makeatletter%
  \providecommand\color[2][]{%
    \errmessage{(Inkscape) Color is used for the text in Inkscape, but the package 'color.sty' is not loaded}%
    \renewcommand\color[2][]{}%
  }%
  \providecommand\transparent[1]{%
    \errmessage{(Inkscape) Transparency is used (non-zero) for the text in Inkscape, but the package 'transparent.sty' is not loaded}%
    \renewcommand\transparent[1]{}%
  }%
  \providecommand\rotatebox[2]{#2}%
  \newcommand*\fsize{\dimexpr\f@size pt\relax}%
  \newcommand*\lineheight[1]{\fontsize{\fsize}{#1\fsize}\selectfont}%
  \ifx\svgwidth\undefined%
    \setlength{\unitlength}{850.39370079bp}%
    \ifx\svgscale\undefined%
      \relax%
    \else%
      \setlength{\unitlength}{\unitlength * \real{\svgscale}}%
    \fi%
  \else%
    \setlength{\unitlength}{\svgwidth}%
  \fi%
  \global\let\svgwidth\undefined%
  \global\let\svgscale\undefined%
  \makeatother%
  \begin{picture}(1,0.83333333)%
    \lineheight{1}%
    \setlength\tabcolsep{0pt}%
    \put(0,0){\includegraphics[width=\unitlength,page=1]{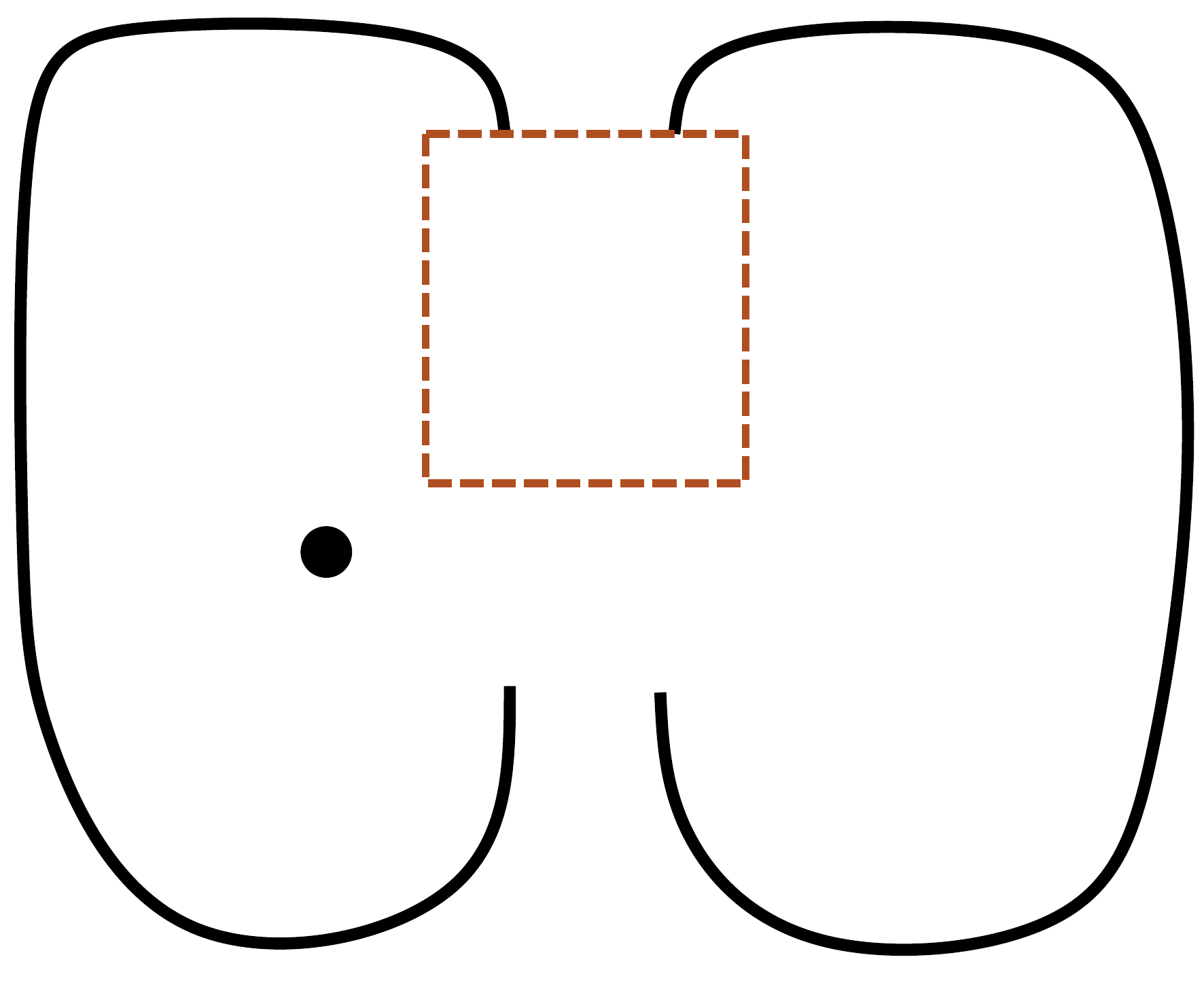}}%
    \put(0.41296056,0.63532405){\color[rgb]{0,0,0}\makebox(0,0)[lt]{\lineheight{1.25}\smash{\begin{tabular}[t]{l}{\footnotesize$n$}\end{tabular}}}}%
    \put(0.39320602,0.5747006){\color[rgb]{0,0,0}\makebox(0,0)[lt]{\lineheight{1.25}\smash{\begin{tabular}[t]{l}{\footnotesize half}\end{tabular}}}}%
    \put(0.3773229,0.51540365){\color[rgb]{0,0,0}\makebox(0,0)[lt]{\lineheight{1.25}\smash{\begin{tabular}[t]{l}{\footnotesize twists}\end{tabular}}}}%
    \put(0.11715,0.40699188){\color[rgb]{0,0,0}\makebox(0,0)[lt]{\lineheight{1.25}\smash{\begin{tabular}[t]{l}{\footnotesize $\tau$}\end{tabular}}}}%
    \put(0,0){\includegraphics[width=\unitlength,page=2]{looping_ringed_n_2_torus_knot_tunnel.pdf}}%
  \end{picture}%
\endgroup%

\caption{Looping of $(\sphere,L_n\cup\tau)$.}
\label{fig:looping_ringed_n_2_torus_knot_tunnel}
\end{subfigure}
\caption{Handlebody-knot exteriors that contain a type $3$-$2$i
annulus.}
\end{figure} 
Next, we consider the union of the $2$-component link $(\sphere,L_n)$ with $n$ odd and the tunnel $\tau$ in Fig.\ \ref{fig:ringed_n_2_torus_knot_tunnel}. 
Then the looping of $\pairLntau$ in Fig.\ \ref{fig:looping_ringed_n_2_torus_knot_tunnel}
induces a handlebody-knot $\pairn$
whose exterior contains a type $3$-$2$i annulus
given by the cabling annulus of the $(n,2)$-torus knot 
component of $(\sphere,L_n)$, 
so we have the following.
\begin{corollary}
The annulus diagram of 
the handlebody-knot $\pairn$ 
obtained by the looping of $\pairLntau$ in Fig.\ \ref{fig:looping_ringed_n_2_torus_knot_tunnel}
is \raisebox{-.2cm}{\includegraphics[scale=.14]{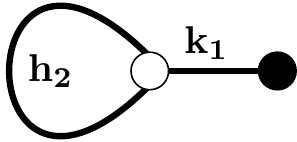}}.
\end{corollary}

Lastly, to produce handlebody-knots
with the annulus diagram 
\raisebox{-.2cm}{\includegraphics[scale=.12]{typetwotwo_ann_1}}, we observe that,
given a handlebody-knot $\pair$
with a type $2$-$2$ annulus $A\subset\ComplHK$, 
the loops $l_+,l_-$ bounds two disks $D_+,D_-$ in $\HKA$, respectively, and $D_+,D_-$ determine 
a spine $\SG_A$ of $\HKA$;
denote by $l_1,l_2$ the constituent loops in $\SG_A$
with $l_2$ disjoint from $D_+$ in $\HKA$ bounded by $l_+$, and orient $l_1,l_2$.
Then we have the following criterion
for the non-uniqueness of $A\subset\ComplHK$.
\begin{lemma}\label{lm:criteria_nonuniqueness_typetwotwo}
\begin{enumerate}[label=\textnormal{(\arabic*)}]
\item \label{itm:typethreetwo}Suppose $\ComplHK$ contains a type $3$-$2$ annulus $A'$.
Then $\lk{l_1}{l_2}\neq \pm 1$.
\item \label{itm:typetwotwo} Suppose $\ComplHK$ contains a type $2$-$2$ annulus $A'$
not isotopic to $A$, and $(\sphere, l_1)$ is a trivial knot. 
Then $(\sphere, l_1\cup l_2)$ is either a trivial link or a Hopf link.
\end{enumerate} 
\end{lemma}
\begin{proof}
\ref{itm:typethreetwo}: \textbf{Case $1$: $A'$ is of type $3$-$2$i.}
Let $W\subset\ComplHK$ be the solid torus 
cut off by $A'$, and $l_w$ an oriented core of $W$.
Note that the core of $A'$ is a $(p,q)$-curve 
on $\partial W$ with
$\vert q\vert>1$ since $A'\subset\ComplHK$ is essential.
If the linking number 
$\lk{l_1}{l_w}$ is $m$, the linking number 
$\lk{l_1}{l_2}$ is $\pm qm\neq \pm 1$. 

\textbf{Case $2$: $A'$ is of type $3$-$2$ii.}
Let $D\subset \HK_A$ be a non-separating disk dual to $l_1$, 
and denote by $V$ the solid torus $\HK_A-\openrnbhd{D}$.
$A'$ cuts $\Compl V$ into two solid tori, one of which, denoted by 
$W$, contains $D$. Note that the core of the annulus 
$W\cap V$ has a slope of $\frac{p}{q}$, $\vert p\vert > 1$,
with respect to $(\sphere,l_2)$.
Let $D_w$ be an oriented meridian disk of $W$.
If the linking number $\lk{l_1}{\partial D_w}=n$,
then the linking number $\lk{l_1}{l_2}=\pm np\neq \pm 1$.

\ref{itm:typetwotwo}: Observe first
that $(\sphere,l_2)$ is trivial by the existence 
of $A'$. Therefore, $(\sphere,l_1\cup l_2)$ is trivial if it is split. Suppose it is non-split.
Then there exists an essential disk $D\subset \Compl{l_2}$ meeting $l_1$ at exactly one point.
Denote by $W$ the $3$-ball $\overline{\Compl{l_2}-\rnbhd{D}}$. Then since $(\sphere, l_1)$ 
is trivial, the ball-arc pair $(W,l_1\cap W)$ is trivial, so $(\sphere,l_1\cup l_2)$ is a Hopf link. 
\end{proof}

\begin{figure}[h]
\begin{subfigure}{0.47\textwidth}
\centering
\def\svgwidth{.6\columnwidth}
\begingroup%
  \makeatletter%
  \providecommand\color[2][]{%
    \errmessage{(Inkscape) Color is used for the text in Inkscape, but the package 'color.sty' is not loaded}%
    \renewcommand\color[2][]{}%
  }%
  \providecommand\transparent[1]{%
    \errmessage{(Inkscape) Transparency is used (non-zero) for the text in Inkscape, but the package 'transparent.sty' is not loaded}%
    \renewcommand\transparent[1]{}%
  }%
  \providecommand\rotatebox[2]{#2}%
  \newcommand*\fsize{\dimexpr\f@size pt\relax}%
  \newcommand*\lineheight[1]{\fontsize{\fsize}{#1\fsize}\selectfont}%
  \ifx\svgwidth\undefined%
    \setlength{\unitlength}{850.39370079bp}%
    \ifx\svgscale\undefined%
      \relax%
    \else%
      \setlength{\unitlength}{\unitlength * \real{\svgscale}}%
    \fi%
  \else%
    \setlength{\unitlength}{\svgwidth}%
  \fi%
  \global\let\svgwidth\undefined%
  \global\let\svgscale\undefined%
  \makeatother%
  \begin{picture}(1,0.83333333)%
    \lineheight{1}%
    \setlength\tabcolsep{0pt}%
    \put(0,0){\includegraphics[width=\unitlength,page=1]{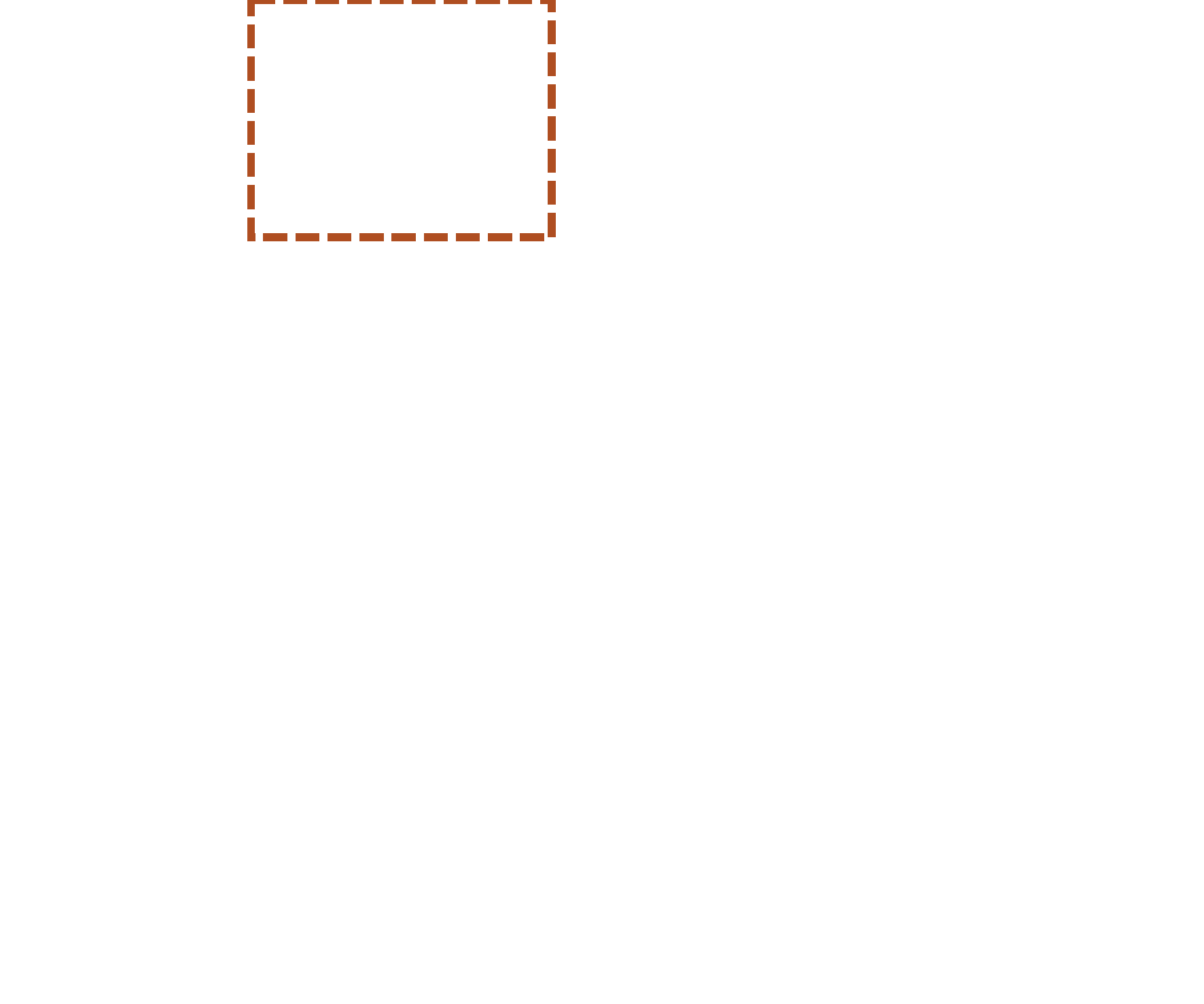}}%
    \put(0.25517252,0.78775804){\color[rgb]{0,0,0}\makebox(0,0)[lt]{\lineheight{1.25}\smash{\begin{tabular}[t]{l}{\footnotesize$n$}\end{tabular}}}}%
    \put(0.23986958,0.72281376){\color[rgb]{0,0,0}\makebox(0,0)[lt]{\lineheight{1.25}\smash{\begin{tabular}[t]{l}{\footnotesize half}\end{tabular}}}}%
    \put(0.22618527,0.65813744){\color[rgb]{0,0,0}\makebox(0,0)[lt]{\lineheight{1.25}\smash{\begin{tabular}[t]{l}{\footnotesize twists}\end{tabular}}}}%
    \put(0.86471466,0.50123159){\color[rgb]{0,0,0}\makebox(0,0)[lt]{\lineheight{1.25}\smash{\begin{tabular}[t]{l}{\footnotesize $\tau$}\end{tabular}}}}%
    \put(0,0){\includegraphics[width=\unitlength,page=2]{link_w_lk_pm1_tunnel.pdf}}%
  \end{picture}%
\endgroup%

\caption{Union of $(\sphere,L_n)$, $n$ even, and a tunnel $\tau$.}
\label{fig:link_w_lk_pm1_tunnel}
\end{subfigure}
\raisebox{.4cm}{$\Longrightarrow$}
\begin{subfigure}{0.47\textwidth}
\centering
\def\svgwidth{.6 \columnwidth}
\begingroup%
  \makeatletter%
  \providecommand\color[2][]{%
    \errmessage{(Inkscape) Color is used for the text in Inkscape, but the package 'color.sty' is not loaded}%
    \renewcommand\color[2][]{}%
  }%
  \providecommand\transparent[1]{%
    \errmessage{(Inkscape) Transparency is used (non-zero) for the text in Inkscape, but the package 'transparent.sty' is not loaded}%
    \renewcommand\transparent[1]{}%
  }%
  \providecommand\rotatebox[2]{#2}%
  \newcommand*\fsize{\dimexpr\f@size pt\relax}%
  \newcommand*\lineheight[1]{\fontsize{\fsize}{#1\fsize}\selectfont}%
  \ifx\svgwidth\undefined%
    \setlength{\unitlength}{850.39370079bp}%
    \ifx\svgscale\undefined%
      \relax%
    \else%
      \setlength{\unitlength}{\unitlength * \real{\svgscale}}%
    \fi%
  \else%
    \setlength{\unitlength}{\svgwidth}%
  \fi%
  \global\let\svgwidth\undefined%
  \global\let\svgscale\undefined%
  \makeatother%
  \begin{picture}(1,0.83333333)%
    \lineheight{1}%
    \setlength\tabcolsep{0pt}%
    \put(0,0){\includegraphics[width=\unitlength,page=1]{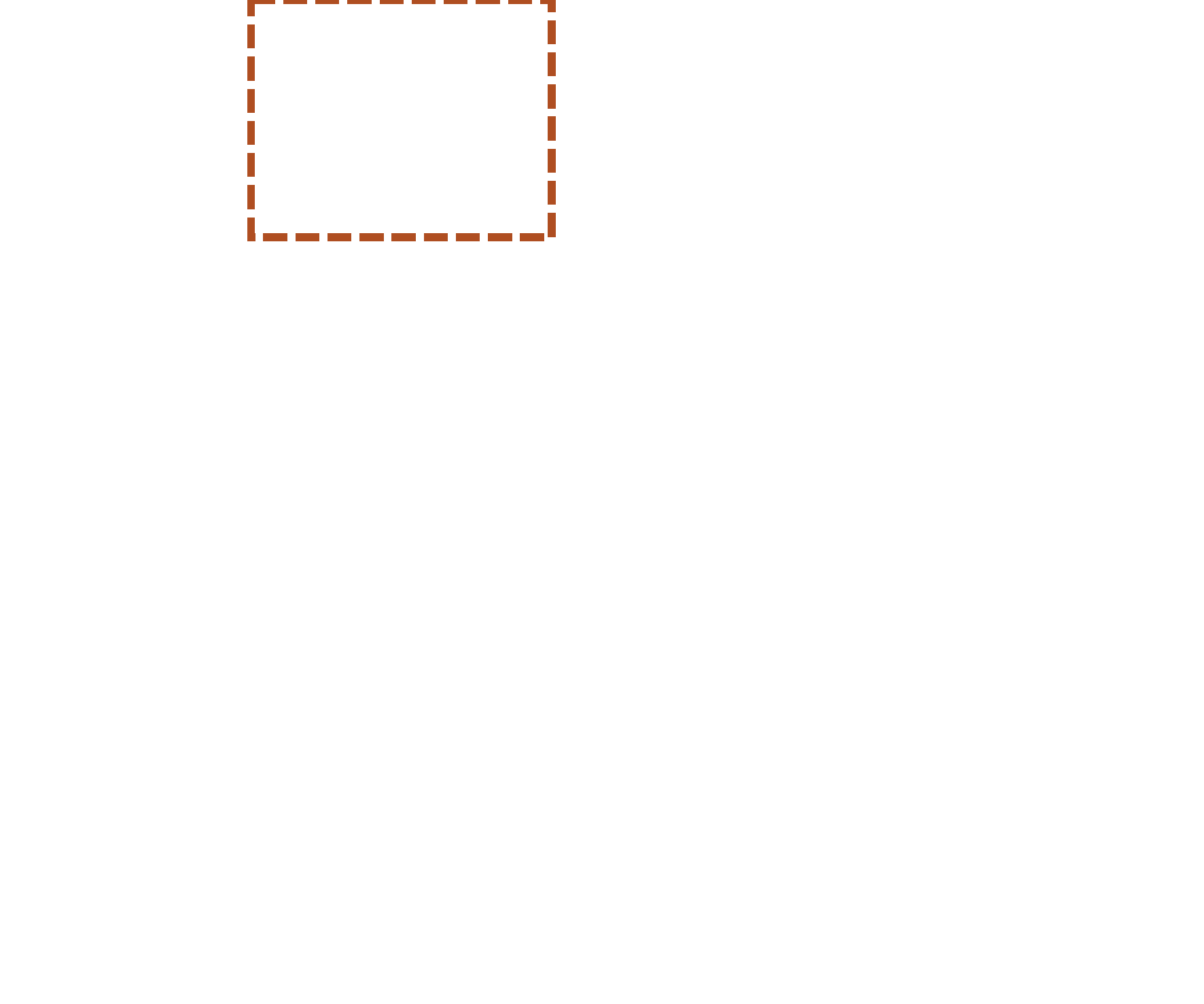}}%
    \put(0.25517252,0.78775804){\color[rgb]{0,0,0}\makebox(0,0)[lt]{\lineheight{1.25}\smash{\begin{tabular}[t]{l}{\footnotesize$n$}\end{tabular}}}}%
    \put(0.23986958,0.72281376){\color[rgb]{0,0,0}\makebox(0,0)[lt]{\lineheight{1.25}\smash{\begin{tabular}[t]{l}{\footnotesize half}\end{tabular}}}}%
    \put(0.22618527,0.65813744){\color[rgb]{0,0,0}\makebox(0,0)[lt]{\lineheight{1.25}\smash{\begin{tabular}[t]{l}{\footnotesize twists}\end{tabular}}}}%
    \put(0,0){\includegraphics[width=\unitlength,page=2]{looping_link_w_lk_pm1_tunnel.pdf}}%
  \end{picture}%
\endgroup%
 
\caption{Looping of $(\sphere,L_n\cup\tau)$.}
\label{fig:looping_link_lk_pm1_tunnel}
\end{subfigure}
\caption{Handlebody-knots with a unique
type $2$ annulus.}
\end{figure}
Consider now the handcuff graph given by 
the union of the $2$-component link $\pairLn$ with $n$ even and
the tunnel $\tau$ in Fig.\ \ref{fig:link_w_lk_pm1_tunnel}.
\begin{corollary}
The handlebody-knot induced by
the looping of $\pairLntau$ in Fig.\ \ref{fig:looping_link_lk_pm1_tunnel}
with even $n\neq 0$
is atoroidal with the annulus diagram
\raisebox{-.2cm}{\includegraphics[scale=.12]{typetwotwo_ann_1}}.
\end{corollary}
\begin{proof}
It follows from Lemmas \ref{lm:looping_nonsplit_link_tunnel} and 
\ref{lm:criteria_nonuniqueness_typetwotwo}
since the linking number of 
$\pairLn$ is $\pm 1$, and it is not a Hopf link, for every even $n\neq 0$.
\end{proof}

Handlebody-knots induced by Figs.\ \ref{fig:looping_hk5_2}, \ref{fig:double_looping_hk5_2}, 
\ref{fig:looping_n_2_torus_link}, \ref{fig:looping_ringed_n_2_torus_knot_tunnel}, 
and \ref{fig:looping_link_lk_pm1_tunnel} imply the following.
\begin{proposition}\label{prop:realization}
Annulus diagrams in Theorem \ref{teo:classification_char_diagram_type_two} can
all be realized.
\end{proposition}




\section*{Acknowledgment}
The author thanks Makoto Sakuma and Yuya Koda 
for the helpful and constructive discussions.
The work was supported by National Sun Yat-sen University and Academia Sinica, and MoST (grant no. 110-2115-M-001-004-MY3), Taiwan. 


\end{document}